\documentclass[11pt]{amsart}
\usepackage[all,cmtip]{xy}
\usepackage{amssymb,amsmath,amsgen,amsxtra,amsfonts} 
\usepackage{times}

\begin{document}
%
%
%
\theoremstyle{definition}
\newtheorem{Definition}{Definition}[section]
\newtheorem*{Definitionx}{Definition}
\newtheorem{Conclusion}{Conclusion}[section]
\newtheorem{Construction}{Construction}[section]
\newtheorem{Example}[Definition]{Example}
\newtheorem{Examples}[Definition]{Examples}
\newtheorem{Exercise}[Definition]{Exercise}
\newtheorem{Exercises}[Definition]{Exercises}
\newtheorem{Remark}[Definition]{Remark}
\newtheorem{Remarks}[Definition]{Remarks}
\newtheorem{Caution}[Definition]{Caution}
\newtheorem{Conjecture}[Definition]{Conjecture}
\newtheorem{Question}[Definition]{Question}
\newtheorem{Questions}[Definition]{Questions}
\newtheorem*{Acknowledgements}{Acknowledgements}
\theoremstyle{plain}
\newtheorem{Theorem}[Definition]{Theorem}
\newtheorem*{Theoremx}{Theorem}
\newtheorem{Proposition}[Definition]{Proposition}
\newtheorem*{Propositionx}{Proposition}
\newtheorem{Lemma}[Definition]{Lemma}
\newtheorem{Corollary}[Definition]{Corollary}
\newtheorem{Fact}[Definition]{Fact}
\newtheorem{Facts}[Definition]{Facts}
\newtheoremstyle{voiditstyle}{3pt}{3pt}{\itshape}{\parindent}%
{\bfseries}{.}{ }{\thmnote{#3}}%
\theoremstyle{voiditstyle}
\newtheorem*{VoidItalic}{}
\newtheoremstyle{voidromstyle}{3pt}{3pt}{\rm}{\parindent}%
{\bfseries}{.}{ }{\thmnote{#3}}%
\theoremstyle{voidromstyle}
\newtheorem*{VoidRoman}{}

%
\newcommand{\prf}{\par\noindent{\sc Proof.}\quad}
\newcommand{\blowup}{\rule[-3mm]{0mm}{0mm}}
\newcommand{\cal}{\mathcal}
\newcommand{\mathds}[1]{{\mathbb #1}}
\newcommand{\Aff}{{\mathds{A}}}
\newcommand{\BB}{{\mathds{B}}}
\newcommand{\CC}{{\mathds{C}}}
\newcommand{\FF}{{\mathds{F}}}
\newcommand{\GG}{{\mathds{G}}}
\newcommand{\HH}{{\mathds{H}}}
\newcommand{\NN}{{\mathds{N}}}
\newcommand{\ZZ}{{\mathds{Z}}}
\newcommand{\PP}{{\mathds{P}}}
\newcommand{\QQ}{{\mathds{Q}}}
\newcommand{\RR}{{\mathds{R}}}
\newcommand{\Liea}{{\mathfrak a}}
\newcommand{\Lieb}{{\mathfrak b}}
\newcommand{\Lieg}{{\mathfrak g}}
\newcommand{\Liem}{{\mathfrak m}}
\newcommand{\ideala}{{\mathfrak a}}
\newcommand{\idealb}{{\mathfrak b}}
\newcommand{\idealg}{{\mathfrak g}}
\newcommand{\idealm}{{\mathfrak m}}
\newcommand{\idealp}{{\mathfrak p}}
\newcommand{\idealq}{{\mathfrak q}}
\newcommand{\idealI}{{\cal I}}
\newcommand{\lin}{\sim}
\newcommand{\num}{\equiv}
\newcommand{\dual}{\ast}
\newcommand{\iso}{\cong}
\newcommand{\homeo}{\approx}
\newcommand{\mm}{{\mathfrak m}}
\newcommand{\pp}{{\mathfrak p}}
\newcommand{\qq}{{\mathfrak q}}
\newcommand{\rr}{{\mathfrak r}}
\newcommand{\pP}{{\mathfrak P}}
\newcommand{\qQ}{{\mathfrak Q}}
\newcommand{\rR}{{\mathfrak R}}
%
%
\newcommand{\dq}{{``}}
\newcommand{\OO}{{\cal O}}
\newcommand{\into}{{\hookrightarrow}}
\newcommand{\onto}{{\twoheadrightarrow}}
\newcommand{\Spec}{{\rm Spec}\:}
\newcommand{\BigSpec}{{\rm\bf Spec}\:}
\newcommand{\Proj}{{\rm Proj}\:}
\newcommand{\Pic}{{\rm Pic }}
\newcommand{\Br}{{\rm Br}}
\newcommand{\NS}{{\rm NS}}
\newcommand{\chit}{\chi_{\rm top}}
\newcommand{\KanDiv}{{\cal K}}
\newcommand{\perdef}{{\stackrel{{\rm def}}{=}}}
\newcommand{\Cycl}[1]{{\ZZ/{#1}\ZZ}}
\newcommand{\Sym}{{\mathfrak S}}
\newcommand{\Xcan}{X_{{\rm can}}}
\newcommand{\Ycan}{Y_{{\rm can}}}
\newcommand{\ab}{{\rm ab}}
\newcommand{\Aut}{{\rm Aut}}
\newcommand{\Hom}{{\rm Hom}}
\newcommand{\Supp}{{\rm Supp}}
\newcommand{\ord}{{\rm ord}}
\newcommand{\divisor}{{\rm div}}
\newcommand{\Alb}{{\rm Alb}}
\newcommand{\Jac}{{\rm Jac}}
\newcommand{\Ig}{{\rm Ig}}
\newcommand{\piet}{{\pi_1^{\rm \acute{e}t}}}
\newcommand{\Het}[1]{{H_{\rm \acute{e}t}^{{#1}}}}
\newcommand{\Hcris}[1]{{H_{\rm cris}^{{#1}}}}
\newcommand{\HdR}[1]{{H_{\rm dR}^{{#1}}}}
\newcommand{\HdRalg}[1]{{H_{\rm dR, alg}^{{#1}}}}
\newcommand{\hdR}[1]{{h_{\rm dR}^{{#1}}}}
\newcommand{\defin}[1]{{\bf #1}}

\title[Supersingular K3 Surfaces]{Lectures on Supersingular K3 Surfaces and the Crystalline Torelli Theorem}


\author{Christian Liedtke}
\dedicatory{Expanded lecture notes of a course given as part of the IRMA master class ``Around Torelli's theorem for K3 surfaces'' held in Strasbourg, October 28 - November 1, 2013}
\address{Zentrum Mathematik - M11, Boltzmannstr. 3, D-85748 Garching bei M\"unchen, Germany}
\date{May 19, 2015}
\email{liedtke@ma.tum.de}

\begin{abstract}
 We survey crystalline cohomology, crystals,
 and formal group laws with an emphasis on geometry.
 We apply these concepts to K3 surfaces, and
 especially to supersingular K3 surfaces.
 In particular, we discuss stratifications of the moduli space
 of polarized K3 surfaces in positive characteristic, 
 Ogus' crystalline Torelli theorem for supersingular K3 surfaces,
 the Tate conjecture, and the unirationality of K3 surfaces. 
\end{abstract}

\maketitle


\section*{Introduction}

In these notes, we cover the following topics
\begin{enumerate}
  \item[--] algebraic de Rham cohomology, crystalline cohomology, and $F$-crystals, 
  \item[--] characteristic-$p$ aspects of K3 surfaces,
  \item[--] Ogus' crystalline Torelli theorem for supersingular K3 surfaces,
  \item[--] formal group laws, and in particular, the formal Brauer group, and
  \item[--] unirationality and supersingularity of K3 surfaces.
\end{enumerate}
We assume familiarity with algebraic geometry, say, at the level of 
the textbooks of Hartshorne \cite{Hartshorne} and Griffiths--Harris \cite{Griffiths Harris}.
\medskip

One aim of these notes is to convince the reader that crystals and crystalline
cohomology are rather explicit objects, and that they are characteristic-$p$
versions of Hodge structures and de Rham cohomology, respectively.
Oversimplifying and putting it a little bit sloppily, the {\em crystalline cohomology} 
of a smooth and proper variety in characteristic $p$ 
is the de Rham cohomology of a lift to characteristic zero.
(Unfortunately, such lifts may not exist, and even if they do, they may not be unique -
it was Grothendieck's insight that something like crystalline cohomology 
exists nevertheless and that it is well-defined.)
Just as complex de Rham cohomology comes with complex conjugation, crystalline
cohomology comes with a Frobenius action, and this latter leads to the notion
of a {\em crystal}.
Therefore, {\em period maps} in characteristic $p$ should take 
values in moduli spaces of crystals.
For example, {\em Ogus' crystalline Torelli theorem} states that moduli spaces of certain 
K3 surfaces, namely, {\em supersingular K3 surfaces}, can be entirely understood 
via a period map to a moduli space of suitably enriched crystals.
Conversely, the classification 
of crystals arising from K3 surfaces gives
rise to a stratification of the moduli space of K3 surfaces in characteristic $p$: for example,
the {\em height stratification} in terms of {\em Newton polygons} arises this way.

A second aim of these notes is to introduce {\em formal group laws}, and,
following Artin and Mazur, to explain how they arise
naturally from algebraic varieties.
For us, the most important examples will be the {\em formal Picard group} and the {\em formal Brauer group}
of a smooth and proper variety in characteristic $p$.
Whereas the former arises as formal completion of the Picard scheme along its origin, 
the latter does not have such a description, and is something new.
Oversimplifying again, we have a sort of crystal associated to a formal group law, namely, its {\em Cartier--Dieudonn\'e module}.
For example, the Cartier--Dieudonn\'e modules of the formal Picard group and the formal
Brauer group give rise to subcrystals inside first and second crystalline cohomology of the variety in question.
And despite their abstract appearance, these formal group laws do have a geometric
interpretation: for example, for supersingular K3 surfaces, the formal Brauer group 
controls certain very special one-parameter deformations, {\em moving torsors},
which are characteristic-$p$ versions of twistor space.
These deformations are the key to the proof that 
{\em supersingular K3 surfaces are unirational}.\footnote{Due to a mistake in the proof, the unirationality of supersingular K3 surfaces remains a conjecture}

And finally, a third aim of these notes is to make some of the more abstract concepts
more accessible, which is why we have put an emphasis on computing everything for K3 surfaces.
Even for a reader who cares not so much for K3 surfaces, these notes may be interesting,
since we show by example, how to perform computations with crystals and formal Brauer groups. 
\medskip

These notes are organized as follows:
\begin{trivlist}
 \item[{\em Section \ref{sec:cohomology}}] We start by discussing de Rham cohomology
   over the complex numbers, and then turn to algebraic de Rham cohomology.
   After a short detour to $\ell$-adic cohomology, we introduce the Witt ring $W$, and survey
   crystalline cohomology.
 \item[{\em Section \ref{sec:K3}}] We define K3 surfaces,
   give examples, and discuss their position within the surface classification.
   Then, we compute their cohomological invariants, and end by introducing
   polarized moduli spaces.
  \item[{\em Section \ref{sec:fcrystals}}] 
   Crystalline cohomology takes it values in $W$-modules, where $W$ denotes the Witt
   ring, and it comes with a Frobenius-action, which leads to the notion of an $F$-crystal.
   After discussing the Dieudonn\'e--Manin classification of $F$-crystals up to isogeny,
   we introduce Hodge and Newton polygons of $F$-crystals.
  \item[{\em Section \ref{sec:supersingularK3}}] 
   The $F$-crystal associated to the second crystalline cohomology group 
   of a K3 surface comes with a quadratic form arising from Poincar\'e duality, 
   which is captured in the notion of a K3 crystal.
   After discussing supersingular K3 crystals and the Tate conjecture, we explicitly
   classify these crystals and construct their moduli space.
  \item[{\em Section \ref{sec:torelli}}]
   Associated to a supersingular K3 surface, we have its supersingular K3 crystal,
   which gives rise to a period map from the moduli space of supersingular K3 surfaces
   to the moduli space of supersingular K3 crystals.
   Equipping supersingular K3 crystals with ample cones, we obtain a new period
   map, which is an isomorphism by Ogus' crystalline Torelli theorem.
  \item[{\em Section \ref{sec:formal groups}}]
   After defining formal group laws and giving examples, we introduce basic
   invariants and state the Cartier--Dieudonn\'e classification.
   Then, we discuss formal group laws arising from algebraic varieties, and relate
   their Cartier--Dieudonn\'e modules to crystalline cohomology.
  \item[{\em Section \ref{sec:unirational}}]
   We show that a K3 surface is unirational if and only if it is supersingular.
   The key idea is to give a geometric construction that uses the formal Brauer group.
  \item[{\em Section \ref{sec:stratification}}]
   Finally, we discuss several stratifications of the moduli space of polarized
   K3 surfaces, which arise from $F$-crystals and formal Brauer groups.
\end{trivlist}

For further reading, we refer the interested reader
\begin{itemize}
 \item[--] to Wedhorn's notes \cite{Wedhorn} for more on de Rham cohomology
    and $F$-zips, as well as to
    Illusie's survey \cite{Illusie Frobenius} for the degeneration
    of the Fr\"olicher spectral sequence with characteristic-$p$ methods,
  \item[--] to Chambert-Loir's survey \cite{CL} for more on crystalline cohomology, 
    and to Katz's article \cite{Katz crystal} for more on $F$-crystals,
  \item[--] to my own survey \cite{Liedtke overview} for more on 
    the classification of surfaces in positive characteristic, 
  \item[--] to Huybrechts' lecture notes \cite{Huybrechts} for more on K3 surfaces,
  \item[--]  to Ogus' original articles \cite{Ogus} and \cite{Ogus PM} for more on 
    supersingular K3 crystals and the Torelli theorem, and finally,
  \item[--] to the books of  Hazewinkel \cite{Hazewinkel} and Zink \cite{Zink} for 
    more on formal group laws and their applications.
\end{itemize}
Of course, this list is by no means complete, but merely represents 
a small part of the literature for some of the topics touched upon in these notes.
\medskip

\begin{Acknowledgements}
These notes grew out of a lecture series I gave as part of the master class
``Around Torelli's theorem for K3 surfaces''
at the Institut de Recherche Math\'ematique Avanc\'ee (IRMA)
in Strasbourg from October 28 to November 1, 2013, which was
organized by Christian Lehn, Gianluca Pacienza, and Pierre Py.
I thank the organizers for the invitation and hospitality.
It was a great pleasure visiting the IRMA and giving these lectures.
I thank Nicolas Addington, Daniel Bragg, Gerard van der Geer, Christopher Hacon, Eike Lau,
Max Lieblich, and the referee  for many helpful comments on earlier versions of these notes.
\end{Acknowledgements}

\section{Crystalline Cohomology}
\label{sec:cohomology}

\subsection{Complex Geometry}
\label{subsec:complex geometry}
Let $X$ be a smooth and projective variety over the
field $\CC$ of complex numbers.
In this subsection we briefly recall how the topological, differentiable and holomorphic
structure on $X$ give rise to extra structure on the 
de~Rham cohomology groups $\HdR{\ast}(X/\CC)$, and how these structures
are related.

First, we consider $X$ as a differentiable manifold, where differentiable
shall always mean differentiable of class $C^\infty$ without further mentioning it.
Let ${\cal A}_X^n$ be the sheaf of $\CC$-valued differentiable $n$-forms
with respect to the analytic topology, and in particular, 
${\cal A}_X^0$ is the sheaf of $\CC$-valued differentiable
functions on $X$.
Then, by the Poincar\'e lemma, the {\em differentiable de~Rham complex}
$$
  0\,\to\,{\cal A}_X^0\,\stackrel{d}{\longrightarrow}
  \,{\cal A}_X^1\,\stackrel{d}{\longrightarrow}\,{\cal A}^2_X\,\stackrel{d}{\longrightarrow}\,...
$$
is a fine resolution of the constant sheaf $\underline{\CC}$.
Taking global sections, we obtain a complex, whose cohomology is  by definition 
de~Rham cohomology.
And since the differentiable de~Rham complex is a fine resolution of $\underline{\CC}$,
we obtain a natural isomorphism 
$$
    H^n(X,\underline{\CC})\,\iso\,\HdR{n}(X/\CC)
$$
for all $n$,  and refer to \cite[Chapter 5]{Warner} for details and proofs.

Next, we consider $X$ as an algebraic variety with the Zariski topology.
Let $\Omega_{X/\CC}^n$ be the coherent $\OO_X$-module of K\"ahler $n$-forms,
and in particular, $\Omega_{X/\CC}^0$ is the structure sheaf $\OO_X$.
As in the differentiable setting, we obtain a complex of sheaves
$$
 0\,\to\,\OO_X\,\stackrel{d}{\longrightarrow}\,
 \Omega_{X/\CC}^1\,\stackrel{d}{\longrightarrow}\,
 \Omega_{X/\CC}^2\,\stackrel{d}{\longrightarrow}\,...,
$$
the {\em algebraic de~Rham-complex}.
Now, the $\Omega_{X/\CC}^n$ are locally free $\OO_X$-modules, but not acyclic in general.
Therefore, we choose injective resolutions of the $\Omega_{X/\CC}^n$ 
that are compatible with the differentials of the algebraic de~Rham complex, and eventually,
we obtain a double complex.
Then, we take global sections of the injective sheaves of this double complex and
take cohomology, which defines {\em algebraic de~Rham cohomology}.
In particular, algebraic de~Rham cohomology arises as hypercohomology of the 
algebraic de~Rham complex.
Since the columns of this double complex compute the coherent cohomology
groups $H^j(X,\Omega_{X/\CC}^j)$,
there is a spectral sequence to algebraic de~Rham cohomology
$$
  E_1^{i,j}\,=\,H^j(X,\Omega_{X/\CC}^i)\,\Rightarrow\,\HdRalg{i+j}(X/\CC),
$$
the {\em Fr\"olicher spectral sequence}.
We refer to \cite{Illusie Frobenius} or \cite{Wedhorn} for details.

Of course, we can also consider $X$ as a complex manifold with the analytic topology.
Then, we can define the holomorphic de~Rham complex, which gives
rise to a  holomorphic Fr\"olicher spectral sequence,
see \cite[Chapter 3.5]{Griffiths Harris} or \cite[Chapitre II.8]{Voisin}.
By the holomorphic Poincar\'e lemma, the holomorphic
de\-Rham complex is a resolution (although not acyclic) of the
constant sheaf $\underline{\CC}$, which is why
the hypercohomology of the holomorphic de~Rham complex
is canonically isomorphic to the differentiable de~Rham cohomology,
see loc. cit.
By Serre's GAGA-theorems \cite{Serre GAGA},
the algebraic and holomorphic $E_1^{p,q}$'s are canonically
isomorphic, and thus, also the hypercohomologies of the  
holomorphic and algebraic de~Rham complex are canonically
isomorphic.
In particular, the differentiable, holomorphic, and algebraic
de~Rham cohomologies are mutually and canonically isomorphic,
and there is only one Fr\"olicher spectral sequence.
We will therefore omit the subscript ${\rm alg}$ from algebraic de~Rham cohomology
in the sequel.
Apart from the already given references, we refer the interested reader  
to \cite{Grothendieck deRham} for more on this subject.

For complex projective varieties, and even for compact
K\"ahler manifolds, the Fr\"olicher spectral sequence degenerates at $E_1$
by the Hodge decomposition theorem, see \cite[Chapter 0.7]{Griffiths Harris}
or \cite[Th\'eor\`eme 8.28]{Voisin},
as well as Exercise \ref{exercise: froelicher} below.
Next, the Fr\"olicher spectral sequence gives rise to a filtration
$$
 0\,=\,F^{n+1}\,\subseteq\,...\,\subseteq\,F^i\,\subseteq\,...\,\subseteq\,F^0=\HdR{n}(X/\CC),
$$
the {\em Hodge filtration}.
Since the Fr\"olicher spectral sequence degenerates at $E_1$, 
we obtain canonical isomorphisms 
$$
   F^i/F^{i+1}\,\iso\, E_1^{i,n-i}\,=\,H^{n-i}(X,\Omega_X^i)\mbox{ \quad for all } i,n \mbox{ with } 0\leq i\leq n.
$$

Next, we consider $X$ again only as a differentiable manifold with its analytic topology.
Then, sheaf cohomology of the constant 
sheaf $\underline{\ZZ}$ is isomorphic to singular cohomology
$$
   H^n(X,\underline{\ZZ})\,\iso\,H^n_{\rm sing}(X,\ZZ),
$$
see, for example, \cite[Chapter 5]{Warner}.
In particular, the inclusion $\ZZ\subset\CC$ gives rise to an isomorphism
$$
  H^{n}_{\rm sing}(X,\ZZ)\otimes_\ZZ\CC\,\iso\,\HdR{n}(X/\CC).
$$
and thus, to an integral structure on de~Rham cohomology. 
Similarly, the inclusion $\RR\subset\CC$ gives rise to a real structure
on de~Rham cohomology.
In particular, we obtain a complex conjugation on $\HdR{n}(X/\CC)$,
which we apply to the Hodge filtration to define a second filtration,
the so-called {\em conjugate filtration}, by setting
$$
  F_{\rm con}^i \,:=\, \overline{F^{n-i}}\,.
$$
It follows from the Hodge decomposition theorem that the Hodge and its conjugate filtration satisfy
$$
F^i \,\cap\, \overline{F^{n-i}}\,\iso\,
E_1^{i,n-i}\,=\,H^{n-i}(X,\Omega^i_{X/\CC}),
$$
see \cite[Chapter 0.6]{Griffiths Harris} or \cite[Remarque 8.29]{Voisin}.
From this, we deduce 
$$
    H^{n-i}(X,\Omega_{X/\CC}^i) \,=\, E_1^{i,n-i} \,\iso\, \overline{E_1^{n-i,i}}
    \,=\, \overline{H^i(X,\Omega_{X/\CC}^{n-i})},
$$
for all $i,n$ with $0\leq i\leq n$.
We stress that this isomorphism is induced from complex conjugation,
making it a priori not algebraic.

Putting all the results and observations made so far together, 
we see that each de~Rham cohomology group $\HdR{n}(X/\CC)$ 
is a finite dimensional $\CC$-vector space together with the following data
\begin{enumerate}
  \item an integral structure, 
  \item a real structure, and in particular, a complex conjugation $^{-}$, and
  \item two filtrations $F^\bullet$ and $F_{\rm con}^\bullet$.
\end{enumerate}
Again, we refer to \cite{Illusie Frobenius} and \cite{Wedhorn} for details and further references.

\subsection{Algebraic de~Rham cohomology}
\label{subsec:algebraic deRham}
Let $X$ be a smooth and proper variety, but now over a field $k$ of arbitrary characteristic.
The definition of algebraic de\-Rham cohomology and
the Fr\"olicher spectral sequence (but not its degeneration at $E_1$) in the previous section 
is purely algebraic, and thus, we have these cohomology groups and spectral sequences
also over $k$.
Let us now discuss $E_1$-degeneration of the Fr\"olicher spectral sequence, as well as
extra structures on de~Rham cohomology in this purely algebraic setting.
Although many aspects of this section are discussed in greater detail
in \cite{Wedhorn}, let us run through the main points needed later on
for the reader's convenience.

\begin{Exercise}
  \label{exercise: froelicher}
  Let $X$ be a smooth and proper variety over a field $k$.
  Show that already the existence of the Fr\"olicher spectral sequence implies
  the inequalities
  $$
      \sum_{i+j=n} h^{i,j}(X)\,\geq\, h^n_{\rm dR}(X)\mbox{ \quad for all }n\geq1,
  $$
  where $h^{i,j}(X)=\dim_k H^j(X,\Omega_{X/k}^i)$ and 
  $h^n_{\rm dR}(X)=\dim_k \HdR{n}(X/k)$.
  Moreover, show that equality for all $n$ is equivalent to the degeneration
  of the Fr\"olicher spectral sequence at $E_1$. 
\end{Exercise}

If $k$ is of characteristic zero, then the Fr\"olicher spectral sequence
of $X$ degenerates at $E_1$ by the following line of reasoning,
which is an instance of the so-called {\em Lefschetz principle}:
namely, $X$ can be defined over a subfield
$k_0\subseteq k$ that is finitely generated over $\QQ$, and then,
$k_0$ can be embedded into $\CC$.
Since cohomology does not change under flat base change,
and field extensions are flat, it suffices to prove 
$E_1$-degeneration for $k=\CC$, where it holds by 
the results of Section \ref{subsec:complex geometry}.

In arbitrary characteristic, $E_1$-degeneration of the  Fr\"olicher spectral sequence 
holds, for example, for curves, Abelian varieties, K3 surfaces, and
complete intersections, see \cite[(1.5)]{Wedhorn}.
For example, for curves it follows from Theorem \ref{thm: Deligne Illusie} below,
and for K3 surfaces, we will show it in 
Proposition \ref{prop: k3 frolicher degenerates}.

On the other hand, Mumford \cite{Mumford Pathologies} gave
explicit examples of smooth and projective surfaces 
$X$ in positive characteristic $p$ with non-closed global $1$-forms.
This means that the exterior derivative 
$d:H^0(X,\Omega_{X/k}^1)\to H^0(X,\Omega_{X/k}^2)$ is non-zero,
which gives rise to a non-zero differential in the Fr\"olicher spectral sequence,
which implies that it does {\em not} degenerate at $E_1$.

In Section \ref{subsec:Witt} below we will recall and construct the ring 
$W(k)$ of Witt vectors associated to a perfect field $k$
of positive characteristic, which is a local, complete, and discrete valuation of 
ring of characteristic zero with residue field $k$, whose unique maximal
ideal is generated by $p$.
In particular, the truncated Witt ring
$W_n(k):=W(k)/(p^n)$ is a flat $\ZZ/p^n\ZZ$-algebra.
For a scheme $X$ over $k$, a {\em lift} of $X$ to $W_n(k)$
is a flat scheme ${\cal X}\to\Spec W_n(k)$ such that
${\cal X}\times_{\Spec W_n(k)}\Spec k\iso X$.
Such lifts may not exist.
However, if they exist, we have the following fundamental result
concerning the degeneration of the Fr\"olicher spectral sequence:

\begin{Theorem}[Deligne--Illusie]
  \label{thm: Deligne Illusie}
  Let $X$ be a smooth and proper variety
  over a perfect field $k$ of characteristic $p\geq\dim(X)$,
  and assume that $X$ admits a lift to $W_2(k)$.
  Then, the Fr\"olicher spectral sequence of $X$ degenerates 
  at $E_1$.
\end{Theorem}

This is the main result of \cite{Deligne Illusie}, and we refer the interested reader
to \cite{Illusie Frobenius} for an expanded version with lots of background information.
In fact, this theorem can be used to obtain a purely algebraic proof
of the $E_1$-degeneration of the Fr\"olicher spectral sequence in
characteristic zero \cite[Theorem 6.9]{Illusie Frobenius}.
This theorem also shows that the above mentioned examples of Mumford
of surfaces with non-closed $1$-forms do not lift to $W_2(k)$.
We refer the interested reader to 
\cite[Section 11]{Liedtke overview} for more about liftings, 
and liftable, as well as non-liftable varieties.

In the previous section, where all varieties were defined over the complex numbers,
we applied complex conjugation to the Hodge filtration to obtain a new
filtration, the conjugate filtration.
Over fields of positive characteristic, there is no complex conjugation.
However, since algebraic de~Rham theory is the hypercohomology of
the de~Rham complex, there exists a second spectral sequence,
the {\em conjugate spectral sequence} 
(this is only a name in analogy with complex geometry: there is no complex
conjugation in positive characteristic)
$$
   E_2^{i,j}\,:=\,H^i(X,{\cal H}^j(\Omega_{X/k}^\bullet))\,\Rightarrow\,\HdR{i+j}(X/k),
$$
see \cite[Chapter 3.5]{Griffiths Harris} or \cite[Section 1]{Wedhorn}.
If $k=\CC$, and when considering $X$ as a complex manifold and holomorphic
differential forms, it follows from the holomorphic Poincar\'e lemma that
the cohomology sheaves ${\cal H}^j(\Omega_{X/\CC}^i)$ are zero for all
$j\geq1$, and thus, the conjugate spectral sequence induces only a trivial
filtration on de~Rham cohomology.

On the other hand, if $k$ is of positive characteristic, then the cohomology
sheaves ${\cal H}^j(\Omega^i_{X/k})$ are usually non-trivial for $q\geq1$.
More precisely, let $F:X\to X'$ be the $k$-linear Frobenius morphism. 
Then, there exists a canonical isomorphism
$$
C^{-1}\,:\,\Omega_{X'/k}^n\,\stackrel{\iso}{\longrightarrow}\,{\cal H}^n(F_\ast\Omega_{X/k}^\bullet)
$$
for all $n\geq0$, the {\em Cartier isomorphism}, and we refer to \cite[Section 3]{Illusie Frobenius}
or \cite[(1.6)]{Wedhorn} for details, definitions, and further references.
Because of this isomorphism, the conjugate spectral sequence is usually non-trivial
in positive characteristic, and gives rise to an interesting second filtration 
on de~Rham cohomology, called the {\em conjugate Hodge filtration}.
The data of de~Rham-cohomology, the Frobenius action, and the two
Hodge filtrations are captured in the following structure.

\begin{Definition}
  \label{def: f-zip}
  An {\em $F$-zip} over a scheme $S$ of positive
  characteristic $p$ is a tuple 
  $(M,C^\bullet,D_\bullet,\varphi_\bullet)$, where $M$ 
  is a locally free $\OO_S$-module of finite rank, $C^\bullet$ is a descending
  filtration on $M$, $D_\bullet$ is an ascending filtration on $M$,
  and $\varphi_\bullet$ is a family of $\OO_S$-linear isomorphisms
  $$
  \varphi_n\,:\,({\rm gr}_C^n)^{(p)}\stackrel{\iso}{\longrightarrow} {\rm gr}_D^n,
  $$
  where ${\rm gr}$ denotes the graded quotient modules, and $^{(p)}$ denotes 
  Frobenius pullback. 
  The function
  $$
  \begin{array}{ccccc}
   \tau &:& \ZZ &\to& \ZZ_{\geq0}\\
   && n &\mapsto&{\rm rank}_{\OO_S}\,{\rm gr}_C^n
  \end{array}
  $$
  is called the {\em (filtration) type} of the $F$-zip.
 \end{Definition}

The category of $F$-zips over $\FF_p$-schemes 
with only isomorphisms as morphisms forms a smooth
Artin-stack $\cal F$ over $\FF_p$, and $F$-zips of type $\tau$
form an open, closed, and quasi-compact substack ${\cal F}^\tau\subseteq{\cal F}$,
see \cite[Proposition 1.7]{Moonen Wedhorn}.
Despite the similarity with period domains for Hodge structures,
this moduli space $\cal F$ is a rather discrete object,
and more of a mod $p$ reduction of such a period domain.
More precisely, if $k$ is an algebraically closed field of
positive characteristic, then the set of
$k$-rational points of ${\cal F}^\tau$ is finite,
and we refer to \cite[Theorem 4.4]{Moonen Wedhorn} and
\cite[Theorem 3.6]{Wedhorn} for precise statements.

Thus, $F$-zips capture discrete invariants arising from de~Rham cohomology
of smooth and proper varieties in positive characteristic.
For example, an $F$-zip $(M,C^\bullet,D_\bullet,\varphi_\bullet)$ is called 
{\em ordinary} if the filtrations $C^\bullet$ and $D_\bullet$ are in opposition,
that is, if the rank $C^i\cap D_j$ is as small as possible for all
$i,j\in\ZZ$.
We refer to \cite{Moonen Wedhorn} and \cite{Wedhorn} for details,
examples, proofs, and further references.
Moreover, one can also consider $F$-zips with additional structure,
such as orthogonal and symplectic forms, see
\cite[Section 5]{Moonen Wedhorn}, as well as
\cite{PWZ} for further generalizations.

Finally, we would like to mention that Vasiu independently developed in \cite{Vasiu1}, \cite{Vasiu2} 
a theory similar to the $F$-zips of Moonen and Wedhorn \cite{Moonen Wedhorn} 
within the framework of Shimura $F$-crystals modulo $p$ and 
truncated Dieudonn\'e modules.
We refer the interested reader to \cite[Sections 3.2.1 and 3.2.9]{Vasiu1}, 
\cite{Vasiu2}, and \cite[Section 2.22]{Xiao} for details, examples, and applications.

\subsection{$\ell$-adic cohomology}
\label{subsec:l-adic cohomology}
The singular cohomology of a smooth and complex projective variety coincides
with cohomology in the locally constant sheaf $\underline{\ZZ}$ 
(with respect to the analytic topology).
Now,  let $X$ be a smooth variety over a field $k$ of characteristic $p\geq0$.
The Zariski topology on $X$ is too coarse to give a good cohomology
theory for locally constant sheaves.
However, for locally constant sheaves of finite Abelian groups, it turns out that \'etale topology,
which is a purely algebraically defined Grothendieck topology that is finer than the Zariski topology, 
has the desired properties, and we shall equip $X$ with this topology.
It goes without saying that then, being locally constant is then meant with respect
to the \'etale topology.
We refer to \cite{SGA 4} and \cite{Milne} for precise definitions.

Let $\ell$ be a prime number, and then, we define {\em $\ell$-adic cohomology}
$$
 \Het{n}(X,\ZZ_\ell)\,:=\,\varprojlim \Het{n}(X,\ZZ/\ell^m\ZZ).
$$
(When considering locally constant sheaves, \'etale cohomology works
best for finite Abelian groups.
So rather than trying to define something like $\Het{n}(-,\ZZ)$ directly, we take the inverse
limit over $\Het{n}(-,\ZZ/\ell^m\ZZ)$, which then results in coefficients
in the $\ell$-adic numbers $\ZZ_\ell$, see also Section \ref{sec:Serres observation} below.
To deal directly with infinite groups, Bhatt and Scholze \cite{Bhatt} introduced the 
pro-\'etale topology, but we will not pursue this here.)
Next, we define
$$
 \Het{n}(X,\QQ_\ell)\,:=\,\Het{n}(X,\ZZ_\ell)\otimes_{\ZZ_\ell}\QQ_\ell.
$$
If $X$ is smooth and proper over an algebraically closed field $k$, 
then this cohomology theory has the following properties:
\begin{enumerate}
 \item $\Het{n}(X,\QQ_\ell)$ is a contravariant functor in $X$.
   The cohomology groups are finite dimensional $\QQ_\ell$-vector spaces, and
   zero if $n<0$ or $n>2\dim(X)$.
 \item There is a cup-product structure
   $$
     \cup_{i,j}\,:\,\Het{i}(X,\QQ_\ell)\,\times\,\Het{j}(X,\QQ_\ell)\,\to\,\Het{i+j}(X,\QQ_\ell).
   $$
   Moreover, $\Het{2\dim(X)}(X,\QQ_\ell)$ is $1$-dimensional, and 
   $\cup_{n,2\dim(X)-n}$ induces a perfect pairing, called {\em Poincar\'e duality}.
 \item $\Het{n}(X,\ZZ_\ell)$ defines an integral structure on $\Het{n}(X,\QQ_\ell)$,
 \item If $k=\CC$, one can choose an inclusion $\QQ_\ell\subset\CC$
     (such inclusions exist using cardinality arguments and the axiom of choice, 
     but they are neither canonical nor compatible with the topologies on these fields), 
     and then, there exist isomorphisms for all $n$ 
     $$
         \Het{n}(X,\QQ_\ell)\otimes_{\QQ_\ell}\CC\,\iso\,H^n(X,\underline{\CC}), 
     $$
     where we consider $X$ as a differentiable manifold on the right hand side
     and $\underline{\CC}$ is locally constant with respect to the analytic topology.
     This comparison isomorphism connects $\ell$-adic cohomology to
     de Rham-, singular, and constant sheaf cohomology if $k=\CC$.
  \item If ${\rm char}(k)=p>0$ and $\ell\neq p$, then the dimension
     $\dim_{\QQ_\ell}\Het{n}(X,\QQ_\ell)$ is independent of $\ell$
     (see \cite{Katz Messing}).
     Thus,
     $$
        b_n(X)\,:=\,\dim_{\QQ_\ell}\Het{n}(X,\QQ_\ell)
     $$
     is well-defined for $\ell\neq p$ and it is called the 
     {\em $n^{th}$ Betti number}.
  \item Finally, there exists a Lefschetz fixed point formula, 
     there are base change formulas,
     there exist cycle classes in $\Het{2q}(X,\QQ_\ell)$
     for codimension $q$ subvarieties,...
\end{enumerate}

We refer to \cite[Appendix C]{Hartshorne} for an overview, and to
\cite{Milne} or \cite{SGA 4} for a thorough treatment.
The following example shows that the assumption $\ell\neq p$ in property (5) above 
is crucial, and gives a hint of the subtleties involved.

\begin{Example}
  \label{example: wrong p-adic cohomology}
  Let $A$ be a $g$-dimensional Abelian variety over an algebraically
  closed field $k$ of positive characteristic $p$.
  For a prime $\ell$, we define
  the {\em $\ell$-torsion subgroup scheme} $A[\ell]$
  to be the kernel of multiplication by $\ell:A\to A$.
  The scheme $A[\ell]$ is a finite flat
  group scheme of length $\ell^{2g}$ over $k$,
  whereas the group of its $k$-rational points depends on $\ell$:
  $$
    A[\ell] (k) \,\iso\, \left\{
    \begin{array}{cl}
      (\ZZ/\ell\ZZ)^{2g} & \mbox{ if }\ell\neq p, \mbox{ and }\\
      (\ZZ/p\ZZ)^r & \mbox{ for some }0\leq r\leq g\mbox{ if }\ell=p.
    \end{array}
    \right.
  $$
  This integer $r$ is called the {\em $p$-rank} of $A$, and we have
  $$
    \dim_{\QQ_\ell} H^1(A,\QQ_\ell) \,=\,
    \left\{
    \begin{array}{cl}
      2g & \mbox{ if }\ell\neq p, \mbox{ and }\\
        r  & \mbox{ if }\ell=p.
    \end{array}
    \right.
  $$
  In particular, the assumption $\ell\neq p$ in property (5)
  is crucial, since we have $r<2g$ in any case.
  The group scheme $A[p]$ is of rank $p^{2g}$
  (although only rank $p^r$ can be ``seen'' via $k$-rational points),
  which should be reflected in the ``correct'' $p$-adic cohomology theory.
  Anticipating crystalline cohomology,
  which we will introduce in Section \ref{sec:crystalline cohomology} below, 
  there exists an isomorphism (see \cite[Th\'eor\`eme II.5.2]{Illusie deRham Witt})
  $$
     \Het{1}(A,\QQ_p)\otimes_{\QQ_p} K\,\iso\, 
     \left( \Hcris{1}(A/W)\otimes_W K\right)_{[0]}
     \,\subset\,\Hcris{1}(A/W)\otimes_W K,
  $$
  where the subscript $[0]$ denotes the slope zero sub-$F$-isocrystal
  ($W$ is the ring of Witt vectors of $k$, $K$ is its field of fractions,
  and we refer to Section \ref{sec:fcrystals} for details).
  Since $\Hcris{1}(A/W)$ is of rank $2g$,
  it gives the ``correct'' answer,
  and even the fact that 
   $\Het{1}(A,\QQ_p)$
  is ``too small'' can be explained using crystalline cohomology.
  In fact, we will see in Section \ref{sec:Serres observation} 
  that there exists no ``correct'' $p$-adic cohomology theory
  with $\QQ_p$-coefficients.
\end{Example}

\subsection{The ring of Witt vectors}
\label{subsec:Witt}
In the next section, we will introduce crystalline cohomology.
Since these cohomology groups are modules over the ring of Witt vectors, 
let us shortly digress on this ring.

Let $k$ be a perfect field of positive characteristic $p$.
For example, $k$ could be a finite field or algebraically closed.
Associated to $k$, there exists a ring $W(k)$, called the
{\em ring of Witt vectors} (or simply, {\em Witt ring}) 
of $k$, such that
\begin{enumerate}
 \item $W(k)$ is a discrete valuation ring of characteristic zero,
 \item the unique maximal ideal $\idealm$ of $W(k)$ is generated by $p$,
  and the residue field $W(k)/\idealm$ is isomorphic to $k$,
  \item $W(k)$ is complete with respect to the $\idealm$-adic topology,
  \item every $\idealm$-adically complete discrete valuation ring
    of characteristic zero with residue field $k$ contains $W(k)$ as subring,
  \item the Witt ring $W(k)$ is functorial in $k$.
\end{enumerate}
Note that property (4) shows that $W(k)$ is unique up to isomorphism.
There is also a more intrinsic characterization of the functor $k\mapsto W(k)$
as left adjoint to reduction modulo $p$, but we shall not pursue this here.

Instead, let us quickly run through an explicit construction of $W(k)$.
We refer to \cite[Chapitre II.6]{Serre Corps Locaux} and \cite[Section 17]{Hazewinkel}
for details, proofs, and generalizations.
Also, we refer to \cite{cd13} and \cite{cd14} for a completely different approach to the ring of
Witt vectors.
Let $p$ be a prime. 
Then, we define the {\em Witt polynomials} (with respect to $p$, which is understood
from the context and omitted) to be
the following polynomials with coefficients in $\ZZ$:
$$
\begin{array}{lcl}
  W_0(x_0) &:=& x_0\\
  W_1(x_0,x_1) &:=&x_0^p+px_1\\
  &...\\
  W_n(x_0,...,x_n) &:=& \sum_{i=0}^np^i x_i^{p^{n-i}} \,=\, x_0^{p^n}+px_1^{p^{n-1}}+...+p^nx_n
\end{array}
$$
Then, one can show that there exist unique polynomials $S_n$ and $P_n$ in
$2n+2$ variables with coefficients in $\ZZ$ such that the following holds true:
$$
 \begin{array}{ccccl}
    W_n(x_0,...,x_n) &+& W_n(y_0,...,y_n) &=& W_n(S_n(x_0,...,x_n,y_0,...,y_n)) \\
    W_n(x_0,...,x_n) &\cdot& W_n(y_0,...,y_n) &=& W_n(P_n(x_0,...,x_n,y_0,...,y_n))  
 \end{array}
$$
For an arbitrary ring $R$, which is not necessarily of characteristic $p$,
we define the {\em truncated Witt ring} $W_n(R)$ to be the 
set $R^n$, whose ring structure is defined to be
$$\begin{array}{rcl}
  (x_0,...,x_{n-1}) &\oplus& (y_0,...,y_{n-1}) \\
    &:=& (S_0(x_0,y_0),...,S_{n-1}(x_0,...,x_{n-1},y_0,...,y_{n-1}))\\
  (x_0,...,x_{n-1}) &\odot& (y_0,...,y_{n-1}) \\
    &:=& (P_0(x_0,y_0),...,P_{n-1}(x_0,...,x_{n-1},y_0,...,y_{n-1}))
\end{array}
$$
It turns out that $W_n(R)$ is indeed a ring with zero $0=(0,...,0)$ and unit $1=(1,0,...,0)$.
For example, we have $S_0(x_0,y_0)=x_0+y_0$ and $P_0(x_0,y_0)=x_0\cdot y_0$, and 
thus, $W_1(R)$ is just $R$ with its usual addition and multiplication.
Next, if $R$ is positive characteristic $p$, we define
$$\begin{array}{cclcl}
  V &:& (x_0,...,x_{n-1}) &\mapsto& (0,x_0,...,x_{n-2}) \\
  \sigma &:& (x_0,...,x_{n-1}) &\mapsto& (x_0^p,...,x_{n-1}^p) .
\end{array}
$$
Then, $V$ is called {\em Verschiebung} (German for ``shift''), and it is an additive map.
Next, $\sigma$ is called {\em Frobenius}, and it is a ring homomorphism.
(In order to avoid a clash of notations when dealing with $F$-crystals, see Section \ref{sec:fcrystals}
below, it is customary to denote the Frobenius on $W(k)$ by $\sigma$ rather than $F$.)
The maps $V$ and $\sigma$ are related to multiplication
by $p$ on $W_n(R)$ by 
$$
   \sigma\circ V \,=\,V\circ\sigma \,=\, p\cdot {\rm id}_{W_n(R)}
$$
We have $W_1(R)=R$ as rings, and for all $n\geq2$ the projection
$W_n(R)\to W_{n-1}(R)$ onto the first $(n-1)$ components is a surjective 
ring homomorphism.
Then, by definition, the {\em ring of Witt vectors} $W(R)$ is the inverse limit 
$$
    W(R)\,:=\,\varprojlim W_n(R),
$$
or, equivalently, the previous construction with respect to  the infinite product
$R^\NN$.

\begin{Exercise}
  For the finite field $\FF_p$, show that $W_n(\FF_p)\iso \ZZ/p^n\ZZ$ and thus, 
  $$
      W(\FF_p)\,\iso\,\varprojlim \ZZ/p^n\ZZ
  $$
  is isomorphic to $\ZZ_p$, the ring of $p$-adic integers. 
  Show that $\sigma$ is the identity and $V$ is multiplication by $p$
  in $W(\FF_p)$.
\end{Exercise} 

\begin{Exercise}
  If $k$ is a field of positive characteristic, 
  show that $W(k)$ is a ring of characteristic zero with residue field $k$.
  Set $\idealm:=\ker(W(k)\to k)$ and show
  that if $k$ is perfect, then $\idealm$ is generated by $p$, 
  that $W(k)$ is $\idealm$-adically complete, and that 
  $W(k)$ is a DVR.
  On the other hand if $k$ is not perfect, and if moreover $k/k^p$ is not a finite field
  extension, show that $\idealm$ is not finitely generated and that 
  $W_n(k)$ is not Noetherian for all $n\geq2$.
\end{Exercise}

If $X$ is a scheme, we can also sheafify the construction of $W_n(R)$ 
to obtain sheaves of rings  $W_n\OO_X$ and $W\OO_X$, respectively.
The cohomology groups
$$
  H^i(X,W_n\OO_X)\mbox{ \quad and \quad }
  H^i(X,W\OO_X)\,:=\,\varprojlim H^i(X,W_n\OO_X)
$$
were introduced and studied by Serre \cite{Serre Mexico}, they are 
called {\em Serre's Witt vector cohomology} groups, and we will come 
back to them in Section \ref{sec:formal groups}.
However, let us already note at this point that the torsion of the $W(k)$-module
$H^i(X,W\OO_X)$ may {\em not} be finitely generated
(for example, this is the case if $i=2$ and $X$ is a supersingular K3 surface), which is 
rather unpleasant.
Let us finally mention that the $W_n\OO_X$ are just the zeroth step
of the {\em de~Rham--Witt complex}
$(W_n\Omega_X^j,d)$ introduced by Illusie in
\cite{Illusie deRham Witt}, and we refer to 
\cite{Illusie Witt light} for an overview.

\subsection{Crystalline cohomology}
\label{sec:crystalline cohomology}
To a complex projective variety $X$, we have its de~Rham cohomology
$\HdR{\ast}(X/\CC)$ and showed in Section \ref{subsec:complex geometry}
that it comes with extra structure.
Now, let $X$ be a smooth and proper variety over a perfect
field $k$ of positive characteristic $p$.
Let us shortly summarize what we achieved so far 
in the algebraic setting:
\begin{enumerate}
\item
In Section \ref{subsec:algebraic deRham}, we associated
to $X$ its algebraic de~Rham-cohomology, which
is a $k$-vector space $\HdR{n}(X/k)$ together with a Frobenius
action and two filtrations, which is captured in the notion of an $F$-zip 
(Definition \ref{def: f-zip}).

On the other hand, there is no integral structure.
Another drawback is the following:
although there exists a Chern map $c_1:{\rm Pic}(X)\to \HdR{2}(X/k)$,
we have for every ${\cal L}\in\Pic(X)$
$$
  c_1({\cal L}^{\otimes p})\,=\,p\cdot c_1({\cal L})\,=\,0\mbox{ \quad in \quad }\HdR{2}(X/k)\,,
$$
giving a zero Chern class even for some very ample line bundles.
Also, counting fixed points via Lefschetz fixed point formulas 
(an important technique when dealing with varieties over finite fields)
gives us these numbers as traces in $k$, and thus, we obtain
the number of fixed points only 
as a congruence modulo $p$.

These observations suggest to look for cohomology theories whose
groups are modules over rings of characteristic zero.
\item
In Section \ref{subsec:l-adic cohomology}, we discussed
$\ell$-adic cohomology $\Het{n}(X,\QQ_\ell)$, which comes with an 
integral structure from $\Het{n}(X,\ZZ_\ell)$, but we have
no Hodge filtrations.

Moreover, if $\ell=p$, then we have seen in
Example \ref{example: wrong p-adic cohomology}
that $\Het{n}(X,\QQ_p)$ does not always give the desired answer.
Even worse, Serre showed that there does not exist
a ``well-behaved'' cohomology with $\QQ_p$-coefficients,
and we refer to Section \ref{sec:Serres observation} for a precise statement.

On the other hand, if $\ell\neq p$, then the fields $k$ and $\QQ_\ell$
usually have little in common, making comparison
theorems between de~Rham- and $\ell$-adic cohomology
even difficult to conjecture.
\end{enumerate}
It is here, that Witt vectors enter the picture.
As we shall now see, crystalline cohomology has all desired features
and provides an answer to all problems just raised.

To explain crystalline cohomology, let us assume for a moment that 
$X$ is smooth and projective over $k$ and that there exists 
a projective {\em lift} of $X$ to $W:=W(k)$, that is, a 
smooth projective scheme ${\cal X}\to\Spec W$
such that its special fiber 
${\cal X}\times_{\Spec W}\Spec k$ is isomorphic to $X$.
Then, for each $n$, the de~Rham-cohomology group
$\HdR{n}({\cal X}/W)$ is a finitely generated $W$-module.
It was Grothendieck's insight \cite{Grothendieck Crystal}
that these cohomology groups are independent of choice 
of lift $\cal X$ of $X$.
In fact, these cohomology groups can even be defined in case 
$X$ does not admit a lift to $W$.
The construction is quite involved, but we refer to
\cite[Section 1.3]{CL} for motivation and to \cite[Section 2]{CL}
for a detailed introduction.

More precisely, for every $m\geq1$, we have 
cohomology groups $\Hcris{\ast}(X/W_m(k))$,
all of which are finitely generated $W_m(k)$-modules.
For $m=1$, we obtain de~Rham-cohomology
$$
   \HdR{n}(X/k) \,\iso\, \Hcris{n}(X/W_1(k))\,=\, \Hcris{n}(X/k) \mbox{ for all }n\geq0,
$$
and, by definition, the limit
$$
  \Hcris{n}(X/W) \,:=\,\Hcris{n}(X/W(k)) \,:=\,\varprojlim \Hcris{n}(X/W_m(k))
$$
is called {\em crystalline cohomology}.
The origin of the name is as follows:
although $X$ may not lift to $W(k)$, its cohomology ``grows'' locally over $W$.
One can make these growths ``rigid'', so to glue and to obtain
a well-defined cohomology theory over $W(k)$.
And thus, growing and being rigid, it is natural to call such an object
a ``crystal''.

If $K$ denotes the field of fractions of $W$, then it has the following
properties:
\begin{enumerate}
 \item $\Hcris{n}(X/W)$ is a contravariant functor in $X$.
   These groups are finitely generated $W$-modules, and
   zero if $n<0$ or $n>2\dim(X)$.
 \item There is a cup-product structure
   $$
     \cup_{i,j}\,:\,\Hcris{i}(X/W)/{\rm torsion}\,\times\,\Hcris{j}(X/W)/{\rm torsion}\,\to\,\Hcris{i+j}(X/W)/{\rm torsion}
   $$
   Moreover, $\Hcris{2\dim(X)}(X/W)\iso W$, and 
   $\cup_{n,2\dim(X)-n}$ induces a perfect pairing, called {\em Poincar\'e duality}.
 \item $\Hcris{n}(X/W)$ defines an integral structure on $\Hcris{n}(X/W)\otimes_W K$,
  \item If $\ell$ is a prime different from $p$, then (see \cite{Katz Messing})
    $$
       b_n(X)\,\stackrel{{\rm def}}{=}\,\dim_{\QQ_\ell}\,\Het{n}(X,\QQ_\ell)\,=\,{\rm rank}_{W}\,\Hcris{n}(X/W),
    $$
    showing that crystalline cohomology computes $\ell$-adic Betti numbers.
   \item If $X$ lifts to $W$, then crystalline cohomology is isomorphic to 
    de~Rham cohomology of a lift,
    from which we deduce
    a universal coefficient formula 
    $$
       0\,\to\,\Hcris{n}(X/W)\otimes_W k\,\to\,\HdR{n}(X/k)\,\to\,{\rm Tor}_1^W(\Hcris{n+1}(X/W), k)\,\to\,0,
    $$
    for all $n\geq0$.
    This formula also holds true if $X$ does not lift.
    In any case, this shows that crystalline cohomology computes de~Rham cohomology.    
  \item Finally, there exists a Lefschetz fixed point formula, 
     there are base change formulas,
     there exist cycle classes in $\Hcris{2q}(X/W)$
     for codimension $q$ subvarieties,...
\end{enumerate}
By functoriality, the absolute Frobenius morphism $F:X\to X$  induces a 
$\sigma$-linear morphism $\varphi:\Hcris{n}(X/W)\to\Hcris{n}(X/W)$ of $W$-modules.
Ignoring torsion, this motivates to consider free $W$-modules together 
with injective $\sigma$-linear maps, which leads to the notion of an {\em $F$-crystal}, 
to which we come back in Section \ref{sec:fcrystals}.

We refer the interested reader to \cite{CL} for a much more detailed introduction to
crystalline cohomology, to \cite{Grothendieck Crystal}, \cite{Berthelot french} and 
\cite{Berthelot} for proofs and technical 
details, as well as to \cite{Illusie deRham Witt}  and \cite{Illusie Witt light} 
for the connection with the de~Rham--Witt complex.

\begin{Exercise}
  \label{exercise:crystalline free}
  Let $X$ be a smooth and proper variety over a perfect field $k$ 
  of positive characteristic $p$, and assume that the Fr\"olicher spectral 
  sequence degenerates at $E_1$.
  Using only the properties of crystalline cohomology mentioned above,
  show that the following are equivalent
  \begin{enumerate}
    \item For all $n\geq0$, the $W$-module $\Hcris{n}(X/W)$ is torsion-free.
    \item We have 
    $$
       \dim_{\QQ_\ell}\Het{n}(X,\QQ_\ell)\,=\,\dim_k \HdR{n}(X/k)
     $$
     for all $n\geq0$ and all primes $\ell\neq p$.
  \end{enumerate}
  Thus, the $p$-torsion of crystalline cohomology measures the deviation
  between $\ell$-adic Betti numbers to dimensions of de~Rham-cohomology.
\end{Exercise}

\begin{Examples}
  \label{example: crystalline cohomology abelian variety}
  Let us give a two fundamental examples.
  \begin{enumerate}
    \item Let $A$ be an Abelian variety of dimension $g$.
             Then, all $\Hcris{n}(A/W)$ are torsion-free
             $W$-modules. More precisely,
             $\Hcris{1}(A/W)$ is free of rank $2g$ and for all $n\geq2$
             there are isomorphisms
             $$
                  \Hcris{n}(A/W)\,\iso\,\Lambda^n \Hcris{1}(A/W). 
             $$
             Let us mention the following connection
             (for those familiar with $p$-divisible groups and Dieudonn\'e modules), 
             which we will not need in the sequel:
            let $A[p^n]$ be the kernel of multiplication
            by $p^n:A\to A$, which is a finite flat group scheme of rank $p^{2gn}$.
            By definition, the limit
            $$
                A[p^\infty] \,:=\, \varinjlim A[p^n]
            $$
            is the {\em $p$-divisible group} associated to $A$.
            Then, the {\em Dieudonn\'e-module} associated to $A[p^\infty]$
            is isomorphic to $\Hcris{1}(A/W)$, compatible with the Frobenius-actions 
            on both sides, see, for example, \cite[Section II.7.1]{Illusie deRham Witt}.
            We will come back to the Frobenius action on $\Hcris{1}(A/W)$
            in Section \ref{sec:fcrystals}.
     \item For a smooth and proper variety $X$, let
            $\alpha:X\to {\rm Alb}(X)$ be its Albanese morphism.
            Then, $\alpha$ induces an isomorphism
            $$
                 \Hcris{1}(X/W) \,\iso\, \Hcris{1}({\rm Alb}(X)/W).
           $$
           In particular, $\Hcris{1}(X/W)$ is always torsion-free.
           From this, we can compute the crystalline cohomology
           of curves via their Jacobians.
           We refer to \cite[Section II.5 and Section II.6]{Illusie deRham Witt}
           for connections of $p$-torsion of $\Hcris{2}(X/W)$ 
           with Oda's subspace of $\HdR{1}(X/k)$, 
           the non-reducedness of  the Picard scheme of $X$,
           as well as non-closed $1$-forms on $X$.
   \end{enumerate}
   In Section \ref{sec:K3}, we will compute the crystalline cohomology of a 
   K3 surface.
\end{Examples}

We already mentioned that Illusie constructed a complex $(W_m\Omega_{X/k}^j,d)$,
the de~Rham--Witt complex, and that it coincides with the de~Rham complex if $m=1$.
This complex gives rise to spectral sequences for all $m\geq1$
$$
  E_1^{i,j}\,:=\,H^j(X,W_m\Omega_{X/k}^i)\,\Rightarrow\,
  \Hcris{i+j}(X/W_m(k)).
$$
For $m=1$, this is the Fr\"olicher spectral sequence.
In the limit $m\to\infty$, this becomes the {\em slope spectral sequence} 
from Hodge--Witt cohomology $H^j(X,W\Omega_{X/k}^i)$ to crystalline cohomology.
Whereas the Fr\"olicher spectral sequence of $X$ may or may not degenerate at $E_1$
if $k$ is of positive characteristic, the slope spectral sequence
modulo torsion always degenerates at $E_1$.
Moreover, the slope spectral sequence (including torsion) degenerates at $E_1$ if and only if
the $p$-torsion of all $H^j(X,W\Omega_{X/k}^i)$ is finitely generated.
For $i=0$, this gives a conceptional framework for finite generation of Serre's
Witt vector cohomology.
We refer to \cite{Illusie deRham Witt} for details.

Finally, reduction modulo $p$ gives a map
$$
  \pi_n\,:\, \Hcris{n}(X/W)\,\to\,\HdR{n}(X/k),
$$
which, by the universal coefficient formula, is onto if and only
if $\Hcris{n+1}(X/W)$ has no $p$-torsion.
Thus, if all crystalline cohomology groups are torsion-free
$W$-modules, then de~Rham-cohomology is crystalline
cohomology modulo $p$.
Next, by functoriality, the Frobenius of $X$ induces a
$\sigma$-linear map $\varphi:\Hcris{\ast}(X/W)\to\Hcris{\ast}(X/W)$.
Under suitable hypotheses on $X$, the Frobenius action determines
the Hodge- and the conjugate filtration on de~Rham-cohomology.
More precisely, we have the following result of Mazur, and refer to
\cite[Section 8]{Berthelot} details, proofs, and further references.

\begin{Theorem}[Mazur]
  \label{thm: mazur's hodge theorem}
  Let $X$ be a smooth and proper variety over a perfect field $k$
  of positive characteristic $p$.
  Assume that $\Hcris{\ast}(X/W)$ has no $p$-torsion,
  and that the Fr\"olicher spectral sequence of $X$ degenerates
  at $E_1$.
  Then,
 $$
   \begin{array}{cclcc}
     \pi_n &\mbox{ maps }& \varphi^{-1}\left(  p^i\, \Hcris{n}(X/W) \right) &\mbox{ onto } & F^i  \\
     \pi_n\circ p^{-i} &\mbox{ maps } & {\rm Im}(\varphi)\cap \left( p^i\,\Hcris{n}(X/W) \right) &\mbox{ onto } & F_{\rm con}^{n-i},
   \end{array}
  $$
  where $F^i$ and $F_{\rm con}^{n-i}$ denote the Hodge- and its conjugate filtration on $\HdR{n}(X/k)$,
  respectively. 
\end{Theorem}

\subsection{Serre's observation}
\label{sec:Serres observation}
In this section, we have discussed $\ell$-adic and crystalline cohomology,
whose groups are $\QQ_\ell$-vector spaces and $W(k)$-modules, respectively.
One might wonder, whether crystalline cohomology
arises as base change from a cohomology theory, whose
groups are $\ZZ_p$-modules,
or even, whether all of the above cohomology theories 
arise from a cohomology theory, whose groups are $\ZZ$-modules 
or $\QQ$-vector spaces.
Now, cohomology theories that satisfy the ``usual'' properties 
discussed in this section are examples of so-called {\em Weil cohomology theories}, 
and we refer to \cite[Appendix C.3]{Hartshorne} for axioms and discussion.

Serre observed that there exists no Weil cohomology theory
in positive characteristic that take values 
in $\QQ$-, $\QQ_p$-, or $\RR$-vector spaces.
In particular, the above question has a negative answer.
Here is the sketch of a counter-example:
there exist supersingular elliptic curves $E$ over
$\FF_{p^2}$ such that ${\rm End}(E)\otimes\QQ$ is a quaternion
algebra that is ramified at $p$ and $\infty$.
By functoriality, we obtain a non-trivial representation of 
${\rm End}(E)$ on $H^1(E)$, which, being a Weil cohomology theory, must be 
$2$-dimensional.
In particular, we would obtain a non-trivial representation of
${\rm End}(E)\otimes\QQ$ in a $2$-dimensional $\QQ_p$- or $\RR$- vector space, 
a contradiction.
We refer to \cite[p. 315]{Grothendieck Crystal} or \cite[Section I.1.3]{CL}
for details.

\section{K3 Surfaces}
\label{sec:K3}

\subsection{Definition and examples}
In this section, we turn to K3 surfaces, and will compute the various
cohomology groups discussed in the previous section for them.
Let us first discuss their position within the classification of surfaces:
let $X$ be a smooth and projective surface over an algebraically closed field $k$
of characteristic $p\geq0$.
Moreover, assume that $\omega_X$ is numerically trivial, that is,
$\omega_X$ has zero-intersection with every curve on $X$.
In particular, $X$ is a minimal surface of Kodaira dimension zero.
By the Kodaira--Enriques classification (if $p=0$) and results of Bombieri and Mumford
(if $p>0$), then $X$ belongs to one of the following classes:
\begin{enumerate}
 \item Abelian surfaces, that is, Abelian varieties of dimension $2$.
 \item (Quasi-)hyperelliptic surfaces.
 \item K3 surfaces.
 \item Enriques surfaces.
\end{enumerate}
We refer to \cite[Chapter VI]{bhpv} for the surface classification over $\CC$,
and to \cite{Liedtke overview} for an overview in positive characteristic.
If $p\neq2,3$, the only surfaces with $\omega_X\iso\OO_X$
are Abelian surfaces and K3 surfaces. 
We refer the interested reader to \cite[Section 7]{Liedtke overview} for some
classes of Enriques surfaces in $p=2$, as well as quasi-hyperelliptic surfaces
in $p=2,3$ that have trivial canonical sheaves, and to
\cite{bm2} and \cite{bm3} for a detailed analysis of these surfaces.
Here, we are mainly interested in K3 surfaces, and recall
the following definition, which holds in any characteristic.

\begin{Definition}
 \label{def: k3}
 A {\em K3 surface} is a smooth and projective surface $X$ 
 over a field such that
 $$
  \omega_X\iso\OO_X\mbox{ \quad and \quad }h^1(X,\OO_X)=0.
 $$
\end{Definition}

\begin{Examples} 
 Let $k$ be an algebraically closed field of characteristic $p\geq0$.
 \begin{enumerate}
  \item If $X$ is a smooth quartic surface in $\PP^3_k$, then
    $\omega_X\iso\OO_X$ by the adjunction formula, and taking cohomology
    in the short exact sequence
    $$
      0\,\to\,\OO_{\PP^3_k}(-4)\,\to\,\OO_{\PP^3_k}\,\to\,\OO_X\,\to\,0
    $$
    we find $h^1(\OO_X)=0$.
    In particular, $X$ is a K3 surface.
   \item Similarly, smooth complete intersections of quadric and cubic hypersurfaces
    in $\PP^4_k$, as well as smooth complete intersections of three quadric hypersurfaces in
    $\PP^5_k$ give examples of K3 surfaces.
   \item If $p\neq 2$ and $A$ is an Abelian surface over $k$, then
    the quotient $A/\pm{\rm id}$ has $16$ singularities of type $A_1$, and its minimal
    resolution ${\rm Km}(A)$ of singularities is a K3 surface, the {\em Kummer surface
    associated to $A$}.
    (We refer the interested to reader to \cite{Shioda Kummer} and \cite{Katsura} to learn
     what goes wrong if $p=2$, and to \cite{Schroeer} how to remedy this.)
  \end{enumerate}
  We note that these three classes differ in size:
  the three example classes in (1) and (2) form $19$-dimensional families, whereas the
  Kummer surfaces in (3) form a $3$-dimensional family.
\end{Examples}

\subsection{Cohomological invariants}
\label{subsec: k3 cohomology}
Let us now compute the $\ell$-adic Betti numbers, the Hodge numbers,
and the crystalline cohomology groups of a K3 surface.
We will give all details so that the interested reader can see
where the characteristic-$p$ proofs are more difficult than the ones in
characteristic zero.

\begin{Proposition}
 \label{prop:Betti numbers}
 The $\ell$-adic Betti numbers of a K3 surface are as follows
 $$
    \begin{array}{c|ccccc}
         i       & 0 & 1 & 2 & 3 & 4 \\
                 \hline
      b_i(X) & 1 & 0 & 22 & 0 & 1
    \end{array}
    $$ 
    In particular, we have $c_2(X)=\sum_{i=0}^4 (-1)^i b_i(X)=24$.
\end{Proposition}

\prf
Since $X$ is a surface, we have $b_0=b_4=1$.
By elementary deformation theory of invertible sheaves,
$H^1(\OO_X)$ is the Zariski tangent space of $\Pic_{X/k}^0$ 
at the origin, see \cite[Section 3.3]{Sernesi}, for example.
Since $h^1(\OO_X)=0$ by definition of a K3 surface, it follows that
$\Pic_{X/k}^0$ is trivial.
Thus, also the Albanese variety ${\rm Alb}(X)$, 
which is the dual of the reduced Picard scheme,
is trivial, and we find $b_1(X)=2\dim{\rm Alb}(X)=0$.
By Poincar\'e duality, we have $b_1=b_3=0$.
Next, from Noether's formula for surfaces 
$$
12 \chi(\OO_X) \,=\, c_1(X)^2\,+\,c_2(X),
$$
we compute $c_2(X)=24$, which, together with the known Betti numbers,
implies $b_2(X)=22$.
\qed\medskip

Next, we recall that the {\em Hodge diamond} of a smooth projective variety $Y$
is given by ordering the dimensions $h^{i,j}(Y)=h^j(Y,\Omega_{Y/k}^i)$ in a rhombus.

\begin{Proposition}
 \label{prop:Hodge numbers}
 The Hodge diamond of a K3 surface is as follows:
 $$
 \begin{array}{c}
   h^{0,0} \\
   h^{1,0} \mbox{ \quad } h^{0,1}\\
   h^{2,0} \mbox{ \quad } h^{1,1} \mbox{ \quad } h^{0,2} \\
   h^{2,1} \mbox{ \quad } h^{1,2}\\
   h^{2,2}
 \end{array}
 \mbox{ \qquad $=$ \qquad }
  \begin{array}{c}
   1 \\
   0 \mbox{ \quad } 0\\
   1 \mbox{ \quad } 20 \mbox{ \quad } 1 \\
   0 \mbox{ \quad } 0\\
   1
 \end{array}
 $$  
\end{Proposition}

\prf
We have $h^{0,0}=h^{2,2}=1$ since $X$ is a surface, and 
 $h^{0,1}=0$ by the definition of a K3 surface.
Next, Serre duality gives $h^{0,1}=h^{2,1}$ and $h^{1,0}=h^{1,2}$.
If $k=\CC$, then complex conjugation induces the Hodge symmetry
$h^{1,0}=h^{0,1}$.
However, in positive characteristic, this Hodge symmetry may fail
in general
(see \cite{Serre Mexico} and \cite{Liedtke uniruled} for examples), 
and thus, we have to compute $h^{1,0}(X)$ another way:
using the isomorphism ${\cal E}^\vee\iso{\cal E}\otimes\det({\cal E})$,
which holds for locally free sheaves of rank $2$
(see \cite[Exercise II.5.16]{Hartshorne}, for example), we find
$T_X\iso\Omega^1_{X/k}$ for a K3 surface, and thus
$$
H^{1,0}(X)\,\stackrel{{\rm def}}{=}\,H^0(\Omega^1_{X/k}) \,\iso\, H^0(T_X)\,.
$$
Now, by a theorem of Rudakov and Shafarevich \cite{Rudakov Shafarevich}, a K3 surface has no
non-zero global vector fields, and thus, these cohomology groups are zero.
Finally, we use the Grothendieck--Hirzebruch--Riemann--Roch theorem to
compute
$$
\begin{array}{cl}
&  \chi(\Omega_{X/k}^1) \\
\,=\,&
{\rm rank}(\Omega_{X/k}^1)\cdot \chi(\OO_X)\,+\,
\frac{1}{2}\left( c_1(\Omega^1_{X/k})\cdot(c_1(\Omega^1_{X/k}) - K_X)  \right) 
\,-\,c_2(\Omega_{X/k}^1)\\
\,=\,&4\,+\,0\,-\,24\,=\,-20,
\end{array}
$$
which implies $h^1(\Omega_{X/k}^1)=20$.
\qed\medskip

As a consequence of this proposition, together with the 
Rudakov--Shafarevich theorem on non-existence of global vector fields
on K3 surfaces, we obtain

\begin{Proposition}
  \label{prop: k3 frolicher degenerates}
  For a K3 surface $X$, the Fr\"olicher spectral sequence
  $$
    E_1^{i,j}\,=\,H^j(X,\Omega_{X/k}^i)\,\Rightarrow\,\HdR{i+j}(X/k)
  $$
  degenerates at $E_1$.
  Moreover, $\Hcris{n}(X/W)$ is a free $W$-module of rank $b_n(X)$ for all $n\geq0$.
\end{Proposition}

\prf
By Proposition \ref{prop:Hodge numbers} and using
the isomorphism $\Omega_X^1\iso T_X$ (seen in the proof
of Proposition \ref{prop:Hodge numbers}), we have
$H^2(T_X)\iso H^2(\Omega_X^1)=H^{1,2}=0$.
Since $H^2(T_X)=0$, deformations of $X$ are unobstructed, and thus, 
$X$ lifts to $W_2(k)$ (see also the discussion in Section \ref{subsec: deformations k3}
and Theorem \ref{thm: K3 deformation theory} below), 
and thus, degeneracy of the Fr\"olicher spectral sequence at $E_1$ 
follows from Theorem \ref{thm: Deligne Illusie}.
From Proposition \ref{prop:Hodge numbers} and Exercise
\ref{exercise: froelicher} we compute the dimensions of the
de~Rham cohomology groups, which then turn out to be the same as the
$\ell$-adic Betti numbers given in Proposition \ref{prop:Betti numbers}.
Thus, by Exercise \ref{exercise:crystalline free},
the crystalline cohomology groups are free $W$-modules of the stated rank.
\qed\medskip

\begin{Remark}
  For a smooth and proper variety $X$ over a perfect field $k$,
  the slope spectral sequence from Hodge--Witt to crystalline cohomology
  degenerates at $E_1$ if and only if all the $W$-modules $H^j(X,W\Omega_{X/k}^i)$
  are finitely generated \cite[Th\'eor\`eme II.3.7]{Illusie deRham Witt}.
  For a K3 surface, this is the case if and only if it
  is not supersingular - we refer to 
  Section \ref{sec:supersingularK3} for definition of supersingularity
  and to \cite[Section II.7.2]{Illusie deRham Witt} for details.
\end{Remark}

\subsection{Deformation theory}
\label{subsec: deformations k3}
Infinitesimal and formal deformations of a smooth and proper variety $X$
over a field $k$ can be controlled by a tangent--obstruction theory arising from
the $k$-vector spaces $H^i(T_X), i=0,1,2$,
see \cite[Chapter 2]{Sernesi} or \cite[Chapter 6]{Fantechi Goettsche} for a
reader-friendly introduction.

Let us recall the most convenient case:
if $H^2(T_X)=0$, then every infinitesimal
deformation of order $n$ can be extended to one of order $n+1$,
and then, the set of all such extensions is an affine space under
$H^1(T_X)$.
In particular, this applies to lifting problems:
if $k$ is perfect of positive characteristic, then
the Witt ring $W(k)$ is a limit of rings $W_n(k)$,
see Section \ref{subsec:Witt}.
Since the kernel of $W_{n+1}(k)\to W_n(k)$ is the ideal generated by 
$p^{n-1}$ and $(p^{n-1})^2=0$, it is a small extension, and thus,
a smooth and proper variety $X$ over a perfect field $k$ of positive
characteristic with $H^2(T_X)=0$ admits a formal lift to 
$W(k)$, and we refer to \cite[Chapter 8.5]{Illusie Existence Theorem}
and \cite[Section 11.2]{Liedtke overview} for details and further references.
Since this most convenient case applies to K3 surfaces,
we have the following result.

\begin{Theorem}
  \label{thm: K3 deformation theory}
  Let $X$ be a K3 surface over a perfect field $k$
  of positive characteristic.
  Then, the formal deformation space ${\rm Def}(X)$ of $X$ is 
  smooth of relative dimension $20$ over $W(k)$, that is,
   $$
       {\rm Def}(X) \,\iso\, {\rm Spf}\,W(k)[[t_1,...,t_{20}]]\,.
   $$ 
   In particular, $X$ formally lifts over $W(k)$.
\end{Theorem}

\prf
By Proposition \ref{prop:Hodge numbers} and using
the isomorphism $\Omega_X^1\iso T_X$ (seen in the proof
of Proposition \ref{prop:Hodge numbers}), we find
$$
h^0(T_X)\,=\,h^2(T_X)\,=\,0\mbox{ \quad and \quad }
h^1(T_X)\,=\,20,
$$
from which all assertions follow from standard results of deformation theory,
see \cite[Chapter 2]{Sernesi} or \cite[Chapter 6]{Fantechi Goettsche}, for
example.
\qed\medskip

If $X$ is a K3 surface over a perfect field $k$ of positive
characteristic, then the previous theorem implies that there exists
a compatible system $\{{\cal X}_n\to\Spec W_n(k)\}_n$ of 
{\em algebraic} schemes ${\cal X}_n$, each flat over $W_n(k)$, 
and with special fiber ${\cal X}_1=X$.
Now, the limit of this system is a 
{\em formal} scheme \cite[Section II.9]{Hartshorne}, and it is not clear
whether it is algebraizable, that is, we do not know, whether this 
limit arises as completion of a scheme over $W(k)$ along its special
fiber (in fact, it is not true in general, 
see Section \ref{sec: moduli spaces k3} below).

By Grothendieck's existence theorem
(see \cite[Theorem 8.4.10]{Illusie Existence Theorem}, for example),
algebraization of formal schemes holds, for example, if one is able to equip 
${\cal X}_n$ with a compatible system ${\cal L}_n$ of invertible sheaves on ${\cal X}_n$
such that ${\cal L}_1$ is ample on ${\cal X}_1=X$.
This poses the question whether a given formal deformation can be equipped
with such a compatible system of invertible sheaves.
The obstruction to deforming an invertible sheaf to a small extension
lies in $H^2(\OO_X)$, which is $1$-dimensional for a K3 surface.
We thus expect that this should impose one non-trivial equation 
to ${\rm Def}(X)$, which is true and made precise by 
the following results of Deligne, \cite[Proposition 1.5]{Deligne K3} 
and \cite[Th\'eor\`eme 1.6]{Deligne K3}.

\begin{Theorem}[Deligne]
  \label{thm: deligne deforms line bundle}
  Let $X$ be a K3 surface over a perfect field $k$ of positive
  characteristic, and let ${\cal L}$ be a non-trivial invertible sheaf on $X$.
  Then, the space ${\rm Def}(X,{\cal L})$
  of formal deformations of the pair $(X,{\cal L})$ is a formal Cartier
  divisor inside ${\rm Def}(X)$, that is,
  $$
      {\rm Def}(X,{\cal L})\,\subset\,{\rm Def}(X),
  $$
  is a formal subscheme defined by one equation.
  Moreover, ${\rm Def}(X,{\cal L})$ is flat and of relative dimension
  $19$ over $W(k)$.
\end{Theorem}

Unfortunately, it is not clear whether ${\rm Def}(X,{\cal L})$
is smooth over $W(k)$, and we refer to \cite[\S2]{Ogus} for an
analysis of its singularities.
In particular, if we pick an ample invertible sheaf $\cal L$
on $X$ in order to construct a formal lift of the pair $(X,{\cal L})$ to
$W(k)$ (in order to apply Grothendieck's existence theorem),
then it could happen that ${\rm Def}(X,{\cal L})$ is flat,
but not smooth over $W(k)$.
Thus, a priori, we only have an algebraic lift of $X$ 
to some finite extension ring $R\supseteq W(k)$.
However, thanks to a refinement of Ogus \cite[Corollary 2.3]{Ogus} 
of Deligne's result, we have

\begin{Theorem}[Deligne, Ogus]
  Let $X$ be a K3 surface over an algebraically closed field of odd 
  characteristic. Then, there exists a projective lift of $X$ to $W(k)$.
\end{Theorem} 

\prf
By \cite[Corollary 2.3]{Ogus}, any nonsuperspecial K3 surface can be
lifted projectively to $W(k)$, and we refer to
\cite[Example 1.10]{Ogus} for the notion  of superspecial K3 surfaces.
Since the Tate-conjecture holds for K3 surfaces in odd characteristic
(see Theorem \ref{thm: Tate conjecture} below),
the only nonsuperspecial K3 surface
is the supersingular K3 surface with Artin invariant $\sigma_0=1$,
see \cite[Remark 2.4]{Ogus}.
However, this latter surface is the Kummer surface associated
to the self-product of a supersingular elliptic curve 
by \cite[Corollary 7.14]{Ogus} and can be
lifted ``by hand'' projectively to $W(k)$.
\qed\medskip

\subsection{Moduli spaces}
\label{sec: moduli spaces k3}
By Theorem \ref{thm: K3 deformation theory}, the formal
deformation space ${\rm Def}(X)$ of a K3 surface
$X$ is formally smooth and $20$-dimensional over $W(k)$.
However, it is not clear (and in fact, not true)
whether all formal deformations are algebraizable.
By a theorem of Zariski and Goodman
(see \cite[Theorem 1.28]{Badescu}, for example),
a smooth and proper algebraic surface is automatically projective,
which applies in particular to K3 surfaces.
Thus, associated to an ample invertible sheaf $\cal L$ 
on an algebraic K3 surface $X$, there is a
formal Cartier divisor ${\rm Def}(X,{\cal L})\subset{\rm Def}(X)$ 
by Theorem \ref{thm: deligne deforms line bundle},
along which $\cal L$ extends.
Since formal and polarized deformations are algebraizable
by Grothendieck's existence theorem (see Section \ref{subsec: deformations k3}), we can algebraize
the $19$-dimensional formal family over ${\rm Def}(X,{\cal L})$.
Using Artin's approximation theorems, this latter family can be descended to a
polarized family of K3 surfaces that is $19$-dimensional and of finite type
over $W(k)$, and one may think of it as an \'etale neighborhood
of $(X,{\cal L})$ inside a moduli space of suitably polarized K3 surfaces.
In fact, this can be made precise to give a rigorous algebraic construction 
for moduli spaces of polarized K3 surfaces, and we refer to
\cite[Section 5]{Artin moduli} and \cite{Rizov} for technical details.
These moduli spaces of polarized K3 surfaces are $19$-dimensional,
whereas the unpolarized formal deformation space is $20$-dimensional.

Now, before proceeding, let us shortly leave the algebraic world:
over the complex numbers, there exists a $20$-dimensional analytic 
moduli space for compact K\"ahler surfaces that are of type K3, and
most of which are not algebraic.
Moreover, this moduli space is smooth, but not Hausdorff.
Inside it, the set of algebraic K3 surfaces is a countable union of analytic divisors.
In fact, these divisors correspond to moduli spaces of algebraic and polarized K3 surfaces.
We refer to \cite[Chapter VIII]{bhpv} for details and further references.
We mention this to convince the reader that it is not possible to 
obtain a $20$-dimensional moduli space of algebraic K3 surfaces,
even over the complex numbers.

Therefore, when considering moduli of algebraic K3 surfaces, 
one usually looks at moduli spaces of (primitively) polarized surfaces.
Here, an invertible sheaf ${\cal L}$ on a variety $X$ is called
{\em primitive} if it is not of the form ${\cal M}^{\otimes k}$ for some $k\geq2$.
Then, for a ring $R$, we consider the functor
${\cal M}_{2d,R}^\circ$
$$
\begin{array}{ccc}
\left(\begin{array}{c}\mbox{schemes}\\ \mbox{over }R\end{array}\right) &\to& 
\left(\mbox{groupoids}\right) \\
S&\mapsto&\left\{
\begin{array}{l}
 \mbox{ flat morphisms of algebraic spaces $({\cal X},{\cal L})\to S$,}\\
 \mbox{ all of whose geometric fibers are K3 surfaces, and }\\
 \mbox{ such that $\cal L$ restricts to a primitive polarization of } \\
 \mbox{ degree (=self-intersection) $2d$ on each fiber }
\end{array}
\right\}
\end{array}
$$
In view of the above discussion it should not be surprising that
this functor is representable by a separated Deligne--Mumford
stack \cite[Theorem 4.3.3]{Rizov}.

Combining \cite[Proposition 4.3.11]{Rizov}, \cite[Section 5]{Maulik}, and
\cite[Corollary 4.16]{Madapusi Pera}, we even have 
the following results on the global geometry of these moduli spaces.
Before stating them, let us mention that (2) and (3) are neither obvious nor straight forward,
and partly rest on the Kuga--Satake construction, which allows to obtain
results about moduli spaces of K3 surfaces from the corresponding results
for moduli spaces of Abelian
varieties.

\begin{Theorem}[Madapusi-Pera, Maulik, Rizov]
  \label{thm: moduli spaces are nice}
  The Deligne--Mumford stack
  \begin{enumerate}
    \item ${\cal M}_{2d,\ZZ[\frac{1}{2d}]}^\circ$ is smooth over $\ZZ[\frac{1}{2d}]$.
    \item ${\cal M}_{2d,\FF_p}^\circ$ is quasi-projective over $\FF_p$ if $p\geq5$ and $p\nmid d$.
    \item ${\cal M}_{2d,\FF_p}^\circ$ is 
        geometrically irreducible over $\FF_p$ if $p\geq3$ and $p^2\nmid d$.
  \end{enumerate}
\end{Theorem}

In Section \ref{sec:stratification} below, we will introduce and discuss 
a stratification of ${\cal M}_{2d,\FF_p}^\circ$,
which only exists in positive characteristic.
It would be interesting to understand the geometry and the 
singularities (if there are any) of ${\cal M}_{2d,\FF_p}$  if $p$ divides $2d$.

\section{F-crystals}
\label{sec:fcrystals}

In Section \ref{sec:crystalline cohomology}, we introduced crystalline cohomology, 
and we computed it for K3 surfaces in Section \ref{subsec: k3 cohomology}.
Just as Hodge structures abstractly capture the linear algebra data coming from de~Rham-cohomology
of a smooth complex projective variety (see the end of Section \ref{subsec:complex geometry}), 
{\em $F$-crystals} capture the semi-linear data coming from crystalline cohomology.
In this section, after introducing $F$-crystals, we associate to them
two polygons: the {\em Hodge polygon} and the {\em Newton polygon}.
Under the assumptions of Theorem \ref{thm: mazur's hodge theorem}
(which hold, for example, for K3 surfaces), 
the $F$-crystal associated to $\Hcris{\ast}(X/W)$
not only computes de~Rham-cohomology, but also the Hodge
and its conjugate filtration.
On the other hand, by a result of Dieudonn\'e and Manin,
we can classify $F$-crystals up to {\em isogeny} in terms of {\em slopes},
which gives rise to a second polygon, the Newton polygon, which
always lies on or above the Hodge polygon.
The deviation between Hodge and Newton polygon gives
rise to new discrete invariants of varieties in positive characteristic
that have no analog in characteristic zero.
In Section \ref{sec:stratification}, we will use these discrete invariants to stratify
the moduli space of K3 surfaces.

\subsection{$\mathbf F$-crystals}
In Section \ref{subsec:Witt}, we introduced and discussed the Witt ring $W:=W(k)$
for a perfect field $k$ of positive characteristic $p$.
We denote by $K$ its field of fractions.
We also recall that the Frobenius morphism $x\mapsto x^p$ of $k$
induces a ring homomorphism $\sigma:W\to W$ by functoriality, and
that there exists an additive map $V:W\to W$ such that
$p=\sigma\circ V=V\circ\sigma$.
In particular, $\sigma$ is injective.

\begin{Definition}
  An {\em $F$-crystal} $(M,\varphi_M)$ over $k$ is a free $W$-module $M$
  of finite rank together with an injective and $\sigma$-linear 
  map $\varphi_M:M\to M$, that is,
  $\varphi_M$ is additive, injective, and satisfies
  $$
     \varphi_M(r\cdot m)\,=\,\sigma(r)\cdot \varphi_M(m)\mbox{ \quad for all }
     r\in W,\, m\in M.
  $$
  An {\em $F$-isocrystal} $(V,\varphi_V)$ is a finite dimensional $K$-vector space $V$
  together with an injective and $\sigma$-linear map $\varphi_V:V\to V$. 
  
  A {\em morphism} $u:(M,\varphi_M)\to(N,\varphi_N)$ of $F$-crystals (resp., $F$-isocrystals) is a
  $W$-linear (resp., $K$-linear) map $M\to N$ such that $\varphi_N\circ u=u\circ\varphi_M$.
  An {\em isogeny} of $F$-crystals is a morphism $u:(M,\varphi_M)\to(N,\varphi_N)$ of $F$-crystals, 
  such that the induced map $M\otimes_W K\to N\otimes_W K$
  is an isomorphism of $F$-isocrystals.
\end{Definition}

Let us give two examples of $F$-crystals, one arising from geometry (and being the prototype of
such an object), the other one purely algebraic (and being crucial for the isogeny classification
later on).

\begin{Example}
  \label{example:geometry}
  Let $X$ be a smooth and proper variety over $k$.
  Then, for every $n\geq0$,
  $$
     H^n\,:=\,\Hcris{n}(X/W)/{\rm torsion}
  $$
  is a free $W$-module of finite rank.
  By functoriality, the absolute Frobenius morphism $F:X\to X$ induces
  a $\sigma$-linear map $\varphi:H^n\to H^n$.
  Next, Poincar\'e duality induces a perfect pairing
  $$
    \langle -,-\rangle\,:\,  H^n \times H^{2\dim(X)-n} \,\to\, H^{2\dim(X)}(X/W)\iso W, 
  $$
  which satisfies the following compatibility with Frobenius
  $$
     \langle \varphi(x),\,\varphi(y)\rangle \,=\, p^{\dim(X)}\cdot\sigma \langle x,y\rangle.
  $$
  Since $\sigma$ is injective on $W$, it follows that also
  $\varphi:H^n\to H^n$ is injective, and thus,
  $(H^n,\varphi)$ is an $F$-crystal.
\end{Example}

\begin{Example}
  \label{example:algebra}
  Let $W_\sigma\langle T\rangle$ be the non-commutative polynomial ring 
  in the variable $T$ over $W$ subject to the relations
  $$
      T\cdot r \,=\, \sigma(r)\cdot T\mbox{ \quad for all \quad }r\in W.
  $$
  Let $\alpha=r/s\in\QQ_{\geq0}$, where $r,s$ are non-negative and coprime
  integers.
  Then,
  $$
     M_\alpha \,:=\, W_\sigma\langle T\rangle/(T^s-p^r)
  $$
  together with $\varphi:m\mapsto T\cdot m$ defines an
  $F$-crystal $(M_\alpha,\varphi)$ of rank $s$.
  The rational number $\alpha$ is called the {\em slope} of $(M_\alpha,\varphi)$.
\end{Example}

The importance of the previous example comes from the following 
result, which classifies $F$-crystals over algebraically closed fields
up to isogeny:

\begin{Theorem}[Dieudonn\'e--Manin]
  \label{thm: dieudonne manin}
  Let $k$ be an algebraically closed field of positive characteristic.
  Then, the category of $F$-crystals over $k$
  up to isogeny is semi-simple and
  the simple objects are the $(M_\alpha,\varphi)$, $\alpha\in\QQ_{\geq0}$
  from Example \ref{example:algebra}.
\end{Theorem}

We note in passing that not every $F$-isocrystal is of the 
form $M\otimes_{W(k)}K$ for some $F$-crystal $(M,\varphi)$.
Those $F$-isocrystals that do are called {\em effective}.

\begin{Definition}
  Let $(M,\varphi)$ be an $F$-crystal over an algebraically closed field $k$
  of positive characteristic, and let 
  $$
     (M,\varphi)\,\sim\, \bigoplus_{\alpha\in\QQ_{\geq0}} M_\alpha^{n_\alpha}
  $$
  be its decomposition up to isogeny according to Theorem \ref{thm: dieudonne manin}.
  Then, the elements in the set
  $$
     \{\, \alpha\,\in\,\QQ_{\geq0} \,|\, n_\alpha\neq0 \,\}
  $$
  are called the {\em slopes} of $(M,\varphi)$.
  For every slope $\alpha$ of $(M,\varphi)$, the integer 
  $$
      \lambda_\alpha \,:=\, n_\alpha\cdot{\rm rank}_W\, M_\alpha
  $$
  is called the {\em multiplicity} of the slope $\alpha$.
  In case $(M,\varphi)$ is an $F$-crystal over a perfect field $k$,
  we define its slopes and multiplicities to be the ones of
  $(M,\varphi)\otimes_{W(k)} W(\overline{k})$, where $\overline{k}$ is an algebraic
  closure of $k$. 
\end{Definition}

\subsection{Newton and Hodge polygons}
\label{subsec: newton hodge polygon}
Let $(M,\varphi)$ be an $F$-crystal.
We order its slopes in ascending order
$$
0\,\leq\,\alpha_1\,<\,\alpha_2\,<\,...\,<\alpha_t
$$
and denote by $\lambda_1,...,\lambda_t$ the respective multiplicities.

Then, the {\em Newton polygon} of $(M,\varphi)$ is defined to be the
graph of the piecewise linear function ${\rm Nwt}_M$ from the interval 
$[0,{\rm rank}\,M]\subset\RR$ to $\RR$, such that 
${\rm Nwt}_M(0)=0$ and whose graph has the following slopes 
$$
\begin{array}{lcl}
 \mbox{ slope } & \alpha_1 & \mbox{ if }0\leq t <\lambda_1,\\
 \mbox{ slope } & \alpha_2 & \mbox{ if }\lambda_1\leq t <\lambda_1+\lambda_2,\\
 \multicolumn{3}{c}{...}
\end{array}
$$
Since we ordered the slopes in ascending order, this polygon is convex.
Next, it follows easily from the definitions that the vertices of this polygon
have integral coordinates.
Clearly, since the Newton polygon is built from slopes, it only depends on the
isogeny class of $(M,\varphi)$.
Conversely, we can read off all slopes and multiplicities from the Newton polygon,
and thus, the Newton polygon actually determines the $F$-crystal up to isogeny.
\medskip 

Next, we define the {\em Hodge polygon} of $(M,\varphi)$, whose definition
is motivated by Theorem \ref{thm: mazur's hodge theorem}, but see also
Theorem \ref{thm:mazur nygaard ogus} below.
Since $\varphi$ is injective, $M/\varphi(M)$ is an Artinian $W$-module,
and thus, there exist non-negative integers $h_i$ and an isomorphism
$$
    M/\varphi(M)\,\iso\,\bigoplus_{i\geq1} (W/p^iW)^{h_i}.
$$
Moreover, we define
$$
   h_0\,:=\,{\rm rank}\, M \,-\,\sum_{i\geq1}h_i.
$$
Then, the Hodge polygon of $(M,\varphi)$ is defined to be the
graph of the piecewise linear function ${\rm Hdg}_M$ from the interval 
$[0,{\rm rank}\,M]\subset\RR$ to $\RR$, such that
${\rm Hdg}_M(0)=0$ and whose graph has the following slopes 
$$
\begin{array}{lcl}
 \mbox{ slope } & 0 & \mbox{ if }0\leq t <h_0,\\
 \mbox{ slope } & 1 & \mbox{ if }h_0\leq t <h_0+h_1,\\
 \multicolumn{3}{c}{...}
\end{array}
$$
As above with the Newton polygon, the Hodge polygon is convex
and its vertices have integral coordinates.

\begin{Example}
 \label{example different hodge polygon}
 Let $M_\alpha$ with $\alpha=r/s\in\QQ_{\geq0}$ be the $F$-crystal 
 from Example \ref{example:algebra}.
 Then,
 $$
   M_\alpha/\varphi(M_\alpha) \,\iso\, W/p^r W,
 $$
 we obtain $h_0=s-1$, $h_r=1$, and $h_i=0$ for $i\neq 0,r$.

 Now, we assume that $0<\alpha<1$, that is, $0<r<s$, and
 we define a $W$-module $N_\alpha$ together with an embedding into
 $M_\alpha\otimes_WK$ as follows: 
 $$
  \begin{array}{ccc}
    N_\alpha \,:=\, W[T,U]/(TU-p, T^{s-r}-U^r) &\to& M_\alpha\otimes_W K \\
    T&\mapsto&T\\
    U&\mapsto& pT^{-1}=p^{1-r}T^{s-1}.
  \end{array}
 $$
 Then, $N_\alpha$ inherits the structure of an $F$-crystal from $M_\alpha\otimes_WK$,
 and it is isogenous to $M_\alpha$.
 Since
 $$
   N_\alpha/\varphi(N_\alpha) \,\iso\, (W/pW)^r,
 $$
 we obtain $h_0=s-r$, $h_1=r$, and $h_i=0$ for $i\neq0,1$.
 In particular, we find the following Hodge polygons (solid lines)
 and Newton polygons (dotted lines):
 
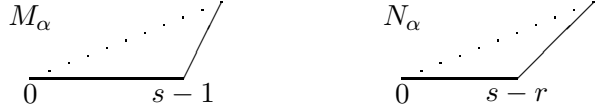
\begin{center}
 \begin{figure}[htbp]
  \setlength{\unitlength}{.1in}
  \centering
  \begin{picture}(10,4)
  \put(0,0.8){\line(1,0){8}}
  \put(8,0.8){\line(1,2){2}}
  \put(0,0.8){\qbezier[10](0,0)(5,2)(10,4)}
  \put(0,0){\makebox(0,0){$0$}}
  \put(8,0){\makebox(0,0){$s-1$}}
  \put(0,4){\makebox(0,0){$M_\alpha$}}
  \end{picture}
  \qquad
  \qquad
  \qquad  
  \begin{picture}(10,4)
  \put(0,0.8){\line(1,0){6}}
  \put(6,0.8){\line(1,1){4}}
  \put(0,0.8){\qbezier[10](0,0)(5,2)(10,4)}
  \put(0,0){\makebox(0,0){$0$}}
  \put(6,0){\makebox(0,0){$s-r$}}
  \put(0,4){\makebox(0,0){$N_\alpha$}}
  \end{picture}
\end{figure}
\end{center}
\end{Example}

This example shows that the Hodge polygon, unlike the Newton polygon, 
is in general {\em not} an isogeny invariant of the $F$-crystal $(M,\varphi)$.
However, the isogeny class of an $F$-crystal, that is,
its Newton polygon, puts restrictions on the possible
Hodge polygons:

\begin{Proposition}
   Let $(M,\varphi)$ be an $F$-crystal. 
   Then, its Newton-polygon lies on or above its Hodge polygon, and both have
   the same startpoint and endpoint:
   $$
   \begin{array}{lcll}
       {\rm Nwt}_M(t) &\geq& {\rm Hdg}_M(t) & \mbox{ for all } t\in [0,{\rm rank}\,M], \mbox{ and } \\
       {\rm Nwt}_M(t) &=&{\rm Hdg}_M(t) & \mbox{ for }t_0=0\mbox{ and }t={\rm rank}\, M.
    \end{array}
   $$
\end{Proposition}

\subsection{$\mathbf{F}$-crystals arising from geometry}
Now, we link the Hodge polygon of an $F$-crystal that arises from crystalline
cohomology of a variety to the Hodge numbers of that variety, which justifies
the terminology:
let $X$ be a smooth and proper variety over $k$ and fix an integer
$n\geq0$.
Then, we consider the Hodge numbers
$$
\widetilde{h}_i\,:=\,h^{i,n-i}\,=\,\dim_k H^{n-i}(X,\Omega_{X/k}^i) 
\mbox{ \quad for all } 0\leq i\leq n
$$
and, as before with the Hodge polygon, we construct from 
these integers a piecewise linear function
$$
\widetilde{{\rm Hdg}}_X^n\,:\,\left[0,\,\widetilde{h}^n_{\rm dR}\right]\,\to\,\RR,
\mbox{ \quad where \quad }
\widetilde{h}^n_{\rm dR} \,=\,\sum_{i=0}^n \widetilde{h}_i,
$$
and whose associated convex polygon is called the
{\em geometric Hodge polygon}.

The following deep and important result shows that under extra hypotheses 
the $F$-crystal associated to
the crystalline cohomology of a smooth and proper variety detects its Hodge numbers.
In fact, part of this result just rephrases Theorem \ref{thm: mazur's hodge theorem}
in terms of Hodge polygons.

\begin{Theorem}[Mazur, Nygaard, Ogus]
  \label{thm:mazur nygaard ogus}
  Let $X$ be a smooth and proper variety over a perfect 
  field $k$ of positive characteristic.
  Fix an integer $n\geq0$ and let
  $$
     H^n\,:=\,(\Hcris{n}(X/W)/{\rm torsion},\varphi)
  $$
  be the associated $F$-crystal.
  Then
  \begin{enumerate}
    \item  For all $t\in[0,{\rm rank}\, H^n]$, we have
      $$
         {\rm Nwt}_{H^n}(t)\,\geq\,\widetilde{{\rm Hdg}}_{X}^n(t).
      $$
    \item If $\Hcris{n}(X/W)$ is torsion-free, and if the Fr\"olicher spectral sequence
      of $X$ degenerates at $E_1$, then for all $t\in[0,{\rm rank}\, H^n]$, we have
      $$
         {\rm Hdg}_{H^n}(t)\,=\,\widetilde{{\rm Hdg}}_X^n(t).
      $$
      In particular, $H^n$ computes all Hodge numbers $h^{i,n-i}$ of $X$.
  \end{enumerate} 
\end{Theorem}

The following exercise shows that there are restrictions on 
the slopes of $F$-crystals that arise as crystalline
cohomology of varieties.

\begin{Exercise}
  \label{exc: slope restriction}
  Let $X$ be smooth and proper variety of dimension $d$,
  and let $(H^n,\varphi)$, $n=0,...,2d$ be the $F$-crystals
  $\Hcris{n}(X/W)/{\rm torsion}$ as above.
  \begin{enumerate}
    \item Using Poincar\'e duality, show that
      $\varphi\circ\varphi^\vee=p^d\cdot{\rm id}$ for all $n\geq0$,
      and deduce that the slopes of $H^n$ lie inside the interval
      $[0,d]$.
    \item Use the hard Lefschetz theorem together with Poincar\'e duality
      to show that the slopes of $H^n$ lie inside the interval
      $$\begin{array}{cl}
        [0,n]    &  \mbox{ if } 0\leq n\leq d \\{}
        [n-d,d] &  \mbox{ if } d\leq n \leq 2d\,.
      \end{array}$$
  \end{enumerate} 
\end{Exercise}

We refer to \cite{Katz crystal} for more about crystals and their slopes.

\subsection{Abelian varieties}
\label{subsec: Abelian varieties crystals}
Let $A$ be an Abelian variety of dimension $g$ over an algebraically 
closed field $k$ of positive characteristic.
As already mentioned in 
Examples \ref{example: crystalline cohomology abelian variety}, 
there exist isomorphisms
$$ 
   \Hcris{n}(A/W) \,\iso\,\Lambda^n \Hcris{1}(A/W)
   \mbox{ \quad for all }n\geq0.
$$
In fact, these isomorphisms are compatible with
Frobenius actions on both sides, and thus, are
isomorphisms of $F$-crystals.
In particular, it suffices to understand the $F$-crystal $\Hcris{1}(A/W)$,
which is a free $W$-module of rank $2g$.
Thus, $\Hcris{\ast}(A/W)$ is torsion-free and since
the Fr\"olicher spectral sequence of $A$ 
degenerates at $E_1$ (see Section \ref{subsec:algebraic deRham}),
the assumptions of Theorem \ref{thm:mazur nygaard ogus}
are fulfilled.
Let us now discuss the two cases $g=1$ and $g=2$
in greater detail.

\subsubsection{Elliptic curves}
If $A$ is an elliptic curve, that is, $g=1$, then
its Hodge polygon is given by the solid polygon

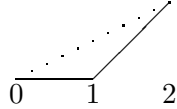
\begin{figure}[htbp]
  \setlength{\unitlength}{.2in}
  \centering
  \begin{picture}(4,2)
  \put(0,0.4){\line(1,0){2}}
  \put(2,0.4){\line(1,1){2}}
  \put(0,0.4){\qbezier[10](0,0)(2,1)(4,2)}
  \put(0,0){\makebox(0,0){$0$}}
  \put(2,0){\makebox(0,0){$1$}}
  \put(4,0){\makebox(0,0){$2$}}
  \end{picture}
\end{figure}

For the Newton polygon, there are two possibilities:
\begin{enumerate}
 \item The Newton polygon is equal to the Hodge polygon, and in this case, $A$ is called 
   {\em ordinary}.
   It follows from the results of 
   Examples \ref{example: crystalline cohomology abelian variety} that $A$ is
   ordinary if and only if $A[p](k)\iso\ZZ/p\ZZ$.
 \item The Newton polygon is equal to the dotted line, and in this case, $A$ is called {\em supersingular}.
   This case is equivalent to 
   $A[p](k)=\{0\}$.
\end{enumerate}
By a result of Deuring, there are roughly $p/12$ supersingular elliptic curves
over an algebraically closed field of positive characteristic $p$, whereas all the other ones 
are ordinary
(see Theorem  \ref{thm: ekedahl van der geer hodge} for a similar count for K3 surfaces).
In particular, a generic elliptic curve in positive characteristic is ordinary, that is, Newton
and Hodge polygon of $\Hcris{1}(A/W)$ coincide.
We refer the interested reader to \cite[Chapter IV.4]{Hartshorne} and \cite[Chapter V]{Silverman}
for more results, reformulations, and background information 
on ordinary and supersingular elliptic curves.

\subsubsection{Abelian surfaces} If $A$ is an Abelian surface, that is, $g=2$,
then its Hodge polygon is given by the solid polygon

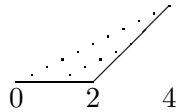
\begin{figure}[htbp]
  \setlength{\unitlength}{.1in}
  \centering
  \begin{picture}(8,4)
  \put(0,0.8){\line(1,0){4}}
  \put(4,0.8){\line(1,1){4}}
  \put(0,0.8){\qbezier[10](0,0)(4,2)(8,4)}
  \put(2,0.8){\qbezier[5](0,0)(2,1)(4,2)}
  \put(0,0){\makebox(0,0){$0$}}
  \put(4,0){\makebox(0,0){$2$}}
  \put(8,0){\makebox(0,0){$4$}}
  \end{picture}
\end{figure}

For the Newton polygon, there are now three possibilities:
\begin{enumerate}
  \item The Newton polygon is equal to the Hodge polygon, that is, $A$ is ordinary, or, equivalently,
   $A[p](k)\iso(\ZZ/p\ZZ)^2$.
 \item The Newton polygon has three slopes (lower dotted line), or, equivalently,
   $A[p](k)\iso\ZZ/p\ZZ$.
 \item The Newton polygon has only one slope (upper dotted line), that is, $A$ is supersingular, or, equivalently,
   $A[p](k)=\{0\}$.
\end{enumerate}
We refer to \cite[Exemples II.7.1]{Illusie deRham Witt} for details and further results.

\subsection{K3 surfaces}
Let $X$ be a K3 surface over $k$.
In Section \ref{subsec: k3 cohomology}, we computed the cohomology groups
of a K3 surface.
In particular, the only interesting crystalline cohomology group 
is $\Hcris{2}(X/W)$, which is free of rank $22$.
Moreover, in loc. cit. we also computed the Hodge numbers
$$
 \widetilde{h}_0\,:=\,h^{0,2}\,=\,1,\mbox{ \quad } \widetilde{h}_1\,:=\,h^{1,1}\,=\,20,
 \mbox{ \quad and \quad } \widetilde{h}_2\,:=\,h^{2,0}\,=\,1
$$
from which we obtain the geometric Hodge polygon.
%
%
In Section \ref{subsec: k3 cohomology}, we have also seen that
the crystalline cohomology groups of $X$ have no $p$-torsion,
and that the Fr\"olicher spectral sequence degenerates at $E_1$.
Thus, by Theorem \ref{thm:mazur nygaard ogus},
the geometric Hodge polygon of $X$ coincides with the 
Hodge polygon of the $F$-crystal $\Hcris{2}(X/W)$.

\begin{Exercise}
  \label{exercise: k3 crystals}
  For a K3 surface $X$, show that there are $12$ possibilities for the Newton polygon
  of the $F$-crystal $\Hcris{2}(X/W)$:
  \begin{enumerate}
    \item The Newton polygon has three slopes and multiplicities as follows:
       $$\begin{array}{lccc}
         \mbox{ slope } & 1-\frac{1}{h} & 1 & 1+\frac{1}{h} \\
         \mbox{ multiplicity } & h & 22-2h & h 
       \end{array}$$
       where $h$ is an integer with 
       $1\leq h\leq 11$.
       In case $h=1$, Hodge and Newton polygon
       coincide, and then, $X$ is called {\em ordinary}.
    \item The Newton polygon is of slope $1$ only
      (upper dotted line), and then, $X$ is called {\em supersingular}.
      In this case we set $h=\infty$.
  \end{enumerate}
  A discussion and details can be found in \cite[Section II.7.2]{Illusie deRham Witt}.
\end{Exercise}
  
Since $X$ is projective, there exists an ample 
line bundle ${\cal L}\in\Pic(X)$ and we will see in
Section \ref{sec: Tate conjecture} that
the $W$-module generated by the Chern class $c_1({\cal L})$
inside $\Hcris{2}(X/W)$ gives
rise to an $F$-crystal of slope $1$
(see also Exercise \ref{exercise: height and picard rank}).
In particular, this shows that the case $h=11$ 
in (1) of Exercise \ref{exercise: k3 crystals}
cannot occur.

In Section \ref{sec:Artin-Mazur}, we will define the formal Brauer group
of a K3 surface.
In Proposition \ref{prop: height determines newton}, we will see that 
the parameter $h$ from Exercise \ref{exercise: k3 crystals} can be interpreted as 
the height of the formal Brauer group.
And eventually, in Section \ref{sec:stratification}, we will see that $h$
gives rise to a stratification of the moduli spaces ${\cal M}_{2d,\FF_p}^\circ$ 
from Section \ref{sec: moduli spaces k3}, and that this stratification 
can also be interpreted in terms of $F$-zips (see Definition \ref{def: f-zip}).

\section{Supersingular K3 Surfaces}
\label{sec:supersingularK3}

For a complex K3 surface $X$, the de~Rham cohomology
group $\HdR{2}(X/\CC)$ comes with an integral and a real structure,
as well as two filtrations (Section \ref{subsec:complex geometry}).
Moreover, Poincar\'e duality equips $\HdR{2}(X/\CC)$ with a non-degenerate
bilinear form.
This linear algebra data is captured in the notion of a {\em polarized Hodge structure
of weight $2$}, and such data is parametrized by their {\em period domain}, which is
an open subset (with respect to the analytic topology) inside some Zariski-closed set
of some Gra\ss mannian.

In the previous section, we associated to a K3 surface $X$ over a perfect field $k$
of positive characteristic the $F$-crystal $(H,\varphi)$ arising from $\Hcris{2}(X/W)$.
This is a module over $W=W(k)$, which may be thought of as an integral structure.
By Theorems \ref{thm: mazur's hodge theorem} and \ref{thm:mazur nygaard ogus},
the two Hodge filtrations of $\HdR{2}(X/k)\iso H/pH$ are encoded in $(H,\varphi)$.
Poincar\'e duality, which also exists for crystalline cohomology 
(Section \ref{sec:crystalline cohomology}), equips $(H,\varphi)$
with a non-degenerate bilinear form.
This resulting structure $(H,\varphi,\langle-,-\rangle)$ is called a {\em K3 crystal}, 
and should be thought of as the characteristic-$p$ version of a
polarized Hodge structure of weight $2$ arising from a K3 surface.

Following Ogus \cite{Ogus}, we will only construct a period domain for
{\em supersingular} K3 crystals.
One crucial technical point is that this will be a Zariski closed subset of some
Gra\ss mannian and thus, projective over $k$.
In general, I would expect period domains for non-supersingular K3 crystals to be
open subsets of Zariski-closed sets of Gra\ss mannians, where open might also be
in the sense of Tate's rigid analytic spaces or Berkovich spaces.

\subsection{K3 crystals}
We start by introducing K3 crystals and their Tate modules, 
and shortly digress on the Tate conjecture.

\begin{Definition}[Ogus]
  \label{def: ogus k3 crystal}
  Let $k$ be a perfect field of positive characteristic $p$ 
  and let $W=W(k)$ be its Witt ring.
  A {\em K3 crystal of rank $n$}
  over $k$ is a free $W$-module $H$ of rank $n$
  together with a $\sigma$-linear injective map
  $\varphi:H\to H$ (that is, $(H,\varphi)$ is an $F$-crystal), 
  and a symmetric
  bilinear form
  $$
    \langle -,-\rangle\,:\,H\otimes_W H\,\to\, W
  $$
  such that
  \begin{enumerate}
    \item $p^2 H\subseteq {\rm im}(\varphi)$,
    \item $\varphi\otimes_W k$ has rank $1$,
    \item $\langle-,-\rangle$ is a perfect pairing,
    \item $\langle\varphi(x),\varphi(y)\rangle=p^2\sigma\langle x,y\rangle$.
  \end{enumerate}
  The K3 crystal is called {\em supersingular}, if moreover
  \begin{enumerate}
    \setcounter{enumi}{4}
    \item the $F$-crystal $(H,\varphi)$ is purely of slope $1$.
  \end{enumerate}
 \end{Definition}

\begin{Example}
  \label{example: supersingular k3 crystal}
   Let $X$ be a K3 surface over $k$.
   By Example \ref{example:geometry}, $H:=\Hcris{2}(X/W)$ 
   with Frobenius $\varphi$ is an $F$-crystal.
   By the results of Section \ref{subsec: k3 cohomology},
   it is of rank $22$.
   \begin{enumerate}
   \item By Exercise \ref{exc: slope restriction}, all slopes
     of $H$ lie in the interval $[0,2]$, which can also 
     be seen from the detailed classification
     in Exercise \ref{exercise: k3 crystals}.
     This implies that condition (1) of Definition \ref{def: ogus k3 crystal}
     is fulfilled.
   \item Poincar\'e duality equips $H$ with a symmetric 
     bilinear pairing $\langle-,-\rangle$, which satisfies conditions
     (3) and (4) of Definition \ref{def: ogus k3 crystal} by general
     properties of Poincar\'e duality.
   \item Since $X$ is a K3 surface, we have $\widetilde{h}_0:=h^2(\OO_X)=1$. 
    By Theorem \ref{thm:mazur nygaard ogus}, we have $h_0=1$
    for the Hodge polygon of $H$. 
    This implies that condition (2) of Definition \ref{def: ogus k3 crystal}
    holds true.
   \end{enumerate}
   Thus, $(H,\varphi,\langle-,-\rangle)$ is a K3 crystal of rank $22$.
   It is supersingular if and only if its Newton
   polygon is a straight line of slope one, which corresponds to case (2)
   in Exercise \ref{exercise: k3 crystals}, that is, $h=\infty$.
\end{Example}

\begin{Exercise}
  \label{exc: abelian surface crystal}
  Let $A$ be Abelian surface, that is, an Abelian variety of dimension $2$.
  Show that Frobenius and Poincar\'e duality turn
  $\Hcris{2}(A/W)$ into a K3 crystal of rank $6$. 
  We refer the interested reader to \cite[Section 6]{Ogus}, where crystals arising from 
  (supersingular) Abelian varieties are discussed in general. 
  For Abelian surfaces, it turns out that these crystals are closely related 
  to K3 crystals of rank $6$, see \cite[Proposition 6.9]{Ogus}.
\end{Exercise}

\subsection{The Tate module}
\label{sec: Tate conjecture}
If $X$ is a smooth and proper variety over $k$, then there exists
a crystalline Chern class map
$$
c_1\,:\,\Pic(X)\,\to\,\Hcris{2}(X/W)\,.
$$
Being a homomorphism of Abelian groups, $c_1$ satisfies
for all ${\cal L}\in\Pic(X)$
$$
   c_1( F^*({\cal L})) \,=\,  c_1({\cal L}^{\otimes p})\,=\, pc_1({\cal L}),
$$
where $F:X\to X$ denotes the absolute Frobenius morphism.
In particular, $c_1(\Pic(X))$  is contained in the Abelian subgroup 
(in fact, $\ZZ_p$-submodule)
of the $F$-crystal $(\Hcris{2}(X/W),\varphi)$ 
of those elements $x$ that satisfy $\varphi(x)=px$.
This observation motivates the following definition.

\begin{Definition}
   Let $(H,\varphi,\langle-,-\rangle)$ be a K3 crystal.
   Then, the {\em Tate module} of $H$ is defined to be
   the $\ZZ_p$-module
   $$
      T_H\,:=\, \{\,x\in H\,|\,\varphi(x)=px\,\}\,.
   $$   
\end{Definition}

Thus, by our computation above, 
we have $c_1(\NS(X))\subseteq T_H$, and it
is natural to ask whether this inclusion is in fact an equality, or, 
at least up to $p$-torsion.
If $X$ is defined over a finite field, 
this is the content of the {\em Tate conjecture}.

\begin{Conjecture}[Tate \cite{Tate Conjecture}]
  \label{conj: Tate conjecture}
  Let $X$ be a smooth and proper surface over a finite field $\FF_q$
  of characteristic $p$.
  Then, the following statements hold true:
  \begin{enumerate}
   \item The first Chern class induces an isomorphism
      $$
         c_1\,:\,\NS(X)\otimes_{\ZZ}\QQ_p \,\stackrel{\iso}{\longrightarrow}\, T_H\otimes_{\ZZ_p}\QQ_p\,.
      $$
   \item For every prime $\ell\neq p$, the first Chern class induces an isomorphism
      $$
         c_1\,:\,\NS(X)\otimes_{\ZZ}\QQ_\ell \,\stackrel{\iso}{\longrightarrow}\, 
         \Het{2}\left(X\times_{\FF_q}\overline{\FF}_q,\,\QQ_\ell(1)\right)^{{\rm Gal}(\overline{\FF}_q/\FF_q)},
      $$
      where the right hand side denotes invariants under the Galois action.
   \item The rank of $\NS(X)$ is the pole order of the zeta function $Z(X/\FF_q,T)$ at $T=q^{-1}$.
  \end{enumerate}
\end{Conjecture}

The equivalences of (1), (2), and (3) follow from the Weil conjectures, more precisely, 
from the Riemann hypothesis, which relates the zeta function to $\ell$-adic and crystalline cohomology, 
see \cite[Appendix C]{Hartshorne} and \cite[Section 9.10]{Liedtke overview}.
In \cite[Theorem 4]{Tate Endomorphisms}, Tate proved this conjecture
for Abelian varieties, as well as for products of curves.
For K3 surfaces, it was established in several steps depending
on the slopes of the $F$-crystal $\Hcris{2}$ - in terms of the notations of
Exercise \ref{exercise: k3 crystals}:
for $h=1$ by Nygaard \cite{Nygaard}, 
for $h<\infty$ by Ogus and Nygaard \cite{Ogus; Nygaard}, and 
in general by Charles \cite{Charles}, Madapusi Pera \cite{Madapusi Pera}, and
Maulik \cite{Maulik}.

\begin{Theorem}[Nygaard, Nygaard--Ogus, Charles, Madapusi-Pera, Maulik]
   \label{thm: Tate conjecture}
  Tate's conjecture holds for K3 surfaces over finite fields of 
   odd characteristic.
\end{Theorem}

Let us mention the following, somewhat curious corollary: namely, Swinnerton-Dyer
observed (see \cite{Artin Supersingular}) that Tate's conjecture for K3 surfaces implies that
the N\'eron--Severi rank of a K3 surface over $\overline{\FF}_p$ is {\em even}.
This was used in \cite{BHT} and \cite{LL} to show that every K3 surface 
of {\em odd} N\'eron--Severi rank contains infinitely many rational curves, and
we refer to \cite{Liedtke overview} and \cite{benoist} for an overview.

\subsection{Supersingular K3 surfaces}
\label{subsec: supersingular K3s}
Let us now discuss supersingular K3 crystals in greater detail, that is, 
K3 crystals that are of slope $1$ only.
It turns out that they are largely determined by their Tate-modules.
In case a supersingular K3 crystal arises as $\Hcris{2}$ of a K3 surface,
the Tate conjecture predicts that the surface has
Picard rank $22$, that is, the K3 surface is {\em Shioda-supersingular}.

First, let us recall a couple of facts on quadratic forms
and their classification, and we refer to \cite[Chapitre IV]{Serre Course}
for details and proofs:
let $R$ be a ring and $\Lambda$ a free $R$-module of finite rank together
with a symmetric bilinear form 
$$
\langle-,-\rangle \,:\Lambda\otimes_R\Lambda\,\to\,\Lambda\,.
$$
We choose a basis $\{e_1,...,e_n\}$ of $\Lambda$,
form the matrix $G:=(g_{ij}:=\langle e_i,e_j\rangle)_{i,j}$,
and define its {\em discriminant}  to be $\det(G)$. 
A different choice of basis of $\Lambda$ changes it by 
an element of $R^{\times2}$, and thus, the class $d(\Lambda)$
of $\det(G)$ in $R/(R^{\times2})$ does not depend on the choice of basis.
The discriminant is zero if and only if
the form is {\em degenerate}, that is, if there exists a
$0\neq v\in\Lambda$ such that $\langle v,w\rangle=0$ 
for all $w\in\Lambda$.
Next, we let $\Lambda^\vee:={\rm Hom}_R(\Lambda, R)$
be the dual $R$-module.
Via $v\mapsto \langle v,-\rangle$, we obtain a natural
map $\Lambda\to\Lambda^\vee$, which is injective
if and only if the form is non-degenerate.
In case this map is an isomorphism, which is the case if and only
if the discriminant is a unit, the form is called {\em perfect}.

Let us now assume that $R$ is a DVR, say, with valuation
$\nu$.
The example we have in mind is
the ring of $p$-adic integers $\ZZ_p$ or, more generally,
the ring $W(k)$ of Witt vectors of a perfect field $k$
of positive characteristic $p$, together
with its $p$-adic valuation ${\rm ord}_p$.
Then, since units have valuation zero, 
$$
  {\rm ord}_p(\Lambda)\,:=\,\nu\left(d(\Lambda)\right)
$$ 
is a well-defined integer, and
the form is perfect if and only if ${\rm ord}_p(\Lambda)=0$.

Finally, we note that quadratic forms over $\QQ_p$ are classified
by their rank, their discriminant, and their so-called Hasse
invariant, see \cite[Chapitre IV]{Serre Course} for proofs and 
details.
These results are the key to the following classification
of Tate modules of supersingular K3 crystals
from \cite{Ogus}.

\begin{Proposition}[Ogus]
  \label{prop: supersingular k3 crystal results}
  Let $(H,\varphi,\langle-,-\rangle)$ be a supersingular K3 crystal
  and let $T_H$ be its Tate module.
  Then,
  $$
      {\rm rank}_W\, H\,=\,{\rm rank}_{\ZZ_p}\,T_H
  $$
  and the bilinear form $(H,\langle-,-\rangle)$ restricted
  to $T_H$ induces a non-degenerate form
  $$
      T_H \,\otimes_{\ZZ_p}\, T_H\,\to\,\ZZ_p,
  $$
  which is not perfect.
  More precisely,
  \begin{enumerate}
    \item ${\rm ord}_p(T_H)=2\sigma_0>0$ for some integer $\sigma_0$,
      called the {\em Artin invariant}.
    \item $(T_H,\langle-,-\rangle)$ is determined up to isometry by $\sigma_0$.
    \item ${\rm rank}_W\, H\geq2\sigma_0$.
    \item There exists an orthogonal decomposition
     $$
        (T_H,\langle-,-\rangle) \,\iso\, (T_0,p\langle-,-\rangle) \perp (T_1,\langle-,-\rangle)
     $$
     where $T_0$ and $T_1$ are $\ZZ_p$-lattices, whose bilinear forms are perfect,
     and of ranks ${\rm rank}\,T_0=2\sigma_0$ and ${\rm rank}\,T_1={\rm rank}_W\,H-2\sigma_0$.
  \end{enumerate}
\end{Proposition}

Combining this proposition with the Tate conjecture (Theorem \ref{thm: Tate conjecture}),
we obtain a characterization of those K3 surfaces whose associated K3 crystal
is supersingular. 
Namely, let us recall from Section \ref{subsec: k3 cohomology} that the second crystalline
cohomology group of a K3 surface is a free $W$-module of rank $22$.
Using that the first crystalline Chern map is injective, this shows that the rank of the N\'eron--Severi
group of a K3 surface can be at most $22$.
This said, we have the following result.

\begin{Theorem}
 \label{thm: two supersingular k3 notions}
  Let $X$ be a K3 surface over an algebraically closed field of odd characteristic.
  Then, the following are equivalent
  \begin{enumerate}
    \item The K3 crystal  $\Hcris{2}(X/W)$ is supersingular.
    \item The N\'eron--Severi group $\NS(X)$ has rank $22$.
  \end{enumerate}
\end{Theorem}

\prf
If $\NS(X)$ has rank $22$, then $c_1(\NS(X))\otimes_\ZZ W$ is a 
sub-$F$-crystal of $\Hcris{2}(X/W)$ of slope $1$, thereby
establishing $(2)\Rightarrow(1)$.
Conversely, assume that $H:=\Hcris{2}(X/W)$ is a supersingular $F$-crystal.
By Proposition \ref{prop: supersingular k3 crystal results}, the Tate module
of $H$ has rank $22$.
If $X$ can be defined over a finite field of odd characteristic 
then Theorem \ref{thm: Tate conjecture} implies that $\NS(X)$ is of rank $22$.
If $X$ is not definable over a finite field, 
then there exists some variety $B$ over some finite field $\FF_q$, such that
$X$ is definable over the function field of $B$.
Spreading out $X$ over $B$ and passing to an open and dense subset of $B$ if necessary, we may
assume that $X$ is the generic fiber of a smooth and projective family ${\cal X}\to B$
of K3 surfaces over $\FF_q$.
Since $H$ is supersingular, all fibers in this family 
also have supersingular $\Hcris{2}$ by \cite[Section 1]{Artin Supersingular},
and since the N\'eron--Severi rank is constant in families of K3 surfaces
with supersingular $\Hcris{2}$ by
\cite[Theorem 1.1]{Artin Supersingular}, this establishes the converse direction 
$(1)\Rightarrow(2)$.
\qed\medskip

\begin{Remark}
  K3 surfaces satisfying (1) are called {\em Artin-supersingular}, see \cite{Artin Supersingular}, 
  where it is formulated in terms of formal Brauer groups, a point of view that we will discuss 
  in Section \ref{sec:formal groups} below.
  K3 surfaces satisfying (2) are called {\em Shioda-supersingular}, see \cite{Shioda 1974}.
  In view of the theorem, 
  a K3 surface in odd characteristic 
  satisfying (1) or (2) is simply called {\em supersingular}.
\end{Remark}

\begin{Examples}
  \label{examples: supersingular K3}
  Let us give examples of supersingular K3 surfaces.
  \begin{enumerate}
   \item Let $A$ be a supersingular Abelian surface in odd characteristic
      (see Section \ref{subsec: Abelian varieties crystals}).
      Then, the Kummer surface
      $X:={\rm Km}(A)$ of $A$ is a supersingular K3 surface.
      Let $\sigma_0$ be the Artin invariant of the Tate module $T_H$
      of the supersingular K3 crystal $H:=\Hcris{2}(X/W)$.
      Then,
      $$
      \sigma_0(T_H) \,=\,\left\{
      \begin{array}{cl}
      1 & \mbox{ if $A=E\times E$, where $E$ is a supersingular}\\
      & \mbox{ elliptic curve (the isomorphism class of $A$}\\
      & \mbox{ does not depend on the choice of $E$),}\\
      2 & \mbox{ else, }
      \end{array}\right.
      $$
      see \cite[Theorem 7.1 and Corollary 7.14]{Ogus}, or
      \cite[Proposition 3.7 and Theorem 4.3]{Shioda supersingular}.
      Conversely, by loc. cit, every supersingular K3 surface 
      in odd characteristic with $\sigma_0\leq2$
      is the Kummer surface of a supersingular
      Abelian surface.
    \item The Fermat quartic
    $$
      X_4 \,:=\, \{  x_0^4+x_1^4+x_2^4+x_3^4\,=\,0 \}\,\subset\,\PP^3_k
    $$
    defines a K3 surface in characteristic $p\neq2$
    and it is supersingular if and only if 
    $p\equiv3\mod 4$ by \cite[Corollary to Proposition 1]{Shioda 1974}.
    Moreover, if $X_4$ is supersingular, then it has $\sigma_0(T_H)=1$
    by \cite[Example 5.2]{Shioda supersingular}, and thus, it is 
    a Kummer surface by the previous example.
  \end{enumerate}
 We note that supersingular Kummer surfaces form a $1$-dimensional
 family, whereas all supersingular K3 surfaces form a $9$-dimensional
 family - we refer to Section \ref{sec:torelli} for moduli spaces.
\end{Examples}

In view of Theorem \ref{thm: two supersingular k3 notions}, we now
identify the N\'eron--Severi lattices arising from supersingular K3 surfaces
abstractly, and classify them in terms of discriminants, which
gives rise to the Artin invariant of such a lattice.

\begin{Definition}
  \label{def: k3 lattice}
  A {\em supersingular K3 lattice} is a free Abelian group $N$ of rank $22$
  with an even symmetric bilinear form $\langle-,-\rangle$ with the following
  properties
  \begin{enumerate}
    \item The discriminant $d(N\otimes_\ZZ\QQ)$ is $-1$ in $\QQ^*/\QQ^{*2}$.
    \item The signature of $(N\otimes_\ZZ\RR)$ is $(1,21)$.
    \item The cokernel of $N\to N^\vee$ is annihilated by $p$.
  \end{enumerate}
\end{Definition}
  
We note that we follow \cite[Definition 1.6]{Ogus PM}, which is slightly 
different from \cite[Definition 3.17]{Ogus} (in the latter article, it 
is stated for $\ZZ_p$-modules rather than $\ZZ$-modules). 
Let us shortly collect the facts:
by \cite[(1.6)]{Ogus PM} and the references given there, the N\'eron--Severi
lattice of a supersingular K3 surface is a supersingular K3 lattice in the
sense of Definition \ref{def: k3 lattice}
(for example, condition (2) follows from the Hodge index theorem).
Next, if $N$ is a supersingular K3 lattice, then its discriminant
$d(N)$, which is an integer,  is equal to $-p^{2\sigma_0}$ for some integer $1\leq\sigma_0\leq10$.

\begin{Definition}
  \label{def: artin invariant lattice}.
  The integer $\sigma_0$ associated to a supersingular K3 lattice is called
  the {\em Artin invariant} of the lattice.
  If $X$ is a supersingular K3 surface, we define its Artin invariant
  to be the Artin invariant of its N\'eron--Severi lattice.
\end{Definition}

This invariant was introduced in \cite{Artin Supersingular}, and 
an important result is the following theorem,
see \cite[Section 1]{Rudakov Shafarevich} and \cite[Section 3]{Ogus}.

\begin{Theorem}[Rudakov--Shafarevich]
 The Artin invariant determines a supersingular K3 lattice up to isometry.
\end{Theorem}

We refer the interested reader to \cite[Section 1]{RS supersingular}  for explicit descriptions 
of these lattices, which do exist for all values $1\leq\sigma_0\leq10$.

Before proceeding, let us shortly digress on quadratic forms over finite fields: 
let $V$ be a $2n$-dimensional vector space over a finite field $\FF_q$ of
odd characteristic.
Let $\langle-,-\rangle:V\times V\to\FF_q$ be a non-degenerate quadratic form.
Two-dimensional examples are the hyperbolic plane $U$,
as well as $\FF_{q^2}$ with the quadratic form arising from the norm.
By the classification of quadratic forms over finite fields, 
$V$ is isometric to $nU$ or to $(n-1)U\perp \FF_{q^2}$.
The form $\langle-,-\rangle$ is called {\em non-neutral} if there 
exists no $n$-dimensional isotropic subspace inside $V$.
By the classification result just mentioned, there is precisely
one non-neutral quadratic space of dimension $2n$ over
$\FF_q$, namely, $(n-1)U\perp\FF_{q^2}$.

Next, for a supersingular K3 lattice $(N,\langle-,-\rangle)$,
we set $N_1:=N/pN^\vee$.
Then, $N_1$ is a $(22-2\sigma_0)$-dimensional $\FF_p$-vector space
and $\langle-,-\rangle$ induces a quadratic form on $N_1$, which is
non-degenerate and non-neutral.
The form $\langle-,-\rangle$ on $pN^\vee\subseteq N$ is divisible by $p$
and dividing it by $p$ we obtain a non-degenerate and non-neutral bilinear 
form on the  $2\sigma_0$-dimensional $\FF_p$-vector space $N_0:=pN^\vee/ pN$.
We refer to \cite[(1.6)]{Ogus PM} for details.
In Section \ref{subsec: char subspace} below, we will use these $\FF_p$-vector spaces
to classify supersingular K3 crystals explicitly as well as to construct their moduli
spaces - the point is that it is easier to deal with $\FF_p$-vector spaces rather 
than $\ZZ$- or $\ZZ_p$- lattices.

Finally, for a supersingular K3 lattice $N$, we set $\Gamma:=N\otimes_\ZZ\ZZ_p$ and denote
the induced bilinear form on $\Gamma$ again by $\langle-,-\rangle$.
Then, we have ${\rm ord}_p(\Gamma)=2\sigma_0$.
By \cite[Lemma 3.15]{Ogus}, non-neutrality of the form induced on $N_0$ is equivalent
to the Hasse invariant of $\Gamma$ being equal to $-1$.
Moreover, since the cokernel of $N\to N^\vee$ is annihilated by $p$,
the same is true for $\Gamma\to\Gamma^\vee$, and thus, 
by \cite[Lemma 3.14]{Ogus}, we obtain an orthogonal
decomposition 
$$
 (\Gamma,\langle-,-\rangle)\,\iso\, (\Gamma_0,p\langle-,-\rangle)\,\perp\,(\Gamma_1,\langle-,-\rangle),
$$
where $\Gamma_0$ and $\Gamma_1$ are perfect $\ZZ_p$-lattices
of ranks $2\sigma_0$ and $22-2\sigma_0$, respectively.
In particular, $\Gamma$ satisfies the conditions of
a supersingular K3 lattice over $\ZZ_p$ as defined in \cite[Definition 3.17]{Ogus}.
We refer to \cite[Corollary 3.18]{Ogus} for details about the classification of 
supersingular K3 lattices over $\ZZ_p$ up to isogeny and up to isomorphism.

\subsection{Characteristic subspaces}
\label{subsec: char subspace}
In order to classify supersingular K3 crystals, we now describe them in terms 
of so-called characteristic subspaces, and then, classify these latter ones.
For a supersingular K3 surface, this characteristic subspace
arises from the kernel of the de~Rham Chern class 
$c_1:\NS(X)\to\HdR{2}(X/k)$.
(Note that in characteristic zero, $c_1$ is injective modulo torsion.)
These considerations stress yet again the close relation between 
crystals and de~Rham cohomology.

\begin{Definition}
 \label{def:non-neutral}
  Let $\sigma_0\geq1$ be an integer,
  let $V$ be a $2\sigma_0$-dimensional $\FF_p$-vector space 
  with $p\neq2$, and let
  $$
    \langle-,-\rangle\,:\,V\times V\,\to\,\FF_p
  $$
  be a non-degenerate and non-neutral quadratic form.
  Next, let $k$ be a perfect field of
  characteristic $p$ and set 
  $\varphi:={\rm id}_V\otimes F_k:V\otimes_{\FF_p} k\to V\otimes_{\FF_p}k$.
  A subspace $K\subset V\otimes_{\FF_p}k$ is
  called {\em characteristic} if
  \begin{enumerate}
    \item $K$ is totally isotropic of dimension $\sigma_0$, and
    \item $K+\varphi(K)$ is of dimension $\sigma_0+1$.
  \end{enumerate} 
  Moreover, a characteristic subspace $K$ is {\em strictly characteristic}
  if moreover
  \begin{enumerate}
     \setcounter{enumi}{2}
     \item
     $$
        V\otimes_{\FF_p}k\,=\,\sum_{i=0}^\infty \varphi^i(K)
     $$
     holds true.
  \end{enumerate}
\end{Definition}

For a perfect field $k$ of odd characteristic, we define 
the categories
$$
    {\rm K3}\,(k) \,:=\,
    \left\{
    \begin{array}{l}
      \mbox{ category of supersingular K3 crystals }\\
      \mbox{ with only isomorphisms as morphisms }
    \end{array} \right\} \,.
$$
and
$$
    \CC3\,(k) \,:=\,
    \left\{
    \begin{array}{l}
      \mbox{ category of pairs $(T,K)$, where $T$ is a supersingular }\\
      \quad \mbox{ K3 lattice over $\ZZ_p$, and where $K\subset T_0\otimes_{\ZZ_p}k$ is a }\\
      \quad \mbox{ strictly characteristic subspace, } \\
      \mbox{ with only isomorphisms as morphisms }
    \end{array} \right\} \,.
$$ 
Finally, we define $\CC3(k)_{\sigma_0}$ to be the subcategory of $\CC3(k)$,
whose characteristic subspaces are $\sigma_0$-dimensional.
 
\begin{Theorem}[Ogus]
  \label{thm: K3=C3}
  Let $k$ be an algebraically closed field of odd characteristic.
  Then, the assignment
  $$
  \begin{array}{ccc}
    {\rm K3}\,(k) &\stackrel{\gamma}{\longrightarrow}& \CC3\,(k) \\
    (H,\varphi,\langle-,-\rangle) &\mapsto&
    \left( T_H,\, \ker (T_H\otimes_{\ZZ_p}k\to H\otimes_{\ZZ_p}k)\subset T_0\otimes_{\ZZ_p}k 
    \right)
  \end{array}
  $$
  defines an equivalence of categories.
\end{Theorem}

The map $\gamma$ from the theorem has the following geometric origin, and we refer to
\cite[Section 2]{Ogus PM} for details:
if $X$ is a supersingular K3 surface, then
$H:=\Hcris{2}(X/W)$ is a supersingular K3 crystal.
Moreover, the Tate module $T_H$ is a supersingular K3 lattice over $\ZZ_p$, 
the first Chern class $c_1$ identifies $\NS(X)\otimes_\ZZ \ZZ_p$ with $T_H$, and 
the characteristic subspace associated to $H$ arises from the kernel of
$c_1:\NS(X)\otimes_\ZZ k\to \HdR{2}(X/k)$.
\medskip

Now, we describe and classify characteristic subspaces over an algebraically
closed field $k$ of odd characteristic $p$ explicitly, and we refer 
to  \cite[p. 33-34]{Ogus} for technical details:
let $V$ be a $2\sigma_0$-dimensional $\FF_p$-vector space with
a non-neutral form $\langle-,-\rangle$, let
$\varphi={\rm id}\otimes F_k:V\otimes_{\FF_p} k\to V\otimes_{\FF_p} k$,
and let $K\subset V\otimes_{\FF_p}k$ be a strictly characteristic subspace.
Then, 
$$
\ell_K\,:=\,K\cap\varphi(K)\,\cap\,...\,\cap\,\varphi^{\sigma_0-1}(K)
$$
is a line inside $V\otimes_{\FF_p}k$.
We choose a basis element $0\neq e\in\ell_K$ and set
$$
  e_i\,:=\,\varphi^{i-1}(e)\mbox{ \quad for \quad }i=1,...,2\sigma_0.
$$
Then, the $\{e_i\}$ form a basis of $V\otimes_{\FF_p}k$.
We have $\langle e,e_{\sigma_0+1}\rangle\neq0$, and changing
$e$ by a scalar if necessary (here, we use that $k$ is algebraically closed),
we may assume $\langle e,e_{\sigma_0+1}\rangle=1$.
We note that this normalization makes $e$ unique up to a $(p^{\sigma_0}+1)$.th root of unity.
Then, we define
$$
  a_i\,:=\,a_i(e,V,K)\,:=\, \langle e, e_{\sigma_0+1+i}\rangle\mbox{ \quad for \quad }
  i=1,...,\sigma_0-1\,.
$$
If $\zeta$ is a $(p^{\sigma_0}+1)$.th root of unity, then, replacing
$e$ by $\zeta e$, transforms the $a_i$ as $a_i\mapsto \zeta^{1-p^i}a_i$.
This said, we denote by $\mu_n$ the group scheme of
$n$.th roots of unity, and then, we have the following 
classification result.

\begin{Theorem}[Ogus]
  \label{thm: explicit coordinates}
  Let $k$ be an algebraically closed field of odd characteristic.
  Then, there exists a bijection
  $$
  \begin{array}{ccc}
     \CC3\,(k)_{\sigma_0} &\to& \Aff_k^{\sigma_0-1}(k)/\mu_{p^{\sigma_0}+1}(k) \\
     K &\mapsto& (a_1,...,a_{\sigma_0-1})
  \end{array}
  $$
  where the $a_i:=a_i(e,V,K)$ are as defined above.
\end{Theorem}

Having described characteristic subspaces over algebraically closed fields,
we now study them in families.

\begin{Definition}
  \label{def: generatrix}
  Let $(V,\langle-,-\rangle)$ be a $2\sigma_0$-dimensional $\FF_p$-vector space
  with a non-neutral quadratic form.
  If $A$ is an $\FF_p$-algebra, a {\em geneatrix} of $V\otimes_{\FF_p}A$
  is a direct summand
  $$
     K\subset V\otimes_{\FF_p}A
  $$
  of rank $\sigma_0$ such that $\langle-,-\rangle$ restricted to $K$ vanishes
  identically. 
  We define the set of geneatrices
  $$
    \underline{\rm Gen}_V(A)\,:=\,
    \left\{ \mbox{ generatrices of $V\otimes_{\FF_p}A$ } \right\}
  $$
  as well as 
  $$
    \underline{\rm M}_V(A)\,:=\,
    \left\{ K\in\underline{{\rm Gen}}_V(A), K+F_A^*(K) \mbox{ is a direct summand of rank $\sigma_0+1$ } \right\},
  $$
  which is the set of characteristic generatrices. 
\end{Definition}

\begin{Proposition}[Ogus]
  \label{prop: moduli of generatrices}
  The functor from $\FF_p$-algebras to sets given by
  $$A\,\mapsto\,\underline{{\rm M}}_V(A)$$
  is representable
  by a scheme $M_V$, which is smooth, projective, and of dimension
  $\sigma_0-1$ over $\FF_p$.
\end{Proposition}

Let $N$ be a supersingular K3 lattice with Artin invariant $\sigma_0$.
At the end of Section \ref{subsec: supersingular K3s} we set $N_0:=pN^\vee/pN$ 
and noted that it is a $2\sigma_0$-dimensional  $\FF_p$-vector space that 
inherits a non-degenerate and non-neutral bilinear form from $N$.
We set ${\cal M}_N:=M_{N_0}$.

\begin{Definition}
   \label{def: mN moduli space}
   ${\cal M}_N$ is called the
   {\em moduli space of $N$-rigidified K3 crystals}.
\end{Definition}

\begin{Examples}
  \label{examples: crystalline moduli}
  If $V$ is $2\sigma_0$-dimensional, then
  \begin{enumerate}
    \item If $\sigma_0=1$, then $M_V\iso\Spec\FF_{p^2}$.
    \item If $\sigma_0=2$, then $M_V\iso\PP^1_{\FF_{p^2}}$.
    \item If $\sigma_0=3$, then $M_V$ is isomorphic to the Fermat surface
      of degree $p+1$ in $\PP^3_{\FF_{p^2}}$.
  \end{enumerate}
  We refer to \cite[Examples 4.7]{Ogus} for details, as well as to
  Theorem  \ref{thm: moduli fibration} for a generalization
  to higher dimensional $V$'s.
\end{Examples}

Anticipating the crystalline Torelli theorem in Section \ref{sec:torelli},
let us comment on the $\sigma_0=1$-case and give a geometric interpretation:
then, we have $M_V\iso\Spec\FF_{p^2}$.
By Theorem \ref{thm: crystalline torelli theorem} or 
Examples \ref{examples: supersingular K3}, there exists precisely one supersingular 
K3 surface with $\sigma_0=1$ up to isomorphism over algebraically 
closed fields of odd characteristic.
More precisely, this surface is the Kummer surface
${\rm Km}(E\times E)$, where $E$ is a supersingular elliptic curve.
Although this surface can be defined over $\FF_p$,
there is no model $X$ over $\FF_p$ such
that all classes of $\NS(X_{\overline{\FF}_p})$ are already 
defined over $\FF_p$.
Models with full N\'eron--Severi group do
exist over $\FF_{p^2}$ - but then, there is a non-trivial
Galois-action of ${\rm Gal}(\FF_{p^2}/\FF_p)$ on 
$\NS(X_{\FF_{p^2}})$.
This explains (via the crystalline Torelli theorem) the Galois
action on $M_V$, as well as the fact that
$M_V\times_{\FF_p}\overline{\FF}_p$ consists of two points,
whereas it corresponds to only one surface.

\section{Ogus' Crystalline Torelli Theorem}
\label{sec:torelli}

We now come to the period map and the crystalline Torelli theorem
for supersingular K3 surfaces.
To state it, we fix a prime $p\geq5$ and a
supersingular K3 lattice $N$ as in Definition \ref{def: k3 lattice}.
Then, there exists a moduli space ${\cal S}_N$ of $N$-marked supersingular K3 surfaces,
which is a scheme that is locally of finite type, almost proper, and smooth of dimension
$\sigma_0(N)-1$ over $\FF_p$.

Associating to an $N$-marked supersingular K3 surface $X$ 
the $F$-zip associated to $\HdR{2}(X/k)$ yields a morphism 
$$
  \pi^{\rm mod\, p}_N\,:\,{\cal S}_N\,\to\,{\cal F}^\tau
$$
(notation as at the end of Section \ref{subsec:algebraic deRham}).
However, ${\cal F}^\tau$ is a rather discrete object,
it essentially only remembers the Artin invariant $\sigma_0(X)$,
and thus, $\pi_N^{\rm mod\, p}$ is more of a ``mod $p$ shadow''
of the saught-after period map.

Associating to an $N$-marked supersingular K3 surface the
$N$-rigidified K3 crystal associated to $\Hcris{2}(X/W)$ yields a morphism
$$
\pi_N\,:\,{\cal S}_N\,\to\,{\cal M}_N
$$
(see Definition \ref{def: mN moduli space}),
which is locally of finite type, \'etale, and surjective, but it is not an
isomorphism.

If we equip $N$-rigidified K3 crystals with ample cones,
we obtain a new moduli space ${\cal P}_N$ that comes
with a forgetful morphism ${\cal P}_N\to {\cal M}_N$, and then, 
$\pi_N$ lifts to a morphism
$$
  \widetilde{\pi}_N\,:\, {\cal S}_N \,\to\, {\cal P}_N,
$$ 
the {\em period map}.
By {\em Ogus' crystalline Torelli theorem} for supersingular
K3 surfaces, the period map is an isomorphism.

\subsection{Moduli of marked supersingular K3 surfaces}
Let $N$ be a supersingular K3 lattice in characteristic $p$
as in Definition \ref{def: k3 lattice}, and let
$\sigma_0$ be its Artin invariant as in Definition \ref{def: artin invariant lattice}.
For an algebraic space $S$ over $\FF_p$, we denote by $\underline{N}_S$ the 
constant group algebraic space defined by $N$ over $S$.
Then, we consider the functor $\underline{{\cal S}}_N$
of {\em $N$-marked K3 surfaces}
$$
\begin{array}{ccc}
 \left(\begin{array}{c}\mbox{Algebraic spaces}\\ \mbox{over $\FF_p$}\end{array}\right)
  &\to&\mbox{(Sets)}\\
  S &\mapsto&
  \left\{
  \begin{array}{l}
    \mbox{ smooth and proper morphisms $f:X\to S$ }\\
    \mbox{ of algebraic spaces, each of whose geometric }\\
    \mbox{ fibers  is a K3 surface, together with a }\\
    \mbox{ morphism of group spaces $\underline{N}_S\to\Pic_{X/S}$ }\\
    \mbox{ compatible with intersection forms. }
  \end{array}
  \right\}
\end{array}
$$
Since an $N$-marked K3 surface $X\to S$ has no 
non-trivial automorphisms \cite[Lemma 2.2]{Ogus PM}, 
it is technical, yet straight forward to prove that 
this functor can be represented by an algebraic space
\cite[Theorem 2.7]{Ogus PM}.
It follows a posteriori from the crystalline Torelli
theorem \cite[Theorem III']{Ogus PM} that this algebraic space
is a scheme.

\begin{Theorem}[Ogus]
  The functor $\underline{{\cal S}}_N$ is represented by a
  scheme ${\cal S}_N$, which is locally of finite type, 
  almost proper, and smooth of 
  dimension $\sigma_0(N)-1$ over $\FF_p$.
\end{Theorem}

Here, a scheme is called {\em almost proper} if satisfies
the surjectivity part of the valuative criterion with DVR's as test rings.
However, this moduli space is {\em not} separated.
This non-separatedness arises from elementary transformations,
which is analogous to degenerations of complex K\"ahler K3 surfaces.
We refer to \cite[Section 7]{Burns Rapoport}, \cite{Morrison}, and
\cite[page 380]{Ogus PM} for details and further discussion.

\subsection{The period map}
Let $(H,\varphi)$ be the K3 crystal associated to $\Hcris{2}(X/W)$ of
a K3 surface $X$ over a perfect field $k$ of
positive characteristic $p$.
Moreover, if $X$ is $N$-marked, then the inclusion
$N\to\NS(X)$ composed with the first crystalline Chern map
yields a map $N\to T_H$, where $T_H$ denotes the
Tate module of the K3 crystal.
Thus, an $N$-marked K3 surface gives rise to
an $N$-rigidified supersingular K3 crystal, 
which gives rise to a morphism of schemes
$$
  \pi\,:\,{\cal S}_N\,\to\,{\cal M}_N
$$
(we refer to \cite[Section 5]{Ogus} for families of crystals).
Although this morphism is \'etale and surjective by
\cite[Proposition 1.16]{Ogus PM}, it is not an isomorphism.
In order to obtain an isomorphism (the period map), we have to enlarge ${\cal M}_N$
by considering 
$N$-marked supersingular K3 crystals together with {\em ample cones}.

\begin{Definition}
  Let $N$ be a supersingular K3 lattice.
  Then, we define its {\em roots} to be the set
  $$
     \Delta_N\,:=\,\{ \delta\in N\,|\,\delta^2=-2 \}\,.
  $$
  For a root $\delta\in\Delta_N$, we define the {\em reflection} in $\delta$
  to be the automorphism of $N$ defined by
  $$
    r_\delta \,:\, x\,\mapsto\, x+\langle x,\delta\rangle\cdot\delta\mbox{ \quad for all }x\in N.
  $$
  We denote by $R_N$ the subgroup of ${\rm Aut}(N)$ generated by all 
  $r_\delta$, $\delta\in\Delta_N$.
  We denote by $\pm R_N$ the subgroup of ${\rm Aut}(N)$ generated by
  $R_N$ and $\pm{\rm id}$.
  Finally, we define the set
  $$
   V_N\,:=\,\left\{ x\in N\otimes\RR\,|\, x^2>0\mbox{ and }\langle x,\delta\rangle\neq0
   \mbox{ for all }\delta\in\Delta_N \right\}
   \,\subset\,
   N\otimes\RR.
  $$
\end{Definition}

Then, the subset $V_N\subset N\otimes\RR$ is open, and each of its connected
components meets $N$.
A connected component of $V_N$ is called an {\em ample cone},
and we denote by $C_N$ the set of ample cones.
Moreover, the group $\pm R_N$ operates simply and transitively on $C_N$.
We refer to \cite[Proposition 1.10]{Ogus PM} for details and proof.

\begin{Definition}
  Let $N$ be a supersingular K3 lattice,
  and let $S$ be an algebraic space over $\FF_p$.
  For a characteristic geneatrix $K\in{\cal M}_N(S)$, that is, a
  local direct factor $K\subset\OO_S\otimes N_0$
  as in Definition \ref{def: generatrix}, 
  we set for each point $s\in S$
  $$
  \begin{array}{ccl}
    \Lambda(s) &:=& N_0\cap K(s)\\
    N(s) &:=& \{ x\in N\otimes\QQ\,|\, px\in N\mbox{ and }\overline{px}\in \Lambda(s)\}\\
    \Delta(s)&:=& \{\delta\in N(s)\,|\,\delta^2=-2\}
  \end{array}
  $$
  An {\em ample cone} for $K$ is an element
  $$
     \alpha\,\in\,\prod_{s\in S}C_{N(s)}
  $$
  such that $\alpha(s)\subseteq\alpha(t)$ whenever $s$ is a specialization of $t$.
\end{Definition}

Having introduced these definitions, we consider the functor $\underline{{\cal P}}_N$
$$
\begin{array}{ccc}
 \left(\begin{array}{c}\mbox{Algebraic spaces}\\ \mbox{over $\FF_p$}\end{array}\right)
  &\to&\mbox{(Sets)}\\
  S&\mapsto&
  \left\{
  \begin{array}{l}
    \mbox{ characteristic spaces $K\in{\cal M}_N(S)$ }\\
    \mbox{ together with ample cones }
  \end{array}
  \right\}
\end{array}
$$
There is a natural forgetful map $\underline{{\cal P}}_N\to{\cal M}_N$, 
given by forgetting the ample cones.
Then, we have the following result, and refer to
\cite[Proposition 1.16]{Ogus PM} for details and proof.

\begin{Theorem}[Ogus]
  The functor $\underline{{\cal P}}_N$ is represented by a
  scheme ${\cal P}_N$, which is locally of finite type, 
  almost proper, and smooth of 
  dimension $\sigma_0(N)-1$ over $\FF_p$.
  The natural map
  $$
     {\cal P}_N\,\to\,{\cal M}_N
  $$
  is \'etale, surjective, and locally of finite type.
\end{Theorem}

We repeat that the morphism ${\cal P}_N\to{\cal M}_N$ is 
neither of finite type nor separated, whereas
${\cal M}_N$ is smooth, projective and of finite type over $\FF_p$
by Proposition \ref{prop: moduli of generatrices}.

Now, for an algebraic space $B$ over $\FF_p$ and a family $X\to B$
of $N$-marked K3 surfaces, that is, an element of ${\cal S}_N(B)$,
we have an associated family of $N$-rigidified
supersingular K3 crystals, that is, an element of ${\cal M}_N(B)$.
This gives rise to a morphism
$$
   \pi_N \,:\, {\cal S}_N \,\to\,{\cal M}_N,
$$
which is surjective, but not an isomorphism.
Now, for every point $b\in B$, there is a unique connected component 
of $V_{\NS(X_b)\otimes\RR}$ that
contains the classes of all ample invertible sheaves of $X_b$, 
thereby equipping the family of K3 crystals with ample cones.
This induces a morphism
$$
\widetilde{\pi}_N\,:\,{\cal S}_N\,\longrightarrow\,{\cal P}_N
$$
that lifts $\pi_N$ from above, and is called the {\em period map}.
By \cite[Theorem III]{Ogus PM}, it is an isomorphism.

\begin{Theorem}[Ogus' Crystalline Torelli Theorem]
 \label{thm: crystalline torelli theorem}
 Let $N$ be a supersingular K3 lattice in characteristic $p\geq5$.
 Then, the period map $\widetilde{\pi}_N$ is an isomorphism.
\end{Theorem}

Since ample cones are sometimes inconvenient to handle,
and since they are also responsible for ${\cal P}_N$ being neither separated
nor of finite type,
let us note the following useful application
of the crystalline Torelli theorem:
if two supersingular K3 surfaces have isomorphic
K3 crystals, then they correspond via $\widetilde{\pi}_N$
to points in the same fiber of ${\cal P}_N\to{\cal M}_N$.
In particular, the two surfaces are abstractly isomorphic,
and we obtain the following result, and refer to
\cite[Theorem 1]{Ogus PM} for details.

\begin{Corollary}
  \label{cor: k3 crystals suffice}
  Two supersingular K3 surfaces in characteristic $p\geq5$
  are isomorphic if and only if their associated K3 crystals are isomorphic.
\end{Corollary}

Theorem \ref{thm: crystalline torelli theorem} is the main result of \cite{Ogus PM}.
Let us roughly sketch its proof: 
The existence of the period map $\widetilde{\pi}_N$ is clear.
Separatedness of $\widetilde{\pi}_N$ follows from a theorem of 
Matsusaka and Mumford \cite{MM}.
Properness of $\widetilde{\pi}_N$ follows from a theorem of Rudakov and Shafarevich 
\cite{RS degeneration} that supersingular K3 surfaces have potential good reduction,
that is, given a supersingular K3 surface $X$ over $K:=k((t))$, there exists a finite extension
$R'\supseteq R:=k[[t]]$, say with field of fractions $K'$, 
and a smooth model of $X\times_K K'$ over $R'$ 
(this result uses that $X$ is supersingular -- in general, K3 surfaces do not have
potential good reduction).
Next, $\widetilde{\pi}_N$ is \'etale, which eventually follows from the fact that
$\pi_N:{\cal S}_N\to{\cal M}_N$ is \'etale \cite[Theorem 5.6]{Ogus}, which in turn 
rests on the description of its derivative \cite[Corollary 5.4]{Ogus}.
Finally, to prove that $\widetilde{\pi}_N$ is an isomorphism, it suffices to find one
point $\zeta\in{\cal P}_N$ such that $\widetilde{\pi}_N^{-1}(\zeta)$ consists of a single
point - this is done by taking $\zeta$ to be the supersingular K3 surface that
is the Kummer surface for the self-product of a supersingular elliptic curve.
We refer to \cite[Section 3]{Ogus PM} for details.

\section{Formal Group Laws}
\label{sec:formal groups}

In this section, we introduce formal group laws, which,
at first sight, looks rather independent from what we studied so far.
Before explaining, why this is not so, let 
us first give the prototype of such an object:
let $G$ be a group scheme, say, of finite type and smooth 
over a field $k$.
If $(\OO_{G,O},\idealm)$ is the local ring at the
neutral element $O\in G$, then $(\idealm/\idealm^2)^\vee$ yields
the Lie algebra $\Lieg$, which captures first order infinitesimal information
of $G$ around $O$.
Using the group structure,
the formal completion of $G$ along $O$
$$
   \widehat{G} \,:=\, {\rm Spf}\, \varprojlim \OO_{G,O}/\idealm^n
$$
becomes a group object in the category of formal
schemes, and this {\em formal group law} lies somewhere
between $\Lieg$ and $G$.
We will see that commutative formal group laws can be classified
via their {\em Cartier--Dieudonn\'e modules}, which are 
$W(k)$-modules that resemble $F$-crystals.

For example, if $X$ is a smooth and proper variety over $k$,
then the Picard scheme $\Pic_{X/k}$ is a group scheme over 
$k$, whose formal completion along the origin 
is called the {\em formal Picard scheme}.
If $k$ is perfect of positive characteristic, then
the Cartier--Dieudonn\'e module of
the formal Picard scheme determines $\Hcris{1}(X/W)$.

As shown by Artin and Mazur \cite{Artin Mazur},
the formal Picard group is just the $n=1$-case of a whole series 
of formal group laws $\Phi_{X/k}^n$ that arise from 
cohomological deformation functors, and we refer to Section \ref{sec:Artin-Mazur}
for precise definitions.

For us, $\Phi^2_{X/k}$, which is called the {\em formal Brauer group}, 
will be most relevant in the sequel:
its Cartier--Dieudonn\'e module controls the
part of $\Hcris{2}(X/W)$ that is of slope less than $1$.
What makes this formal group law so fascinating is that
despite its appearance it does in general {\em not} arise as formal completion 
of some group scheme associated to $X$.
Later, in Section \ref{sec:unirational}, we will use the formal Brauer
group to construct non-trivial deformations of supersingular K3 surfaces,
which ultimately proves their unirationality.

\subsection{Formal group laws}
We start with a short introduction to commutative formal group laws
and their classification,
and refer to \cite{Hazewinkel} for the general theory 
of formal group laws, and especially to \cite{Zink} for
the theory of Cartier--Dieudonn\'e modules.

\begin{Definition}
  An $n$-dimensional {\em formal group law} over a ring $k$
  consists of $n$ power series $\vec{F}=(F_1,...,F_n)$
  $$
      F_i(x_1,...,x_n,y_1,...y_n)\,\in\,k[[x_1,...,x_n,y_1,...,y_n]],
      \mbox{ \quad }i=1,...,n
  $$
  such that for all $i=1,...,n$
  $$\begin{array}{lcl}
     F_i(\vec{x},\vec{y}) &\equiv&x_i+y_i\mbox{ modulo terms of degree $\geq2$, and } \\
     F_i ( \vec{x}, F_i(\vec{y},\vec{z})) &=& F_i ( F_i(\vec{x},\vec{y}),\vec{z}),
   \end{array}$$
  where we use the notation $\vec{x}$ to denote $(x_1,...,x_n)$, etc.
  A formal group law $\vec{F}$ is called {\em commutative} if
  $$
  \begin{array}{lcl}
    F_i(\vec{x},\vec{y}) &=&F_i(\vec{y},\vec{x}) \mbox{ \quad for all }i=1,...,n.
  \end{array}
  $$ 
  A {\em homomorphism} $\vec{\alpha}:\vec{F}\to\vec{G}$
  from an $n$-dimensional formal group law $\vec{F}$ 
  to an $m$-dimensional formal group law $\vec{G}$
  consists of $m$ formal power series $\vec{\alpha}=(\alpha_1,...,\alpha_m)$
  in $n$ variables such that $\alpha_i(\vec{x})\equiv0\mbox{ mod degree $1$}$ 
  and
  $$
  \begin{array}{lcl}
     \alpha_i(\vec{x}) &\equiv&0\mbox{ modulo terms of degree $\geq1$, and} \\
     \vec{\alpha}(\vec{F}(\vec{x},\vec{y})) &=& \vec{G}(\vec{\alpha}(\vec{x}),\vec{\alpha}(\vec{y})).
  \end{array}
  $$
  It is called an {\em isomorphism} if there exists
  a homomorphism $\vec{\beta}:\vec{G}\to\vec{F}$
  such that $\vec{\alpha}(\vec{\beta}(\vec{x}))=\vec{x}$ and
  $\vec{\beta}(\vec{\alpha}(\vec{y}))=\vec{y}$.
  An isomorphism $\vec{\alpha}$ is called {\em strict}
  if $\alpha_i(\vec{x})\equiv x_i$ modulo terms of degree $\geq2$ for all $i=1,...,n$.
\end{Definition}

For example, for an integer $n\geq1$, we define {\em multiplication by $n$} to be
$$
 [n]\,(\vec{x})\,:=\,\underbrace{\vec{F}(\vec{F}(...,\vec{x}),\vec{x})}_{n\, {\rm times} }
 \,\in\,k[[\vec{x}]]\,.
$$
If $\vec{F}$ is a commutative formal group law, then
$[n]:\vec{F}\to\vec{F}$ is a homomorphism of formal group laws.
The group laws arising from algebraic varieties, which we will introduce in 
Section \ref{sec:Artin-Mazur} below, will all be commutative.
Let us also mention the following result, which will not
need in the sequel, and refer to \cite[Section 6.1]{Hazewinkel} for details.

\begin{Theorem}
 A $1$-dimensional formal group law over a reduced ring 
 is commutative.
\end{Theorem}

\begin{Examples}
  Here are two basic examples of $1$-dimensional formal group laws.
  \begin{enumerate}
   \item The {\em formal additive group} $\widehat{\GG}_a$ is defined
      by $F(x,y)=x+y$.
   \item The {\em formal multiplicative group} $\widehat{\GG}_m$ is defined
      by $F(x,y)=x+y+xy$.
  \end{enumerate}
  Both are commutative, and their names will be explained in
  Example \ref{example:formal completion}.
\end{Examples}

Over $\QQ$-algebras, all commutative formal group laws of the 
same dimension are mutually isomorphic:

\begin{Theorem}
  \label{thm: logarithm}
  Let $\vec{F}$ be an $n$-dimensional commutative
  formal group law over a $\QQ$-algebra.
  Then, there exists a unique strict isomorphism 
  $$
    \log_{\vec{F}} \,:\, \vec{F}(\vec{x},\vec{y})\,\to\,\widehat{\GG}_a^n,
  $$
  called the {\em logarithm} of $\vec{F}$.
\end{Theorem}

\begin{Example}
  The logarithm of the formal multiplicative group $\widehat{\GG}_m$ 
  over a $\QQ$-algebra $k$ is explicitly given by
  $$
     \log_{\widehat{\GG}_m}\,:\, x\,\mapsto\, \sum_{n=1}^\infty (-1)^{n+1} \frac{x^n}{n}
  $$
  which also motivates the name.
  Note that this power series only makes sense in
  rings that contain $\frac{1}{n}$
  for all integers $n\geq1$, that is, the base ring must be a $\QQ$-algebra.
\end{Example}

On the other hand, if $k$ is an $\FF_p$-algebra,
then Theorem \ref{thm: logarithm} no longer holds true.
In fact, there exist many $1$-dimensional commutative formal group
laws over algebraically closed fields of positive characteristic that are
{\em not} isomorphic to $\widehat{\GG}_a$.
The following discrete invariant is crucial - we will only define
it for $1$-dimensional formal group laws and refer to
\cite[(18.3.8)]{Hazewinkel} for its definition
for higher dimensional formal group laws.

\begin{Definition}
  Let $F=F(x,y)$ be a $1$-dimensional commutative
  formal group law over a field $k$ of positive 
  characteristic $p$, and let $[p]:F\to F$ be multiplication by $p$.
  Then, the {\em height} $h=h(F)$ of $F$ is defined to 
  be
  \begin{itemize}
   \item[-] $h:=\infty$ in case $[p]=0$, and one also says that
    $F$ is {\em unipotent}.
   \item[-] Else, there exists an integer $s\geq1$ and
       a $0\neq a\in k$, such that
       $$
          [p](x)\,=\,a\,\cdot\, x^{p^s} + \mbox{ higher order terms },
       $$
       in which case we set $h:=s$.
       If $F$ is of finite height, one also says that it is {\em $p$-divisible}.
 \end{itemize}
\end{Definition}

\begin{Remark}
  Let us give more characterizations of the height
  of a $1$-dimensional formal group law $F=F(x,y)$
  over a field $k$ of positive characteristic $p$.
  To do so, we define $F^{(p)}=F^{(p)}(x,y)$ to be the
   formal group law obtained from $F$
   by raising each of the coefficients to the $p^{\rm th}$ power.
   Then, $z\mapsto z^p$ defines a homomorphism
   of formal group laws over $k$
   $$
        \sigma\,:\,F\,\to\,F^{(p)},
    $$
    which is called {\em Frobenius}.

  \begin{enumerate}
   \item If finite, the height is characterized as being the
     largest integer $h$ such that 
     there exists a power series $\beta\in k[[z]]$ with
     $$
       [p](z) \,=\, \beta(z^{p^h}).
     $$
     The series $\beta$ defines a homomorphism $F^{(p^h)}\to F$,
     which leads to the following reformulation:
    \item If finite, the height is characterized as being the
     largest integer $h$ such that 
     there exists a factorization
     $$
     \xymatrix{
       F \ar[r]^{\sigma^{\circ h}}  \ar@/_/[rr]_{[p]}
       & F^{(p^h)} \ar@{.>}[r]^{\exists} & F, }
     $$
     where $\sigma^{\circ h}$ denotes the $h$-fold 
     composition of $\sigma$ with itself.
   \item Finally, if the height is finite, then $[p]$ is an {\em isogeny},
     and then, there exists an integer $m\geq0$ and a 
     homomorphism of formal group laws $\psi:F\to F^{(p^{m})}$ such that
     $\psi\circ[p]=\sigma^{\circ m}$ .
     In particular, if finite, the height is characterized as being the
     smallest integer $h$ such that 
     there exists a factorization
     $$
     \xymatrix{
       F \ar[r]^{[p]}  \ar@/_/[rr]_{\sigma^{\circ h}}
       & F \ar@{.>}[r]^{\exists} & F^{(p^h)}. 
      }
     $$
   \end{enumerate}
  We refer to \cite[Section 18.3]{Hazewinkel} for details and generalizations.
\end{Remark}

\begin{Exercise}
  Over fields of positive characteristic, show that
  $$
    h(\widehat{\GG}_m)=1\mbox{ \quad and \quad }
    h(\widehat{\GG}_a)=\infty.
  $$
  In particular, they are not isomorphic.
\end{Exercise}

The importance of the height lies in the following classifcation
result.

\begin{Theorem}[Lazard]
  Let $k$ be an algebraically closed field of
  positive characteristic.
  \begin{enumerate}
   \item For every integer $h\geq1$ or $h=\infty$
   there exists a $1$-dimensional formal group
   law of height $h$ over $k$.
   \item Two $1$-dimensional formal group laws 
   over $k$ are isomorphic if and only if they
   have the same height.
  \end{enumerate}
\end{Theorem}

\begin{Example}
 \label{example:formal completion}
 Let us show how formal group laws arise from group schemes, which also 
 justifies some of the terminology introduced above.
 Thus, let $G$ be a smooth (commutative) group scheme of dimension $n$
 over a field $k$.
 Let $(\OO_{G,O},\idealm)$ be the local ring at the neutral element $O\in G$.
 Using smoothness, we compute the $\idealm$-adic completion to be
 $$
     \widehat{\OO}_{G,O}\,:=\, \varprojlim \OO_{G,O}/\idealm^m
     \,\iso\,k[[x_1,...,x_n]]\,=\,k[[\vec{x}]]\,,
 $$
 and note that 
 $$
   \widehat{G}\,:=\,{\rm Spf}\,\widehat{\OO}_{G,O}
 $$ 
 is the {\em formal completion of $G$ along $O$}.
 The multiplication $\mu:G\times G\to G$ induces a morphism
 of formal schemes
 $\widehat{\mu}:\widehat{G}\times\widehat{G}\to\widehat{G}$
 that turns $\widehat{G}$ into a (commutative) group object
 in the category of formal schemes.
 Explicitly, $\widehat{\mu}$ corresponds to
 a homomorphism of $k$-algebras
 $$
    \widehat{\mu}^{\#}\,:\,k[[\vec{x}]]\,\to\,k[[\vec{y}]]\,\widehat{\otimes}\,k[[\vec{z}]]
    \,\iso\,k[[\vec{y},\vec{z}]]\,,
 $$
 where $\widehat{\otimes}$ denotes the completed tensor product.
 Clearly, $\widehat{\mu}^{\#}$ is uniquely determined by the images of the generators $x_i$, that
 is, by the $n$ formal power series
 $$
   \vec{G} \,:=\,(G_1,...,G_n)\,:=\,\widehat{\mu}^\#(x_1,...,x_n).
 $$
 It is not difficult to see that $\vec{G}$ is an $n$-dimensional (commutative)
 formal group law, and that it encodes $\widehat{G}$.
 Thus, $\widehat{G}$ carries the information of all infinitesimal 
 neighborhoods of $O$ in $G$, and 
 in particular, of the tangent space at $O$, that is, the Lie algebra $\Lieg$ of $G$.
 Put a little bit sloppily, $\widehat{G}$ lies between $G$ and $\Lieg$.
 Here are some standard examples
 \begin{enumerate}
  \item The completion of the multiplicative group scheme $\GG_m\iso\Spec k[x,x^{-1}]$ is
    the formal multiplicative group law $\widehat{\GG}_m$.
  \item The completion of the additive group scheme $\GG_a\iso\Spec k[x]$ is 
    the formal additive group law $\widehat{\GG}_a$.
  \item If $E$ is an elliptic curve over a field $k$, then the completion
    $\widehat{E}$ is a commutative $1$-dimensional formal group law.
    If $k$ is of positive characteristic, then its height is equal to
    $$
       h(\widehat{E})\,=\,\left\{
       \begin{array}{cl}
         1 & \mbox{ if $E$ is ordinary, and } \\
         2 & \mbox{ if $E$ is supersingular.}
       \end{array}
       \right.
    $$
    Thus, if $k$ is algebraically closed and $E$ is ordinary,
    then $\widehat{E}\iso\widehat{\GG}_m$.
    We refer to \cite[Chapter IV]{Silverman} for more about this formal
    group law.
 \end{enumerate} 
\end{Example}

In order to classify $1$-dimensional 
commutative formal group laws over perfect fields $k$ of positive characteristic 
that are not algebraically closed, or even
higher dimensional commutative formal group laws over $k$,
the height is not sufficient.
To state the general classification result, which is in terms of
modules over some non-commutative ring ${\rm Cart}(k)$,
we first have to define this ring.
Be definition, ${\rm Cart}(k)$ is the non-commutative 
ring $W(k)\langle\langle V\rangle\rangle\langle F\rangle$ 
of power series in $V$ and polynomials in $F$ modulo the relations
$$
FV\,=\,p,\mbox{ \quad }
VrF=V(r),\mbox{ \quad }
Fr\,=\,\sigma(r)F,\mbox{ \quad}
rV\,=\,V\sigma(r) \mbox{\quad for all }r\in W(k),
$$
where $\sigma(r), V(r)\in W(k)$ denote Frobenius and Verschiebung of $W(k)$. 

\begin{Theorem}[Cartier--Dieudonn\'e]
  \label{thm:dieudonne cartier modules}
  Let $k$ be a perfect field of positive characteristic.
  Then, there exists a covariant equivalence of categories
  between 
  \begin{enumerate}
    \item The category of commutative formal group laws over $k$.
    \item The category of left ${\rm Cart}(k)$-modules $M$ such that 
      \begin{enumerate}
         \item $V$ is injective, 
         \item $\cap_i V^i(M)=0$, that is, $M$ is $V$-adically separated, and 
         \item $M/VM$ is a finite-dimensional $k$-vector space.
       \end{enumerate}
    \end{enumerate}
  The left ${\rm Cart}(k)$-module associated to a formal group $\vec{G}$
  under this equivalence
  is called the {\em Dieudonn\'e--Cartier module} of $\vec{G}$,
  and is denoted $D\vec{G}$.
\end{Theorem}

Following \cite[Section 1]{Mumford Formal Groups}, let us sketch
the direction $(1)\rightarrow(2)$ of this equivalence:
let $\vec{G}=(G_1,...,G_n)$ be an $n$-dimensional commutative
formal group law over $k$.
Then, $\vec{G}$ defines a  functor
$$
\begin{array}{ccccc}
\Phi_{\vec{G}} &:& (\mbox{$k$-algebras}) &\to& (\mbox{Abelian groups})\\
&&R&\mapsto&\left\{ (x_1,...,x_n)\in R^n\,| \mbox{ each } x_i \mbox{ nilpotent } \right\}
\end{array}
$$
where we define the group structure on the right by setting
$\vec{x}\oplus_{\Phi_{\vec{G}}}\vec{y}:= \vec{G}(\vec{x},\vec{y})$.
Similarly, we define the functor
$$
\begin{array}{ccccc}
\Phi_{\vec{W}} &:& (\mbox{$k$-algebras}) &\to& (\mbox{Abelian groups})\\
&&R&\mapsto&\left\{ (x_0,x_1,...)\,\left |\, 
\begin{array}{l}
   x_i\in R, \mbox{ each }x_i \mbox{ nilpotent, } \\
   \mbox{and almost all }x_i=0 
 \end{array} \right. \right\}
\end{array}
$$
using $\vec{W}:=(W_0,W_1,...)$ to define the group structure,
where the $W_i$ are the Witt polynomials from Section \ref{subsec:Witt}.
We note that $\vec{W}$ is an example of an infinite dimensional
formal group law.
Next, we define
$$
 D\vec{G} \,:=\, {\rm Hom}_{\rm group\, functors/k}(\Phi_{\vec{W}},\Phi_{\vec{G}}).
$$
Multiplication by elements of $W(k)$ gives rise to endomorphisms of $\Phi_{\vec{W}}$,
and one can define Frobenius and Verschiebung for $\Phi_{\vec{W}}$.
Equipped with these operations, ${\rm Hom}(\Phi_{\vec{W}},\Phi_{\vec{W}})$ becomes
a non-commutative $W(k)$-algebra that is isomorphic to the opposite ring of
${\rm Cart}(k)$ defined above.
In particular, this turns $D\vec{G}$ into a left ${\rm Cart}(k)$-module,
which turns out to satisfy the conditions in (2) of
Theorem \ref{thm:dieudonne cartier modules}.
We refer to \cite[Section 1]{Mumford Formal Groups} and \cite[Chapter V]{Hazewinkel} 
for details, generalizations, as well as different approaches.

\begin{Examples}
  \label{examples dieudonne}
  Let $k$ be an algebraically closed field of positive characteristic,
  and let $G=G(x,y)$ be a commutative $1$-dimensional formal group law over $k$.
  \begin{enumerate}
     \item If $G$ is of finite height $h$,
     then $DG\iso {\rm Cart}(k)/(F-V^{h-1})$, which 
     is a free $W(k)$-module of rank $h$.
   \item If $G$ is of infinite height, then $G\iso\widehat{\GG}_a$, and then
     $D\widehat{\GG}_a\iso k[[x]]$ with $F=0$ and $V(x^n)=x^{n+1}$. 
  \end{enumerate}
  In particular, $h$ is equal to the minimal number of generators of $DG$
  as $W(k)$-module.
\end{Examples}

\begin{Exercise}
 \label{exercise: height and slope}
  Let $G$ be the $1$-dimensional commutative formal group law of finite height
  $h$ over an algebraically closed field $k$ of positive characteristic.
  Show that $DG$, considered as a
  $W(k)$-module with the $\sigma$-linear action defined by $\varphi:=F$ is
  an $F$-crystal.
  Show that there exists an isomorphism of $F$-crystals
  $$
     DG\,\iso\, N_\alpha\mbox{ \quad with }\alpha\,=\,(h-1)/h,
  $$
  where $N_\alpha$ is as  in Example \ref{example different hodge polygon}.
  In particular, it is of rank $h$, of slope $\alpha=1-\frac{1}{h}$, and 
  has Hodge numbers $h_0=1$, $h_1=h-1$, and $h_i=0$ if $i\neq0,1$.
\end{Exercise}

\subsection{Formal groups arising from algebraic varieties}
\label{sec:Artin-Mazur}
We now explain how Artin and Mazur \cite{Artin Mazur} associated 
formal group laws to algebraic varieties:
let $X$ be a smooth and proper variety over a field $k$, 
and let  $n\geq1$ be an integer.
Then, we consider the following functor from the category ${\rm Art}_k$ 
of local Artinian $k$-algebras with residue field $k$ to
Abelian groups:
$$
\begin{array}{ccccc}
 \Phi_X^n &:& {\rm Art}_k &\to& (\mbox{Abelian groups}) \\
 && R&\mapsto& \ker\left(
   \Het{n}(X\times_k R,\OO_{X\times_k R}^\times)\stackrel{{\rm res}}{\longrightarrow}
   \Het{n}(X,\OO_X^\times)
 \right)
\end{array}
$$
where $\OO_X^\times$ denotes the sheaf of invertible
elements of $\OO_X$ with respect to multiplication
(we could also write $\GG_m$).
The pro-representability of this functor is studied in
\cite{Artin Mazur}, and there, also a tangent-obstruction theory for
it with tangent space $H^n(X,\OO_X)$ and 
obstruction space $H^{n+1}(X,\OO_X)$
is established.

\begin{Example}
  The case $n=1$ is easy to explain:
  Let $R\in{\rm Art}_k$.
  We identify $\Het{1}(X\times_k R,\OO_{X\times_k R}^\times)$
  with the group of invertible sheaves of $X\times_k R$, and
  then, $\Phi_X^1(R)$ becomes the group of invertible sheaves
  on $X\times_k R$, whose restriction to $X$ is trivial.
  Thus, elements of $\Phi_X^1(R)$ are in bijection to 
  morphisms $\Spec R\to\Pic_{X/k}$, such that the closed point of $\Spec R$ 
  maps to zero, that is, the class of $\OO_X$.
  Thus, if $\widehat{\Pic}_{X/k}$ denotes the completion
  of the Picard scheme $\Pic_{X/k}$ along its zero as
  in Example \ref{example:formal completion}, then we obtain
  an isomorphism of functors 
  $$
   \widehat{\Pic}_{X/k}\,\iso\,\Phi_X^1.
  $$ 
  In particular, $\Phi_X^1$ is pro-representable by a 
  commutative formal group law
  if and only if $\Pic^0_{X/k}$ is smooth over $k$, that is, if and only
  if $\Pic^0_{X/k}$ is an Abelian variety.
  In this case, it is called the {\em formal Picard group},
  and it is of dimension $h^1(\OO_X)$.
\end{Example}

Let us now turn to $\Phi_X^2$, which classifies deformations
of $0\in\Het{2}(X,\OO_X^\times)$.
The group $\Het{2}(X,\OO_X^\times)$ is called the 
{\em (cohomological) Brauer group} of $X$, and we refer to \cite{Gille Szamuely} 
for its algebraic aspects, and to \cite{Grothendieck Brauer} for the
more scheme-theoretic side of this group.
Unlike the Picard group, there is in general no Brauer scheme,
whose points parametrize elements of the Brauer group of $X$.
In particular, $\Phi_X^2$, unlike $\Phi_X^1$,
does {\em not} seem to arise as completion of some group 
scheme associated to $X$.
Nevertheless, we can still study the functor $\Phi^2_X$.
For example, if $h^1(\OO_X)=h^3(\OO_X)=0$ (this holds, for example, for
K3 surfaces), then $\Phi_X^2$ is pro-representable by 
a commutative formal group law of dimension $h^2(\OO_X)$, 
the {\em formal Brauer group}, which is denoted
$$
  \widehat{{\rm Br}}_X\,:=\,\Phi_X^2.
$$
In Section \ref{sec:unirational}, we will use this formal group law to 
construct non-trivial $1$-dimensional deformations of supersingular
K3 surfaces, which is the key to proving their unirationality.

For $n\geq3$, the functors $\Phi^n_X$ are far less understood.
However, we refer the interested reader to \cite{Katsura Geer CY} for 
an analysis of $\Phi_X^n$ in case $X$ is 
an $n$-dimensional Calabi--Yau variety.

\subsection{The connection to Witt vector and crystalline cohomology}
In Section \ref{subsec:Witt}, we introduced Serre's Witt vector cohomology groups
$H^n(W\OO_X)$, which, by functoriality of the Witt vector construction,
carry actions of Frobenius and Verschiebung.
In particular, all these cohomology groups are left ${\rm Cart}(k)$-modules.
The following result from \cite[Proposition (II.2.13)]{Artin Mazur}
and \cite[Corollary (II.4.3)]{Artin Mazur}
links this ${\rm Cart}(k)$-module structure to the Cartier--Dieudonn\'e modules 
of commutative formal group laws associated to the $\Phi^n_X$.

\begin{Proposition}[Artin--Mazur]
 \label{prop: formal group witt vector cohomology}
  Let $X$ be a proper variety over a perfect field $k$
  and assume that $\Phi_X^n$ is pro-representable
  by a formal group law $\vec{F}$
  (for example, this holds true if $h^{n-1}(\OO_X)=h^{n+1}(\OO_X)=0$).
  Then, there exists an isomorphism of left ${\rm Cart}(k)$-modules
  $$
      D\vec{F}\,\iso\, H^n(X,W\OO_X).
  $$
\end{Proposition}

To link the formal group law $\Phi_X^n$ to crystalline cohomology, we use
the slope spectral sequence from Section \ref{sec:crystalline cohomology}.
As mentioned there, it degenerates at $E_1$ if and only if the torsion
of the Hodge--Witt cohomology groups is finitely generated.
However, combining the previous proposition with 
Examples \ref{examples dieudonne}, we see that if 
$\Phi_X^n\iso\widehat{\GG}_a$ (for example, if $n=2$ and $X$ is a
supersingular K3 surface),  then $H^n(W\OO_X)$ will not be finitely
generated and the slope spectral sequence will not degenerate at 
$E_1$, see \cite[Th\'eor\`eme II.2.3.7]{Illusie deRham Witt}.
On the other hand, the slope spectral sequence always degenerates
at $E_1$ modulo torsion, and from this, we obtain 
an isomorphism of $F$-isocrystals
$$
  H^n(X,W\OO_X)\otimes_W K\,\iso\,(\Hcris{n}(X/W)\otimes_W K)_{[0,1[},
$$
where the right hand denotes the direct sum of sub-$F$-isocrystals
of slope strictly less than $1$.
(Here, a little bit of background: 
the point is that the $H^j(W\Omega_X^i)\otimes_W K$ are finite-dimensional
$K$ vector spaces and their sets of slopes are disjoint.
The slope spectral sequence degenerates at $E_1$ after tensoring with $K$,
and since it can be made compatible with the Frobenius actions on both sides, 
the isogeny decomposition of the $F$-isocrystal $\Hcris{n}(X/W)\otimes_W K$ 
can be read off from the isogeny decomposition of the $F$-isocrystals
$H^j(W\Omega_X^i)\otimes_W K$, where $i+j=n$.
From this, it is not so difficult to see that all $F$-isocrystals of slope 
strictly less than $1$ in $\Hcris{n}(X/W)\otimes_W K$ 
arise from $H^n(W\OO_X)\otimes_W K$ via the slope spectral sequence.)

\begin{Example}
\label{example: slope less than one}
 Assume that $\Phi_X^n$ is pro-representable by a 
 $1$-dimensional formal group law of {\em finite} height $h$.
 By Exercise \ref{exercise: height and slope} and
 Proposition  \ref{prop: formal group witt vector cohomology},
 we have isomorphisms 
 $$
   N_{(h-1)/h}\otimes_W K\,\iso\,D\Phi_X^n \otimes_W K \,\iso\,(\Hcris{n}(X/W)\otimes_W K)_{[0,1[}.
 $$
 In particular, $1-\frac{1}{h}$ is the only slope of $\Hcris{n}(X/W)$ less
 than $1$.
\end{Example}

This example applies to the case where $X$ is a K3 surface,
and $n=2$, that is, $\Phi_X^n$ is the formal Brauer group.
In this special case, we can say more.
Before doing so, we remind the reader that we classified
the possible slopes and Newton polygons
of $\Hcris{2}(X/W)$ in Exercise \ref{exercise: k3 crystals}
in terms of some parameter $h$.

\begin{Proposition}
  \label{prop: height determines newton}
  Let $X$ be a K3 surface over a perfect field $k$
  of positive characteristic and let
  $h$ be the height of its formal Brauer group.
  \begin{enumerate}
   \item If $h<\infty$, then the slopes and multiplicities of $\Hcris{2}(X/W)$ are
     as follows
      $$\begin{array}{lccc}
         \mbox{ slope } & 1-\frac{1}{h} & 1 & 1+\frac{1}{h} \\
         \mbox{ multiplicity } & h & 22-2h & h 
       \end{array}$$
   \item If $h=\infty$, then $\Hcris{2}(X/W)$ is of slope $1$ with multiplicity $22$.
  \end{enumerate}
  In particular, $h$ determines the Newton polygon of $\Hcris{2}(X/W)$.
\end{Proposition}

\prf
By Proposition \ref{prop: formal group witt vector cohomology},
$\Phi_X^2$ is pro-representable by a $1$-dimensional formal group law
$\widehat{{\rm Br}}_X$, and let $h$ be its height.

If $h=\infty$, then $D\widehat{{\rm Br}}_X\otimes_W K=0$
by Examples \ref{examples dieudonne}, and so, there are no slopes 
of $\Hcris{2}(X/W)$ less than $1$.
Thus, by Exercise \ref{exercise: k3 crystals}, $\Hcris{2}(X/W)$ is of slope $1$
with multiplicity $22$.

If $h<\infty$, then the considerations of Example \ref{example: slope less than one}
show that the only slope less than $1$ of $\Hcris{2}(X/W)$ is equal to $1-\frac{1}{h}$.
By Exercise \ref{exercise: k3 crystals}, the Newton
polygon of $\Hcris{2}(X/W)$ has the stated slopes and multiplicities.
\qed\medskip

\begin{Exercise}
  \label{exercise: height and picard rank}
  Let $X$ be a K3 surface over an algebraically closed field $k$ of positive characteristic,
  let $h$ be the height of its formal Brauer group, and let $\rho$ be the
  rank of its N\'eron--Severi group.
  \begin{enumerate}
   \item In case $h<\infty$ use the previous proposition and the fact that the image of the crystalline
      Chern class has slope $1$ to deduce the
      {\em Igusa--Artin--Mazur inequality}
      $$
        \rho \,\leq\, 22\,-\,2h.
      $$
      (In fact, this inequality can be generalized to arbitrary smooth projective varieties, 
      see  \cite[Proposition (II.5.12)]{Illusie deRham Witt}.)
      Since $X$ is projective, we have $\rho\geq1$, and therefore,
      $h=11$ is impossible.
    \item In case $\rho=22$ show that $h=\infty$.
      In fact, the Tate conjecture (see Conjecture \ref{conj: Tate conjecture}
      and Theorem \ref{thm: two supersingular k3 notions})
      predicts the equivalence of $\rho=22$ and
      $h=\infty$.
    \item Assuming the Tate conjecture, show that there exist no
     K3 surfaces with $\rho=21$
     (this observation is due to Swinnerton-Dyer, see \cite{Artin Supersingular}).
   \end{enumerate}
\end{Exercise}

We refer to \cite[Appendix B]{Hazewinkel} for more results on formal group laws
arising in algebraic geometry, as well as further references.

\section{Unirational K3 Surfaces}
\label{sec:unirational}

Over algebraically closed fields, a curve is rational
if and only if it is unirational by L\"uroth's theorem.
By a theorem of Castelnuovo, this is also true for surfaces
in characteristic zero.
On the other hand, Zariski constructed examples
of unirational surfaces that are not rational over algebraically
closed fields of positive characterstic.
Since the characterization of unirational surfaces in 
positive characteristic is still unclear,
K3 surfaces provide an interesting testing ground.

The first examples of unirational K3 surfaces were constructed by 
Rudakov and Shafarevich in characteristic $2$ and $3$, as well as by 
Shioda in arbitrary characteristic.
Also, Artin and Shioda showed that unirational K3 surfaces are 
supersingular (using different notions of supersingularity, now
known to be equivalent by the established Tate-conjecture).
Conversely, Artin, Rudakov, Shafarevich, and Shioda conjectured
that supersingular K3 surfaces are unirational, which would
give a cohomological characterization of unirationality.

In this section, we will use the formal Brauer group to associate
to a supersingular K3 surface $X$ over an algebraically
closed field $k$ of positive characteristic 
a $1$-dimensional family ${\cal X}\to\Spec k[[t]]$ of supersingular
K3 surfaces with special fiber $X$, such that 
generic and special fiber are related by a
{\em purely inseparable isogeny}.
In particular, the generic fiber is unirational if and only if
the special fiber is.
Then, we fill up the moduli space of supersingular K3 surfaces
with such families, which implies that every two supersingular K3
surfaces are purely inseparably isogenous.
Since Shioda established the existence of some unirational
K3 surfaces in every positive characteristic, the 
existence of these isogenies implies that all of them are,
thereby establishing the Artin--Rudakov--Shafarevich--Shioda
conjecture.

\subsection{The L\"uroth problem}
We start by recalling some definitions and classical facts
concerning the (uni-)rationality of curves and surfaces.
In this section, the ground field $k$ is always assumed to 
be algebraically closed, in order to avoid the distinction between 
(uni-)rationality over $k$ and $\overline{k}$. 

\begin{Definition}
  An $n$-dimensional variety $X$ over an algebraically closed field $k$
  is called {\em unirational}, if there exists a dominant and rational
  map
  $$
   \PP_k^n\,\dashrightarrow\,X\,.
  $$
  Moreover, $X$ is called {\em rational} if this map can be chosen to be
  birational.
\end{Definition}

Equivalently, $X$ is rational if and only if $k(X)\iso k(t_1,...,t_n)$, and
$X$ is unirational if and only if $k(X)\subseteq k(t_1,....t_n)$ of finite index.
In particular, rational varieties are unirational,
which motivates the following question.

\begin{Question}[L\"uroth]
  Are unirational varieties rational?
\end{Question}

L\"uroth showed that this is true for curves, which nowadays is
an easy consequence of the Riemann--Hurwitz formula, see for example,
\cite[Example IV.2.5.5]{Hartshorne}.
Next, Castelnuovo (characteristic zero) and Zariski (positive characteristic) 
showed that a surface is rational if and only if 
$h^1(\OO_X)=h^0(\omega_X^{\otimes 2})=0$,
that is, we have a cohomological criterion for rationality.
Using this, one can show that L\"uroth's question also has a positive
answer for surfaces in characteristic zero.
We refer to \cite[Theorem VI.3.5]{bhpv} or \cite[Chapter V]{beauville book}
for details and proofs in characteristic zero, and to
\cite[Section 9]{Liedtke overview} for a discussion in positive
characteristic.
In particular, a K3 surface in characteristic zero cannot be unirational.

On the other hand, Zariski \cite{Zariski Castelnuovo} gave examples of unirational surfaces
over algebraically closed fields of positive characteristic
that are not rational, see also \cite[Section 9]{Liedtke overview} for more examples
and discussion, as well as \cite{Liedtke Schuett} for results that show that unirationality
is quite common even among simply connected surfaces of general type in positive
characteristic.

Finally, there are $3$-dimensional unirational varieties over algebraically closed
fields of characteristic zero that are not rational by results of
Iskovskih and Manin \cite{Manin}, Clemens and Griffiths \cite{Clemens}, as well
as Artin and Mumford \cite{Artin Mumford}.

\subsection{Unirational and supersingular surfaces}
By Castelnuovo's theorem, a surface in characteristic zero is 
rational if and only if it is unirational.
Although this is not true in positive characteristic, we have
at least the following necessary condition for unirationality.

\begin{Theorem}[Shioda $+\varepsilon$]
  Let $X$ be a smooth, proper, and unirational surface over an algebraically closed field $k$
  of positive characteristic.
  \begin{enumerate}
   \item The Picard rank $\rho$ is equal to the second Betti number $b_2$, that is, $X$ is
     Shioda-supersingular.
   \item The crystal $\Hcris{2}(X/W)$ is of slope $1$, that is,
     $X$ is Artin-supersingular.
  \end{enumerate}
\end{Theorem}

\prf
Assertion (1) is shown in \cite[Section 2]{Shioda 1974}.
Assertion (2) follows from (1), since $c_1(\NS(X))\otimes W$
defines an $F$-crystal of slope $1$ and rank $\rho$ inside 
$\Hcris{2}(X/W)$, which is of rank $b_2$,
see the discussion in Section \ref{sec: Tate conjecture}
and the proof of Theorem  \ref{thm: two supersingular k3 notions}.
\qed\medskip

The proof also shows that we have the following implications
for surfaces in positive characteristic
$$
 \mbox{ unirational } \,\Rightarrow\, \mbox{ Shioda-supersingular }
 \,\Rightarrow\, \mbox{ Artin-supersingular}\,.
$$
For K3 surfaces in odd characteristic, both notions of supersingularity are equivalent by the
Tate conjecture, see 
Theorem \ref{thm: two supersingular k3 notions}.
Thus, it is natural to ask whether all converse implications hold, at least, for
K3 surfaces.

\begin{Conjecture}[Artin, Rudakov, Shafarevich, Shioda]
  \label{conj: Shioda}
  A K3 surface is unirational if and only if it is supersingular.
\end{Conjecture}

Let us note that Shioda \cite[Proposition 5]{Shioda supersingular unirational} gave
an example of a Shioda-supersingular Godeaux surface that is not unirational.
However, his examples have a fundamental group of order $5$ and the non-unirationality
of them is related to congruences of the characteristic modulo $5$.
Therefore, in loc. cit. he asks whether simply connected and supersingular surfaces are 
unirational, which is wide open.
Coming back to K3 surfaces, Conjecture \ref{conj: Shioda} 
holds in the following cases:

\begin{Theorem}
 \label{thm: previous unirationality results}
  Let $X$ be a K3 surface characteristic $p>0$.
  Then, $X$ is unirational in the following cases:
  \begin{enumerate}
    \item $X$ is the Kummer surface ${\rm Km}(A)$ for a supersingular Abelian
       surface and $p\geq3$ by Shioda \cite{Shioda 1977}. 
    \item $X$ is Shioda-supersingular and 
      \begin{enumerate}
         \item $p=2$ by Rudakov and Shafarevich \cite{RS supersingular}.
         \item $p=3$ and $\sigma_0\leq6$ by Rudakov and Shafarevich \cite{RS supersingular}.
         \item $p=5$ and $\sigma_0\leq3$ by Pho and Shimada \cite{Pho Shimada}.
         \item $p\geq3$ and $\sigma_0\leq2$ by (1), since these surfaces 
           are precisely Kummer surfaces for supersingular Abelian surfaces
           (see Examples \ref{examples: supersingular K3}).
       \end{enumerate}
  \end{enumerate}
\end{Theorem}

In particular, we have examples in every positive characteristic that support
Conjecture \ref{conj: Shioda}.
Let us comment on the methods of proof:
\begin{enumerate}
 \item  Shioda showed (1) by dominating Kummer surfaces by Fermat surfaces,
  that is, surfaces of the form $\{x_0^n+...+x_3^n=0\}\subset\PP^3_k$, 
  and explicitly constructed unirational parametrizations of these latter surfaces
  in case there exists a $\nu$ such that $p^\nu\equiv-1\mod n$, see
  \cite{Shioda 1974}.
  \item Rudakov and Shafarevich \cite{RS supersingular} showed their unirationality
   result using quasi-elliptic fibrations, which can exist only if $2\leq p\leq 3$.
\end{enumerate}
We refer to \cite[Section 9]{Liedtke overview} for more details and further
references.

\subsection{Moving torsors}
An interesting feature of supersingular K3 surfaces is that all of them come with elliptic
fibrations.
In this section we will show that those that admit an elliptic fibration with section
admit very particular $1$-dimensional deformations:
namely, the generic and special fiber are related by purely inseparable isogenies, which implies
that one is unirational if and only if the other one is.

\begin{Definition}
  A {\em genus-$1$ fibration} from  a surface is a proper morphism
  $$
     f\,:\,X\,\to\, B
  $$ 
  from a normal surface $X$ onto a normal curve $B$ with $f_\ast\OO_X\iso\OO_B$
  such that the generic fiber is integral of arithmetic genus $1$.
  In case the geometric generic fiber is smooth, the fibration is called {\em elliptic},
  otherwise it is called {\em quasi-elliptic}.
  In both cases, the fibration is called {\em Jacobian} if it admits a secion.
\end{Definition}

If $X\to B$ is a K3 surface together with a fibration,
then $B\iso\PP^1$
for otherwise, the Albanese map of 
$X$ would be non-trivial, contradicting $b_1(X)=0$.
Moreover, if $F$ is a fiber, then $F^2=0$, and since $\omega_X\iso\OO_X$,
the adjunction formula yields $2p_a(F)-2=F^2+K_XF=0$, that is, the 
fibration is of genus $1$.
If the geometric generic fiber is singular, then, by Tate's theorem 
\cite{Tate Genus}, 
the characteristic $p$ of the ground field $k$ satisfies $2\leq p\leq 3$.
In particular, if $p\geq5$, then every genus-$1$ fibration is
generically smooth, that is, elliptic.

\begin{Theorem}
   Let $X$ be a supersingular K3 surface in odd characteristic $p$, or,
   Shioda-supersingular if $p=2$.
   \begin{enumerate}
    \item $X$ admits an elliptic fibration.
    \item If $p=2$ or $p=3$ and $\sigma_0\leq6$, then $X$
    admits a quasi-elliptic fibration.
    \item If $p\geq5$ and $\sigma_0\leq9$, then $X$ admits
      a Jacobian elliptic fibration.
   \end{enumerate}
\end{Theorem}

\prf
Since indefinite lattices of rank $\geq5$ contain non-trivial isotropic vectors,
the N\'eron--Severi lattice $\NS(X)$ of a supersingular K3 surface
contains a class $0\neq E$ with $E^2=0$.
By Riemann--Roch on K3 surfaces, $E$ or $-E$ is effective, say $E$.
Then, the Stein factorization of 
$X\dashrightarrow\Proj\bigoplus_n H^0(X,\OO_X(nE))$ eventually
gives rise to a (quasi-)elliptic fibration, see, for example, \cite[Section 3]{RS K3} or the proof of 
\cite[Proposition 1.5]{Artin Supersingular} for details.
Moreover, if $X$ admits a quasi-elliptic fibration, it also admits 
an elliptic fibration, see \cite[Section 4]{RS K3}.

Assertion (2) follows from the explicit classification of N\'eron--Severi lattices of 
supersingular K3 surfaces and numerical criteria for the existence of
a quasi-elliptic fibration, see \cite[Section 5]{RS K3}.

Assertion (3) follows again from the explicit classification of N\'eron--Severi
lattices of supersingular K3 surfaces, 
see \cite[Proposition 3.9]{Liedtke Unirational K3} and \cite[Remark 3.11]{Liedtke Unirational K3}.
\qed\medskip

Let us now explain how to deform Jacobian elliptic fibrations on {\em supersingular} K3 surfaces
to non-Jacobian elliptic fibrations:
thus, let $X\to\PP^1$ be a Jacobian elliptic fibration from
a K3 surface.
Contracting the components of the fibers that do not meet the
zero section, we obtain the {\em Weierstra\ss\ model} $X'\to\PP^1$.
This fibration has irreducible fibers, $X'$ has at worst
rational double point singularities, and we let
$A\to\PP^1$ be the smooth locus of $X'\to\PP^1$.
Then, $A\to\PP^1$ is a relative group scheme - more precisely,
it is the identity component of the N\'eron model of $X\to\PP^1$.
Now, one can ask for commutative diagrams of the form
$$
\begin{array}{ccccc}
 X & \supseteq & A & \to & {\cal A} \\
 \downarrow && \downarrow &&  \downarrow\\
 \PP^1_k & = & \PP^1_k & \to & \PP^1_S\\
 \downarrow && \downarrow && \downarrow \\
 \Spec k & = & \Spec k & \to & S
\end{array}
$$
where $S$ is the (formal) spectrum of a local complete and Noetherian $k$-algebra
with residue field $k$, the right squares are Cartesian,
where ${\cal A}\to\PP^1_S$ is a family of elliptic fibrations with
special fiber $A\to\PP^1_k$, 
and where all elliptic fibrations are torsors (principal homogeneous 
spaces) under the Jacobian elliptic fibration $A\to\PP^1_k$.
To bound the situation and in order to algebraize, we assume that these
{\em moving torsors} come with a relative invertible sheaf on
${\cal A}\to\PP^1_S$ that is of some fixed degree $n$.
In \cite[Section 3.1]{Liedtke Unirational K3}, we showed that such 
families correspond to $n$-torsion elements in the formal
Brauer group, that is, elements in the Abelian group
$$
   \widehat{{\rm Br}}_{X}(S)[n]\,.
$$
To obtain non-trivial moving torsors over $S:=\Spec k[[t]]$, this group
must be non-zero, and for this to be the case,
it simply follows from the fact that $\widehat{{\rm Br}}_X$ is a $1$-dimensional
formal group law, that $k$ must be of positive characteristic $p$,
that $\widehat{{\rm Br}}_X$ must be isomorphic to $\widehat{\GG}_a$,
and that $p$ has to divide $n$, see \cite[Lemma 3.3]{Liedtke Unirational K3}.
In particular, $X$ must be a supersingular K3 surface.
In case $p=n$, one can even find a degree-$p$ multisection $D\to Y_S$, 
that is purely inseparable over the base $\PP^1_S$, 
see \cite[Proposition 3.5]{Liedtke Unirational K3}.\footnote{Unfortunately, this is not true: one can find a degree $p$ multisection, but it is not clear whether it is purely inseparable over the base $\PP^1_S$, see also \cite{Liedtke Erratum}. In general, one cannot find such a purely inseparable section as the results of Bragg and Lieblich \cite{BL} show.}.
Next, we use these considerations to construct a family
of supersingular K3 surfaces.

\begin{Theorem}
  \label{thm: moving torsor}
  Let $X$ be a supersingular K3 surface in characteristic $p\geq5$.
  \begin{enumerate}
    \item If $\sigma_0\leq9$, then $X$ admits a Jacobian elliptic fibration.
    \item Associated to a Jacbobian elliptic fibration on $X$, there exists a smooth projective
       family of supersingular K3 surfaces
       $$
       \begin{array}{ccc}
         X &\to& {\cal X}\\
         \downarrow && \downarrow\\
         \Spec k &\to& \Spec k[[t]].
       \end{array}
       $$
       This family has the following properties:
       \begin{enumerate}
        \item
       The Artin invariants of special and geometric generic fiber satisfy
       $$
          \sigma_0({\cal X}_{\overline{\eta}})\,=\,\sigma_0(X)\,+\,1\,.
       $$
       In particular, this family has non-trivial moduli.
       \item There exist dominant and rational maps   
       $$
          {\cal X}_{\overline{\eta}}^{(1/p)}\,\dashrightarrow\,X\times_k \overline{\eta} \,\dashrightarrow\,
          {\cal X}_{\overline{\eta}}^{(p)},
       $$
       both of which are purely inseparable of degree $p^2$.\footnote{This is not true in general. One can show that the two surfaces are related by a sequence of very special correspondences.}
      \end{enumerate}
  \end{enumerate}
\end{Theorem}

The idea is to compactify the moving torsor associated to the Jacobian
elliptic fibration $X\to\PP^1$ that arises from a nontrivial $p$-torsion element of $\widehat{{\rm Br}}_X(S)$
with $S=\Spec k[[t]]$.
After resolving the rational double point singularities in families, we obtain 
a family ${\cal X}\to\PP^1_S$ of supersingular K3 surfaces with special fiber $X$.
Specialization induces
an injection of N\'eron--Severi groups $\NS({\cal X}_{\overline{\eta}})\to\NS(X)$, whose
cokernel is generated by the class of the zero-section of $X\to\PP^1$.
From this, the assertion on Artin invariants follows, 
and we refer to \cite[Theorem 3.6]{Liedtke Unirational K3} for details.
By \cite[Proposition 3.5]{Liedtke Unirational K3}, the elliptic fibration on the 
generic fiber admits a degree-$p$ multisection $D$ that is purely inseparable over the base.\footnote{See the previous footnotes.}
Then, the assertions on dominant and rational maps follow from 
base-changing the family ${\cal X}_{{\eta}}\to\PP^1_{{\eta}}$ 
to $D\to\PP^1_{{\eta}}$, which trivializes the moving torsor,
see \cite[Theorem 3.6]{Liedtke Unirational K3}.

Spreading out the family ${\cal X}\to S$ of Theorem \ref{thm: moving torsor}
to some curve of finite type over $k$, and using the theorem of Rudakov and Shafarevich
on potential good reduction of supersingular K3 surfaces
(see the remarks after Corollary \ref{cor: k3 crystals suffice}), 
we obtain
a smooth and projective family ${\cal Y}\to B$, where $B$ is a smooth 
and projective curve over $k$, with $X$ as some fiber, and ${\cal X}$ as
fiber over the completion $\widehat{\OO}_{B,\eta}$.
Associating to such a family their rigidified K3 cystals, we obtain 
the following statement about moduli spaces of rigidified K3 crystals.

\begin{Theorem}
  \label{thm: moduli fibration}
  Let $N$ and $N_+$ be the supersingular K3 lattices
  in characteristic $p\geq5$ of Artin invariants $\sigma_0$ and $\sigma_0+1$,
  respectively.
  Then, there exists a fibration of moduli spaces
  of rigidified K3 crystals
  $$
    {\cal M}_{N_+}\,\longrightarrow\,{\cal M}_N,
  $$
  whose geometric fibers are rational curves.
\end{Theorem}

We note that the fibers of this fibration correspond to the moving torsor families 
(at least, an open and dense subset does), and refer to 
\cite[Theorem 4.3]{Liedtke Unirational K3} and \cite[Theorem 4.5]{Liedtke Unirational K3}
for details.

As explained in and after Examples \ref{examples: crystalline moduli},
there is only one supersingular K3 surface with Artin invariant
$\sigma_0=1$.
In particular, by Theorem \ref{thm: moving torsor} and Theorem \ref{thm: moduli fibration}
together with an induction on Artin invariants, we can first relate every supersingular K3
surface to the unique supersingular K3 surface with $\sigma_0=1$ via dominant rational maps, 
and ultimately obtain the following structure result for supersingular K3 surfaces.

\begin{Theorem}
  Let $X$ and $Y$ be supersingular K3 surfaces in characteristic $p\geq5$.
  Then, there exist  dominant and rational maps   
  $$
          X\,\dashrightarrow\,Y\,\dashrightarrow\,X
   $$
   which are purely inseparable, that is, the surfaces are 
   {\em purely inseparably isogenous}.\footnote{The best one can say at the moment is that $X$ and $Y$ are related by a sequence of correspondences.}
\end{Theorem}

By a theorem of Shioda \cite{Shioda 1977}, 
supersingular Kummer surfaces in odd characteristic are unirational 
(see also Theorem \ref{thm: previous unirationality results} above), and combining 
this with the previous theorem, this implies Conjecture \ref{conj: Shioda}.

\begin{Theorem}
   Supersingular K3 surfaces in characteristic $p\geq5$ are unirational.\footnote{Since it is unclear whether unirationality is preserved along the correspondences mentioned in the previous footnote, the unirationality of supersingular K3 surfaces remains a conjecture.} 
\end{Theorem}

In particular, a K3 surface in characteristic $p\geq5$ is supersingular in the
sense of Artin if and only if it is supersingular in the sense of Shioda if and only
if it is unirational.

\subsection{Unirationality of moduli spaces}
It follows from Theorem \ref{thm: explicit coordinates} (see
\cite[Proposition 4.10]{Ogus} for details),
as well as Theorem \ref{thm: moduli fibration}, that the
moduli space of $N$-rigidified K3-crystals ${\cal M}_N$ is unirational.
Thus, also the moduli space ${\cal P}_N$ of $N$-rigidified K3-crystals together
with ample cones is in some sense unirational 
(this space is neither separated nor of finite type, but obtained
by glueing open pieces of ${\cal M}_N$).
Moduli spaces of polarized K3 surfaces are much better behaved,
see Theorem \ref{thm: moduli spaces are nice}.
In fact, constructing families of supersingular K3 surfaces using
the formal Brauer group (similar to the moving torsor construction
above, but more general), the supersingular loci inside
moduli spaces of polarized K3 surfaces are rationally connected
- this is a forthcoming result of Lieblich, see
\cite[Section 9]{Lieblich overview} for announcements of
some of these results.

\section{Beyond the Supersingular Locus}
\label{sec:stratification}

In Section \ref{sec: moduli spaces k3}, we introduced and discussed the moduli stack 
${\cal M}_{2d,\FF_p}^\circ$ of degree-$2d$ primitively polarized K3 surfaces over $\FF_p$.
In the previous sections, we focused on supersingular K3 surfaces, that
is, K3 surfaces, whose formal Brauer groups are of infinite height.
In this final section, we collect and survey a couple of results on this moduli space
beyond the supersingular locus.
We stress that these results are just a small outlook, as well as deliberately
a little bit sketchy.

\subsection{Stratification}
 \label{sec:height Artin}
Associated to a K3 in positive characteristic $p$, we associated the following discrete invariants:
\begin{enumerate}
 \item The height $h$ of the formal Brauer group, which satisfies $1\leq h\leq10$ or $h=\infty$,
  see Proposition \ref{prop: height determines newton} and Exercise \ref{exercise: height and picard rank}.
 \item  If $h=\infty$ then ${\rm disc}(\NS(X))=-p^{2\sigma_0}$ for some integer $1\leq\sigma_0\leq10$,
   the Artin-invariant, see Definition \ref{def: artin invariant lattice}.
\end{enumerate}
These invariants allow us to define the following loci 
(just as a set of points, we do not care about scheme structures at the moment) 
inside the moduli space ${\cal M}_{2d,\FF_p}^\circ$ of degree-$2d$ polarized K3 surfaces
$$
\begin{array}{lcl}
 {\cal M}_i  &:=& \left\{ \mbox{ surfaces with $h\geq i$ } \right\} \\
 {\cal M}_{\infty, i} &:=& \left\{ \mbox{ surfaces with $h=\infty$ and $\sigma_0\leq i$ }\right\}\,.
\end{array}
$$
Thus, at least on the set-theoretical level, we obtain inclusions
$$
\begin{array}{rl}
  {\cal M}_{2d,\FF_p}^\circ \,=\,{\cal M}_1\,\supset\,...\,\supset\,{\cal M}_{10}\,\supset\,&{\cal M}_{\infty}\\
  & || \\
  & {\cal M}_{\infty,10}\,\supset\,{\cal M}_{\infty,9}\,\supset\,...\,\supset\,{\cal M}_{\infty,1}.
\end{array}
$$
This stratification was introduced by Artin \cite{Artin Supersingular}, and studied in detail
by Ekedahl, van der Geer, Katsura \cite{Katsura Geer}, \cite{Katsura Geer Brauer},
\cite{Ekedahl van der Geer}, and Ogus \cite{Ogus height}.
It turns out that each ${\cal M}_{i+1}$ is a closed subset of ${\cal M}_i$
and that ${\cal M}_{\infty,i}$ is closed in ${\cal M}_{\infty,i+1}$.
For example, the first closedness assertion can be deduced easily from
the following result \cite[Theorem 5.1]{Katsura Geer}.

\begin{Proposition}[van~der~Geer--Katsura]
  Let $X$ be a K3 surface over an algebraically closed field of positive
  characteristic $p$.
  Then 
  $$
     h\,=\,\min
     \left\{
     n\geq1\,:\, \left(F:H^2(X,W_n(\OO_X))\to H^2(X,W_n(\OO_X))\right)\neq0
     \right\}.
  $$
\end{Proposition}

In fact, this generalizes to higher dimensions:
if $X$ is a Calabi--Yau variety of dimension $n$,
then the height of the one-dimensional formal group law
associated to $\Phi_X^n$
(notation as in Section \ref{sec:Artin-Mazur}) can be
characterized as in the previous proposition, and we refer to
\cite[Theorem 2.1]{Katsura Geer CY} for details. 

\subsection{Stratification via Newton polygon}
\label{sec:height=Newton}
By Proposition \ref{prop: height determines newton},
the height $h$ of the formal Brauer group of a K3 surface
determines the smallest slope of the Newton polygon
of the $F$-crystal $\Hcris{2}$.
Moreover, by Exercise \ref{exercise: k3 crystals},
the smallest slope determines this Newton polygon completely.
In particular, the height stratification (the first part of the stratification
introduced above)
$$
  {\cal M}_{2d,\FF_p}^\circ \,=\,{\cal M}_1\,\supset\,...\,\supset\,{\cal M}_{10}\,\supset\,{\cal M}_{\infty}
$$
coincides with the stratification 
by the Newton polygon associated to the $F$-crystal $\Hcris{2}$.
This stratification also illustrates Grothendieck's theorem that
the Newton polygon goes up under specialization.

\subsection{Stratification via $\mathbf{F}$-zips}
Let $(X,{\cal L})$ be a primitively polarized K3 surface 
over an algebraically closed field $k$ of positive characteristic $p$ such
that $p$ does not divide ${\cal L}^2=:2d$.
The cup product induces a non-degenerate quadratic form on
$\HdR{2}(X/k)$.
(Since we assumed $p\nmid 2d$, we have $p\neq2$ and thus,
we do not have to deal with subtleties of quadratic forms
in characteristic $2$.)
Then, we define the {\em primitive cohomology} to be
$$
   M\,:=\,c_1({\cal L})^\perp \,\subset\,\HdR{2}(X/k).
$$
We note that the condition $p\nmid 2d$ ensures that 
$c_1({\cal L})$ is non-zero, and thus, $M$ is a
$21$-dimensional $k$-vector space.
The Hodge and its conjugate filtration on $\HdR{2}(X/k)$
give rise to two filtrations $C^\bullet$ and $D_\bullet$
on $M$, and the Cartier isomorphism
induces isomorphisms
$\varphi_n:({\rm gr}_C^n)^{(p)}\to{\rm gr}_D^n$,
see Section \ref{subsec:algebraic deRham}.
Next, the quadratic form on $\HdR{2}(X/k)$
induces a non-degenerate quadratic form $\psi$ on 
$M$,
and it turns out that the filtrations are orthogonal
with respect to $\psi$.
Putting this data together, we obtain
an  {\em orthogonal $F$-zip}
$$
  \left( \,
  M,\,C^\bullet,\,D_\bullet,\,\varphi_\bullet,\,\psi
  \,\right)\,
$$
of filtration type $\tau$ with $\tau(0)=\tau(2)=1$ and $\tau(1)=19$,
see Definition \ref{def: f-zip}.
We refer to \cite[Section 5]{Moonen Wedhorn} or \cite{PWZ}
for $F$-zips with additional structure.

As already mentioned in Section \ref{subsec:algebraic deRham}
and made precise by \cite[Theorem 4.4]{Moonen Wedhorn},
$F$-zips of a fixed filtration type form an Artin stack
that has only a finite number of points.
More precisely, orthogonal $F$-zips of type $\tau$ as above
are discussed in detail in \cite[Example (6.18)]{Moonen Wedhorn}.
Let us sketch their results:
if $(V,\psi)$ is an orthogonal space of dimension $21$
over $\FF_p$, then $SO(V,\psi)$ has a root system
of type $B_{10}$.
After a convenient choice of roots, and with appropriate identifications,
the Weyl group $W$ of $SO(V,\psi)$ becomes a subgroup
of the symmetric group $\Sym_{21}$ as follows
$$
 W\,\iso\, \left\{\rho\in \Sym_{21}\,|\, \rho(j)+\rho(22-j)\,=\,22\mbox{ for all }j \right\}.
$$
We set $W_J:=\{\rho\in W\,|\,\rho(1)=1\}$ and it is easy to
see that the set of cosets $W_J\backslash W$ consists of $20$ elements.
As shown in \cite[Example (6.18)]{Moonen Wedhorn},
there exists a bijection between the set of isomorphism classes 
of orthogonal $F$-zips of type $\tau$ over $\overline{\FF}_p$
and $W_J\backslash W$.
Moreover, the Bruhat order on $W$ induces a {\em total} order on
this set of cosets, that is, we can find representatives
$\overline{w}_1>...>\overline{w}_{20}$.
Using this bijection and the representatives, we define 
$$
  {\cal M}^{(i)} \,:=\, \left\{
  \mbox{ surfaces whose associated orthogonal $F$-zip corresponds to $\overline{w}_i$ }
  \right\}.
$$
This gives a decomposition of ${\cal M}_{2d,\FF_p}^\circ$ into $20$
disjoint subsets.
This decomposition is related to the stratification from Section \ref{sec:height Artin}
as follows
$$
\begin{array}{lcll}
  {\cal M}_i \,\backslash\, {\cal M}_{i+1} &=& {\cal M}^{(i)} & \mbox{ for } 1\leq i\leq 10 \\
  {\cal M}_{\infty, 21-i} \,\backslash\, {\cal M}_{\infty, 20-i} &=& {\cal M}^{(i)} & \mbox{ for } 11\leq i\leq 20,
\end{array}
$$
where it is convenient to set ${\cal M}_{11}:={\cal M}_\infty$, and ${\cal M}_{\infty,0}:=\emptyset$.
In particular, both decompositions eventually give rise to the same stratification
of the moduli space.
Again, we refer to \cite[Example (6.18)]{Moonen Wedhorn}
for details.

\subsection{Singularities of the strata}
In \cite{Katsura Geer} and \cite{Ogus height}, Katsura, van der Geer, and 
Ogus found a beautiful description of the
singularities of the height strata ${\cal M}_i$, which, by Section
\ref{sec:height=Newton} coincide with the Newton-strata.

\begin{Theorem}[van der Geer--Katsura, Ogus]
  Let $\{{\cal M}_i\}_{i\geq1}$ be the height stratification of
  ${\cal M}_{2d,\FF_p}^\circ$.
  Still assuming $p\nmid 2d$,
  we have an equality of sets
  $$
     ({\cal M}_i)^{\rm sing} \,=\, {\cal M}_{\infty,i-1}\mbox{ \quad for all \quad } 1\leq i\leq 10,
  $$
  where ${-}^{\rm sing}$ denotes the singular locus of the corresponding height stratum.
\end{Theorem}

\subsection{Cycle classes}
To state the next result, we let $\pi:{\cal X}\to{\cal M}_{2d,\FF_p}^\circ$ be the universal
polarized K3 surface, and we still assume $p\nmid 2d$.
Then, we define the {\em Hodge class} to be the first Chern class
$$
  \lambda_1 \,:=\, c_1 \left( \, \pi_\ast\Omega^2_{{\cal X}/{\cal M}^\circ_{2d,\FF_p}} \, \right).
$$
Let us recall that ${\cal M}^\circ_{2d,\FF_p}$ is $19$-dimensional,
and that each stratum of our stratification is of codimension $1$ inside the next larger stratum.
The following result \cite[Theorem A]{Ekedahl van der Geer} describes the 
cycle classes of these strata in terms of the Hodge class.
We note that the moduli spaces ${\cal M}_{2d,\FF_p}$ are non-complete, but that the formulas still make sense 
on an appropriate compactification.

\begin{Theorem}[Ekedahl--van der Geer]
 \label{thm: ekedahl van der geer hodge}
 In terms of the Hodge class $\lambda_1$, the cycle classes of the strata
 inside ${\cal M}_{2d,\FF_p}$ are as follows
 $$
 \begin{array}{ccl}
    \displaystyle [{\cal M}_i] &\displaystyle=& 
    \displaystyle 
    (p-1)(p^2-1)\cdots(p^{i-1}-1)\lambda_1^{i-1} \\ {}
    [{\cal M}_\infty] &=& \frac{1}{2} (p-1)(p^2-1)\cdots(p^{10}-1)\lambda_1^{10}\\ {} 
    \displaystyle [{\cal M}_{\infty, i}] &\displaystyle=& 
    \displaystyle
    \frac{1}{2}\frac{(p^{2(11-i)}-1)(p^{2(12-i)}-1)\cdots(p^{20}-1)}{(p+1)(p^2+1)\cdots(p^i+1)}\lambda_1^{20-i} 
 \end{array}
 $$
 for $1\leq i\leq 10$.
\end{Theorem}

The appearance of the factor $1/2$ is related to the fact that the formulas of \cite[Theorem 14.2 and Section 15]{Katsura Geer}
count the supersingular stratum doubly, see also \cite{Katsura Geer Brauer} and \cite{Ekedahl van der Geer}.

Theorem  \ref{thm: ekedahl van der geer hodge} measures the ``size'' of these strata, and can be thought of
as a generalization of a theorem of Deuring for elliptic curves:
namely, in characteristic $p$, elliptic curves form a $1$-dimensional moduli 
space over $\FF_p$.
The formal Picard group of an elliptic curve either has height $1$ (ordinary elliptic
curve) or height $2$ (supersingular elliptic curve),
see Example  \ref{example:formal completion}.
Ordinary elliptic curves form an open and dense set, and thus,
the number of supersingular elliptic curves in a fixed characteristic $p$ 
is finite.
Now, theorem of Deuring gives a precise answer:
classically, it is phrased by saying that there are 
$[p/12]+\varepsilon_p$ 
supersingular elliptic curves for some $0\leq\varepsilon_p\leq2$
depending on the congruence class of $p$ modulo $12$, 
see, for example, \cite[Theorem V.4.1]{Silverman}.
However, if we count supersingular elliptic curves and weight each one of them
with respect to their automorphism group
(which can be thought of as counting them on the moduli stack rather than
the coarse moduli space) we obtain the following, much more beautiful formula
$$
 \sum_{E\, {\rm supersingular}/\iso} \frac{1}{\#{\rm Aut}(E)} \,=\, \frac{p-1}{24}.
$$
Theorem  \ref{thm: ekedahl van der geer hodge} is a generalization of
this way of counting to K3 surfaces.

We refer to \cite{Ekedahl van der Geer} for more about the singularities of the strata,
as well as to \cite{Geer Arbeitstagung} for irreducibility results of the strata.

\subsection{A Torelli theorem via Shimura varieties}
We end our survey with a very sketchy discussion of a Torelli
theorem for K3 surfaces in positive characteristic
beyond the supersingular ones.

For curves, Abelian varieties, K3 surfaces,... the classical period map
associates to such a variety some sort of linear algebra data as explained at
the end of Section \ref{subsec:complex geometry}.
Over the complex numbers, this linear algebra data is parametrized as points inside some
Hermitian symmetric domain modulo automorphisms.
Thus, the period map can be interpreted as a morphism from the moduli space
of these varieties to a Hermitian symmetric domain.
In this setting, a Torelli theorem is the statement that this period map is an immersion,
or, at least \'etale.

First, we set up some Shimura varieties, which will serve as the Hermitian symmetric
domain modulo automorphisms:
let $U$ be the hyperbolic plane over $\ZZ$, and set $N:=U^{\oplus3}\oplus E_8^{\oplus2}$.
Thus, abstractly, $N$ is isometric to $H^2(X,\ZZ)$ of a K3 surface with the cup-product pairing
(Poincar\'e duality), also known as the {\em K3 lattice}.
Let $e,f$ be a basis for the first copy of $U$ in $N$.
Then, for $d\geq1$, we define
$$
  L_d\,:=\, \langle e-df\rangle^\perp\,\subseteq\,N\,,
$$
which is modelled on the primitive cohomology $P^2(X):=c_1({\cal L})^\perp\subset H^2(X,\ZZ)$ 
of a polarization ${\cal L}$ of self-intersection $2d$.
We note that $G_d:=SO(L_d)$ is a semi-simple algebraic group over $\QQ$.
Next, let $K_{L_d}\subset G_d(\Aff_f)$ be the largest subgroup of $SO(L_d)(\widehat{\ZZ})$
that acts trivially on the discriminant ${\rm disc}(L_d):=L_d^\vee/L_d$.
Finally, let $Y_{L_d}$ be the space of oriented negative definite planes in
$L_d\otimes\RR$.
Associated to this data, we have the Shimura variety ${\rm Sh}(L_d)$.
It is a smooth Deligne--Mumford stack over $\QQ$ such that, as complex orbifolds,
its $\CC$-valued points are given by the double quotient
$$
   {\rm Sh}(L_d)(\CC)\,=\, 
   G_{L_d}(\QQ)\backslash \left( Y_{L_d}\times G_d(\Aff_f) \right) / K_{L_d},
$$
and we refer to \cite[Section 4.1]{Madapusi Pera} for details.

Let us now return to K3 surfaces:
in Section \ref{sec: moduli spaces k3} we introduced
the moduli space ${\cal M}_{2d,\ZZ[\frac{1}{2d}]}^\circ$ of degree-$2d$ primitively
polarized K3 surfaces over $\ZZ[\frac{1}{2d}]$.
Rather than working with this moduli space, we will add {\em spin structures} first:
let ${\cal L}$ be a primitive polarization with ${\cal L}^2=2d$.
Let $P^2_\ell(X)$ be the primitive $\ell$-adic cohomology of $X$, that is,
the orthogonal complement of $c_1({\cal L})$ inside $\Het{2}(X,\ZZ_\ell)$.
For us, a spin structure is a choice of isometric isomorphism
$$
 \det(L_d)\otimes\ZZ_2 \,\stackrel{\iso}{\longrightarrow}\, \det( P^2_2(X) ), 
$$
and we denote by
$\widetilde{{\cal M}}_{2d,\ZZ[\frac{1}{2d}]}^\circ$ the moduli space of primitively
polarized K3 surfaces together with a choice of spin structure.
Forgetting the spin structure induces a morphism
$$
  \widetilde{{\cal M}}_{2d,\ZZ[\frac{1}{2d}]}^\circ\,\to\,{\cal M}^\circ_{2d,\ZZ[\frac{1}{2d}]},
$$
which is \'etale of degree $2$, and we
refer to \cite[Section 4.1]{Madapusi Pera} and
\cite[Section 6]{Rizov} for details and precise definitions.

Now, over the complex numbers, there classical period map can 
be interpreted as a morphism
$$
\imath_\CC\,:\,\widetilde{{\cal M}}_{2d,\CC}\,\to\,{\rm Sh}(L_d)_\CC.
$$
Obviously, the left side can be defined over $\QQ$ by considering families
of K3 surface with polarization and spin structures over $\QQ$.
The right side possesses a canonical model over $\QQ$.
And then, as shown by Rizov \cite{Rizov CM}, the period map $\imath_\CC$ 
descends to a map $\imath_\QQ$ over $\QQ$.

Thus, by what we have said above, the following result \cite[Theorem 5]{Madapusi Pera}
is a Torelli type theorem for K3 surfaces in positive and mixed
characteristic

\begin{Theorem}[Madapusi Pera]
 There exists a regular integral model ${\cal S}(L_d)$ for ${\rm Sh}(L_d)$
 over $\ZZ[\frac{1}{2}]$ such that $\imath_\QQ$ extends to an
 \'etale map
 $$
   \imath_{\ZZ[\frac{1}{2}]}\,:\, \widetilde{\cal M}^\circ_{2d,\ZZ[\frac{1}{2}]}\,\to\,
   {\cal S}(L_d)\,.
 $$
\end{Theorem}

When adding level structures, one can even achieve a period map that
is an open immersion \cite[Corollary 4.15]{Madapusi Pera}.
As explained in \cite[Section 1]{Madapusi Pera},  the construction of
this map over $\ZZ[\frac{1}{2d}]$ is essentially due to Rizov \cite{Rizov Kuga Satake}, and 
another construction is due to Vasiu \cite{Vasiu}.
Finally, let us also mention that Nygaard \cite{Nygaard Torelli} proved a Torelli-type
theorem for ordinary K3 surfaces using the theory of canonical lifts for such surfaces
and then applying the Kuga--Satake construction.

\end{document}